\documentclass{gtart_h}

\def\ifplaintex{\expandafter\ifx\csname documentclass\endcsname\relax}

\def\gtp{{\mathsurround=0pt\it $\cal G\mskip-2mu$eometry \&\ 
$\cal T\!\!$opology $\cal P\!$ublications}}  

\def\recd{{\small Received:\qua\receiveddate\ifx\reviseddate\relax
\else\qquad Revised:\qua\reviseddate\fi\par}} 

\def\lognumber#1{\def\thelognumber{#1}}
\def\volumenumber#1{\def\thevolumenumber{#1}}
\def\volumeyear#1{\def\thevolumeyear{#1}}
\def\papernumber#1{\def\thepapernumber{#1}}
\def\pagenumbers#1#2{\def\startpage{#1}\def\finishpage{#2}}
\def\published#1{\def\publishdate{#1}}

\def\received#1{\def\receiveddate{#1}}
\def\revised#1{\def\reviseddate{#1}}
\def\accepted#1{\def\accepteddate{#1}}

\def\asciiaddress#1{\def\theasciiaddress{#1}}
\def\asciiemail#1{\def\theasciiemail{#1}}

\long\def\asciiabstract#1{\long\def\theasciiabstract{#1}}
\def\asciikeywords#1{\def\theasciikeywords{#1}}

\let\\\par\let\thelognumber\relax\let\thevolumenumber\relax
\let\thepapernumber\relax\let\thevolumeyear\relax\let\startpage\relax
\let\finishpage\relax\let\publishdate\relax\let\receiveddate\relax
\let\reviseddate\relax\let\accepteddate\relax\let\theasciititle\relax
\let\theasciiauthors\relax\let\theasciiaddress\relax
\let\theasciiabstract\relax\let\theasciikeywords\relax

\let\theasciiemail\relax

\ifplaintex
\font\logobig=cmssbx10 scaled 3836
\font\logomed=cmssbx10 scaled 2557
\else
\font\logobig=cmssbx10 scaled 4200
\font\logomed=cmssbx10 scaled 2800
\fi

\long\def\makeagttitle{   
\count0=\startpage
\agt\hfill      
\hbox to 45truept{\vbox to 0pt{\vglue -13truept{\logomed A\kern -.37em{\logobig 
T}\kern -.38em G}\vss}\hss}
\break
{\small Volume \thevolumenumber\ (\thevolumeyear)
\startpage--\finishpage\nl
Published: \publishdate}

\vglue .25truein

{\parskip=0pt\leftskip 0pt plus
1fil\def\\{\par\smallskip}{\Large\bf\thetitle}\par\medskip} \vglue
0.05truein

{\parskip=0pt\leftskip 0pt plus 1fil\def\\{\par}{\sc\theauthors}
\par\medskip}%
 
\vglue 0.03truein 

{\small\leftskip 25truept\rightskip 25truept{\bf Abstract}\stdspace\theabstract

{\bf AMS Classification}\stdspace\theprimaryclass
\ifx\thesecondaryclass\relax\else; \thesecondaryclass\fi\par
{\bf Keywords}\stdspace \thekeywords\par}\vglue 7truept

}   

\ifplaintex
\font\phead=cmsl9 scaled 950
\font\pnum=cmbx10 scaled 913
\font\pfoot=cmsl9 scaled 950
\headline{\vbox to 0pt{\vskip -4.5mm\line{\small\phead\ifnum
\count0=\startpage ISSN 1472-2739 (on-line) 1472-2747 (printed)
\hfill {\pnum\folio}\else\ifodd\count0\def\\{ }%
\ifx\theshorttitle\relax\thetitle\else\theshorttitle\fi\hfill{\pnum\folio}
\else\def\\{ and }{\pnum\folio}\hfill\ifx\theshortauthors\relax\theauthors
\else\theshortauthors\fi\fi\fi}\vss}}
\footline{\vbox to 0pt{\vglue 0mm\line{\small\pfoot\ifnum\count0=\startpage
\copyright\ \gtp\hfill\else
\agt, Volume \thevolumenumber\ (\thevolumeyear)\hfill\fi}\vss}}
\else
\font=cmsl9 scaled 1050
\font\lnum=cmbx10 
\font=cmsl9 scaled 1050
\makeatletter
\def\@oddhead{{\small\lhead\ifnum\count0=\startpage ISSN 1472-2739 
(on-line) 1472-2747 (printed)\hfill {\lnum\number\count0}\else\ifodd\count0
\def\\{ }\ifx\theshorttitle\relax \thetitle \else\theshorttitle\fi\hfill
{\lnum\number\count0}\else\def\\{ and }{\lnum\number\count0}
\hfill\ifx\theshortauthors\relax 
\theauthors\else\theshortauthors\fi\fi\fi}}\def\@evenhead{\@oddhead}
\def\@oddfoot{\small\lfoot\ifnum\count0=\startpage\copyright\ \gtp\hfill\else
\agt, Volume \thevolumenumber\ (\thevolumeyear)\hfill\fi}
\def\@evenfoot{\@oddfoot}
\makeatother
\fi
\let\maketitlepage\makeagttitle
\let\makeshorttitle\maketitlepage
\let\maketitle\maketitlepage

\newwrite\gtoutfile
\long\gdef\makeheadfile{  
{\def\\{, }\def\s{ }
\immediate\openout\gtoutfile head.xxx
\immediate\write\gtoutfile{Proxy-for: \ifx\theasciiauthors\relax
\theauthors\else\theasciiauthors\fi\s<\ifx\theasciiemail\relax\theemail\else\theasciiemail\fi>}
\immediate\write\gtoutfile{\noexpand\\}
\immediate\write\gtoutfile{Authors: \ifx\theasciiauthors\relax
\theauthors\else\theasciiauthors\fi}
{\def\\{ }\immediate\write\gtoutfile{Title: \ifx\theasciititle\relax
\thetitle\else\theasciititle\fi}}
\immediate\write\gtoutfile{Subj-class: GT or SG, GR etc}
\immediate\write\gtoutfile{MSC-class: \theprimaryclass\ifx\thesecondaryclass\relax\else, \thesecondaryclass\fi}
\immediate\write\gtoutfile{Journal-ref: Algebr. Geom. Topol. \thevolumenumber\s
(\thevolumeyear) \startpage-\finishpage}
\immediate\write\gtoutfile{Comments: Published by Algebraic and
Geometric Topology at}
\immediate\write\gtoutfile{\s\s\s  http://www.maths.warwick.ac.uk/agt/AGTVol\thevolumenumber/agt-\thevolumenumber-\thepapernumber.abs.html}
\immediate\write\gtoutfile{\noexpand\\}
\immediate\write\gtoutfile{}
\ifx\theasciiabstract\relax
\immediate\write\gtoutfile{\theabstract}\else
\immediate\write\gtoutfile{\theasciiabstract}\fi
\immediate\write\gtoutfile{}
\immediate\write\gtoutfile{\noexpand\\}
\immediate\write\gtoutfile{}
\immediate\closeout\gtoutfile}}  

\def\maketitlepage{\makeagttitle\makeheadfile}
\let\makeshorttitle\maketitlepage
\let\maketitle\maketitlepage

\lognumber{49}
\volumenumber{5}
\volumeyear{2005}
\papernumber{49}
\pagenumbers{1223}{1290}
\received{16 July 2004} 
\revised{21 September 2005}
\accepted{29 September 2005}
\published{5 October 2005}

\usepackage{amssymb, amsmath, xy} 
\xyoption{all}

\newcommand\THH{{\it THH}}

\newcommand{\CalC}{\mathcal C}
\newcommand{\CalN}{\mathcal N}
\newcommand{\CalO}{\mathcal O}
\newcommand{\Ch}{\operatorname{Ch}}
\newcommand{\cok}{\operatorname{cok}}
\newcommand{\cotensor}{\mathbin\square}
\newcommand{\Ext}{\operatorname{Ext}}
\newcommand{\F}{\mathbb F}
\newcommand{\sgn}{\operatorname{sgn}}
\newcommand{\tmf}{tm\!f}
\newcommand{\Tor}{\operatorname{Tor}}
\newcommand{\Z}{\mathbb Z}
\renewcommand{\:}{\colon}

\theoremstyle{plain}
\newtheorem{thm}[equation]{Theorem}
\newtheorem{lem}[equation]{Lemma}
\newtheorem{prop}[equation]{Proposition}
\newtheorem{cor}[equation]{Corollary}

\theoremstyle{definition}
\newtheorem{defn}[equation]{Definition}
\newtheorem{rem}[equation]{Remark}

\newtheoremstyle{bold}{14pt plus6pt minus6pt}%
{6pt plus3pt minus3pt}{\bf}{}{\bf}{}{1em}%
{\thmname{#1}\thmnumber{#2}\thmnote{\bf\stdspace[#3]}}

\theoremstyle{bold}
\newtheorem{sshead}[equation]{}

\numberwithin{equation}{section}

\begin{document}

\title{Hopf algebra structure on\\topological Hochschild homology}
\authors{Vigleik Angeltveit\\John Rognes}
\address{Department of Mathematics, Massachusetts Institute of
	Technology\\Cambridge, MA 02139-4307, USA}
\secondaddress{Department of Mathematics, University of Oslo\\Blindern
NO-0316, Norway}

\asciiaddress{Department of Mathematics, Massachusetts Institute of
Technology\\Cambridge, MA 02139-4307, USA\\and\\Department of 
Mathematics, University of Oslo\\Blindern
NO-0316, Norway}
\asciiemail{vigleik@math.mit.edu, rognes@math.uio.no }
\gtemail{\mailto{vigleik@math.mit.edu}{\rm\qua 
and\qua}\mailto{rognes@math.uio.no}}

\begin{abstract}   
The topological Hochschild homology $\THH(R)$ of a commutative
$S$-algebra ($E_\infty$ ring spectrum) $R$ naturally has the structure
of a commutative $R$-algebra in the strict sense, and of a Hopf algebra
over~$R$ in the homotopy category.  We show, under a flatness
assumption, that this makes the B{\"o}kstedt spectral sequence
converging to the mod~$p$ homology of $\THH(R)$ into a Hopf algebra
spectral sequence.  We then apply this additional structure to the
study of some interesting examples, including the commutative
$S$-algebras $ku$, $ko$, $\tmf$, $ju$ and~$j$, and to calculate the
homotopy groups of $\THH(ku)$ and $\THH(ko)$ after smashing with suitable
finite complexes.  This is part of a program to make systematic
computations of the algebraic $K$-theory of $S$-algebras, by means of
the cyclotomic trace map to topological cyclic homology.
\end{abstract}

\asciiabstract{%
The topological Hochschild homology THH(R) of a commutative S-algebra
(E_infty ring spectrum) R naturally has the structure of a commutative
R-algebra in the strict sense, and of a Hopf algebra over R in the
homotopy category.  We show, under a flatness assumption, that this
makes the Boekstedt spectral sequence converging to the mod p homology
of THH(R) into a Hopf algebra spectral sequence.  We then apply this
additional structure to the study of some interesting examples,
including the commutative S-algebras ku, ko, tmf, ju and j, and to
calculate the homotopy groups of THH(ku) and THH(ko) after smashing
with suitable finite complexes.  This is part of a program to make
systematic computations of the algebraic K-theory of S-algebras, by
means of the cyclotomic trace map to topological cyclic homology.}




\asciikeywords{Topological Hochschild homology, commutative S-algebra,
coproduct, Hopf algebra, topological K-theory, image-of-J spectrum,
Boekstedt spectral sequence, Steenrod operations, Dyer-Lashof
operations.}

\maketitle 

%

\section{ Introduction }

The topological Hochschild homology $\THH(R)$ of an $S$-algebra $R$ (or
an $A_\infty$ ring spectrum, or a functor with smash product) was
constructed in the mid-1980's by B{\"o}kstedt \cite{Bo1}, as the
natural promotion of the classical Hochschild homology of an algebra in
the category of vector spaces (equipped with tensor product) to one in
the category of spectra (equipped with smash product).  It is the
initial ingredient in the construction by B{\"o}kstedt, Hsiang and
Madsen \cite{BHM93} of the topological cyclic homology $TC(R; p)$,
which in many cases closely approximates the algebraic $K$-theory
$K(R)$ \cite{Mc97}, \cite{Du97}, \cite{HM97}.

When $R = H\CalO_E$ is the Eilenberg--Mac\,Lane spectrum of the
valuation ring in a local number field, systematic computations of the
topological cyclic homology of $R$ were made in \cite{HM03}, thereby
verifying the Lichtenbaum--Quillen conjectures for the algebraic
$K$-theory of these fields.  Particular computations for other
commutative $S$-algebras, like connective complex $K$-theory~$ku$ and
its $p$-local Adams summand~$\ell$, have revealed a more general
pattern of how algebraic $K$-theory creates a ``red-shift'' in
chromatic filtration and satisfies a Galois descent property
\cite{AR02}, \cite{Au}.

These results indicate that the algebraic $K$-theory of a commutative
$S$-algebra is governed by an associated ``$S$-algebro-geometric''
structure space, for which the smashing localizations of the chromatic
filtration are related to a Zariski topology, and {\'e}tale covers and
Galois extensions are related to an {\'e}tale topology.  There are also
variant geometries associated to the less strictly commutative
$S$-algebras known as $E_n$ ring spectra, for $1 \le n \le \infty$.
The extreme cases $n=1$ and $n=\infty$ correspond to the associative
and the commutative $S$-algebras, respectively.  We shall encounter
examples of the intermediate $E_n$ ring spectrum structures in
Section~5, when we consider the Brown--Peterson spectrum $BP$ (which is
an $E_4$ ring spectrum according to Basterra and Mandell) and the
Johnson--Wilson spectra $BP\langle m{-}1\rangle$.  By convention, we
let $BP = BP\langle\infty\rangle$.

Further systematic computations of the topological Hochschild homology,
topological cyclic homology and algebraic $K$-theory of $S$-algebras
can be expected to shed new light on these geometries and on the
algebraic $K$-theory functor.  The present paper advances the algebraic
topological foundations for such systematic computations, especially by
taking into account the Hopf algebra structure present in the
topological Hochschild homology of commutative $S$-algebras.  This
program is continued in \cite{BR05}, which analyzes the differentials
in the homological homotopy fixed point spectral sequence that
approximates the cyclic fixed points of topological Hochschild
homology, and in \cite{L-N05}, which identifies the action by Steenrod
operations on the (continuous co-)homology of these fixed- and homotopy
fixed point spectra.

\medskip

When $R$ is a commutative $S$-algebra (or an $E_\infty$ ring spectrum,
or perhaps a commutative FSP), there is an equivalence $\THH(R) \simeq R
\otimes S^1$ of commutative $R$-algebras, due to McClure, Schw{\"a}nzl
and Vogt \cite{MSV97}, where $S^1$ is the topological circle and $(-)
\otimes S^1$ refers to the topologically tensored structure in the
category of commutative $S$-algebras.  The pinch map $S^1 \to S^1 \vee
S^1$ and reflection $S^1 \to S^1$ then induce maps $\psi \: \THH(R) \to
\THH(R) \wedge_R \THH(R)$ and $\chi \:  \THH(R) \to \THH(R)$, which make
$\THH(R)$ a Hopf algebra over~$R$ in the homotopy category
\cite[IX.3.4]{EKMM97}.  See also Theorem~\ref{t3.9}.

We wish to apply this added structure for computations.  Such
computations are usually made by starting with the simplicial model
$[q] \mapsto R^{\wedge(q+1)} = R \otimes S^1_q$ for $\THH(R)$, where now
$[q] \mapsto S^1_q$ is the simplicial circle.  The resulting skeleton
filtration on $\THH(R)$ gives rise to the B{\"o}kstedt spectral sequence
in homology
$$
E^2_{**}(R) = HH_*(H_*(R; \F_p))
\Longrightarrow H_*(\THH(R); \F_p) \,,
$$
as we explain in Section~4.  Compare \cite[\S2]{BM94} and
\cite[3.1]{MS93}.  But the coproduct and conjugation maps are not
induced by simplicial maps in this model, so some adjustment is needed
in order for these structures to carry over to the spectral sequence.
This we arrange in Section~3, by using a doubly subdivided simplicial
circle to provide an alternative simplicial model for $\THH(R)$, for
which the Hopf algebra structure maps can be simplicially defined.  The
verification of the Hopf algebra relations then also involves a triply
subdivided simplicial circle.

In Section~4 we transport the Hopf algebra structure on $\THH(R)$ to the
B{\"o}kstedt spectral sequence, showing in Theorem~\ref{t4.5} that if
its initial term $E^2_{**}(R)$ is flat as a module over $H_*(R; \F_p)$,
then this term is an $A_*$-comodule $H_*(R; \F_p)$-Hopf algebra and the
$d^2$-differentials respect this structure.  Furthermore, if every
$E^r$-term is flat over $H_*(R; \F_p)$, then the B{\"o}kstedt spectral
sequence is one of $A_*$-comodule $H_*(R; \F_p)$-Hopf algebras.  We
also discuss, in Proposition~\ref{p4.3}, the weaker algebraic structure
available in the B{\"o}kstedt spectral sequence when $R$ is just an
$E_2$- or $E_3$ ring spectrum.

Thereafter we turn to the computational applications.  The mod~$p$
homology of the topological Hochschild homology of the
Eilenberg--Mac\,Lane $S$-algebras $H\F_p$ and $H\Z$ was already
computed by B{\"o}kstedt \cite{Bo2}.  The first non-algebraic example,
namely the topological Hochschild homology of the Adams summand $\ell =
BP\langle1\rangle$ of $p$-local connective topological $K$-theory, was
computed for $p$ odd by McClure and Staffeldt in \cite{MS93}.

We show in Section~5 how to extend these computations to include the
case of $ku_{(2)} = BP\langle1\rangle$ at $p=2$, and the more general
Brown--Peterson and Johnson--Wilson $S$-algebras $BP\langle
m{-}1\rangle$, for primes $p$ and symbols $0 \le m \le \infty$ such
that these are $E_3$ ring spectra.  As a sample result we have part of
Theorem~\ref{t5.12}(a):
$$
H_*(\THH(BP); \F_2) \cong H_*(BP; \F_2) \otimes E(\sigma\bar\xi_k^2
\mid k\ge1) \,.
$$
Here and below, $E(-)$ and $P(-)$ denote the exterior and polynomial
algebras on the indicated generators, respectively, and the classes
$\bar\xi_k$ are the conjugates of Milnor's generators for the dual
Steenrod algebra $A_*$.

To resolve the multiplicative extensions in the B{\"o}kstedt spectral
sequence we provide a proof of B{\"o}kstedt's formula saying that the
suspension map $\sigma \:  \Sigma R \to \THH(R)$ takes the Dyer--Lashof
operations $Q^k$ on the homology of the commutative $S$-algebra $R$
compatibly to the corresponding operations on the homology of the
commutative $S$-algebra $\THH(R)$.  See Proposition~\ref{p5.9}, which
also makes precise what happens when $R$ is just an $E_{n+1}$ ring
spectrum.

In Section~6 we do the same for the higher real commutative
$S$-algebras $ko$ and $\tmf$ at $p=2$.  As sample results we have
Corollary~\ref{c5.14}(a) and Theorem~\ref{t6.2}(a):
$$
H_*(\THH(ku); \F_2) \cong H_*(ku; \F_2) \otimes
E(\sigma\bar\xi_1^2, \sigma\bar\xi_2^2) \otimes P(\sigma\bar\xi_3)
$$
and
$$
H_*(\THH(ko); \F_2) \cong H_*(ko; \F_2) \otimes E(\sigma\bar\xi_1^4,
\sigma\bar\xi_2^2) \otimes P(\sigma\bar\xi_3) \,.
$$

In the more demanding Section~7 we proceed to the $p$-local real and
complex image-of-$J$ spectra $j$ and $ju$, which are connective,
commutative $S$-algebras.  At odd primes the two are homotopy
equivalent.  We identify the mod~$p$ homology algebra of $ju$ at $p=2$
and of $j=ju$ at odd primes in Proposition~\ref{p7.12}(a) and~(b), and
make essential use of our results about Hopf algebra structures to show
that the corresponding B{\"o}kstedt spectral sequences for $\THH$
collapse at the $E^2$- and $E^p$-terms, respectively, in
Proposition~\ref{p7.13}(a) and~(b).  In Theorem~\ref{t7.15} we resolve
the algebra extension questions to obtain $H_*(\THH(ju); \F_p)$ as an
algebra, both for $p=2$ and for $p$ odd.  This proof involves a
delicate comparison with the case of $\THH(ku)$ for $p=2$, and with
$\THH(\ell)$ for $p$ odd.  Again as a sample result, we have
Theorem~\ref{t7.15}(a):
$$
H_*(\THH(ju); \F_2) \cong H_*(ju; \F_2) \otimes E(\sigma\bar\xi_1^4,
\sigma\bar\xi_2^2) \otimes P(\sigma\bar\xi_3) \otimes \Gamma(\sigma b)
\,.
$$
Here $\Gamma(\sigma b) = E(\gamma_{2^k}(\sigma b) \mid k\ge0)$ is the
divided power algebra on a class $\sigma b$ in degree~$4$.

The algebra structure of $H_*(j; \F_2)$ is described as a split
square-zero extension of $(A/\!/A_2)_*$ in Proposition~\ref{p7.12}(c):
$$
0 \to A_* \cotensor_{A_{2*}} \Sigma^7 K_* \to H_*(j; \F_2) \to
(A/\!/A_2)_* \to 0 \,.
$$
Here $A_2 = \langle Sq^1, Sq^2, Sq^4 \rangle \subset A$, and $K_*
\subset A_{2*}$ is dual to a cyclic $A_2$-module $K$ of rank~$17$ over
$\F_2$.  The symbol $\cotensor$ denotes the cotensor product of
comodules \cite[2.2]{MiMo65}.  The $A$-module structure of $H^*(j;
\F_2)$ was given in \cite{Da75}, but this identification of the algebra
structure seems to be new.  The $E^2$-term of the B{\"o}kstedt spectral
sequence for $j$ is described in Proposition~\ref{p7.13}(c), but it is
not flat over $H_*(j; \F_2)$, so the coproduct on $H_*(\THH(j); \F_2)$
is not conveniently described by this spectral sequence.  We have
therefore not managed to evaluate the homology of $\THH(j)$ at $p=2$ by
these methods.

Next we consider the passage from the homology of $\THH(R)$ to its
homotopy, with suitably chosen finite coefficients.  This has been a
necessary technical switch in past computations of topological cyclic
homology $TC(R; p)$, since $TC$ is defined as the homotopy inverse
limit of a diagram of fixed-point spectra derived from $\THH$, and the
interaction between inverse limits and homology was thought to be
difficult to control.  The homotopy groups of an inverse limit are much
better behaved.  Nonetheless, it may be that future computations of the
topological cyclic homology of $S$-algebras will follow a purely
homological approach, see \cite{BR05} and \cite{L-N05}.

In Section~8 we follow the strategy of \cite{MS93} to compute the
homotopy groups of $\THH(ku) \wedge M$, where $M = C_2$ is the mod~$2$
Moore spectrum, and of $\THH(ko) \wedge Y$, where $Y = C_2 \wedge
C_\eta$ is the $4$-cell spectrum employed by Mahowald \cite{Mah82}.  The
results appear in Theorems~\ref{t8.13} and~\ref{t8.14}, respectively.
In each case the method is to use the Adams spectral sequence to pass
from homology to homotopy, and to use a comparison with a Morava
$K(1)$-based B{\"o}kstedt spectral sequence to obtain enough
information about the $v_1$-periodic towers in the abutment to
completely determine the differential structure of the Adams spectral
sequence.

\medskip

The present paper started out as the first author's Master's thesis
\cite{An02} at the University of Oslo, supervised by the second
author.  Both authors are grateful to the referee for his careful
reading of the paper.

\section{ Hochschild and topological Hochschild homology }

Let $k$ be a graded field, i.e., a graded commutative ring such that
every graded $k$-module is free, and $\Lambda$ a graded $k$-algebra
(always unital and associative).  We recall the definition of the
Hochschild homology of $\Lambda$, e.g.~from \cite[X.4]{Mac75}.  The
Hochschild complex $C_*(\Lambda) = C^k_*(\Lambda)$ is the chain complex
of graded $k$-modules with $C_q(\Lambda) = \Lambda^{\otimes (q+1)}$ in
degree~$q$ (all tensor products are over~$k$) and boundary
homomorphisms $\partial \: C_q(\Lambda) \to C_{q-1}(\Lambda)$ given by
\begin{multline*}
\partial(\lambda_0 \otimes \dots \otimes \lambda_q) = \\
\sum_{i=0}^{q-1} (-1)^i \lambda_0 \otimes \dots \otimes \lambda_i
\lambda_{i+1} \otimes \dots \otimes \lambda_q + (-1)^{q+\epsilon}
\lambda_q \lambda_0 \otimes \dots \otimes \lambda_{q-1}
\end{multline*}
where $\epsilon = |\lambda_q| ( |\lambda_0| + \dots + |\lambda_{q-1}|
)$.  The Hochschild homology $HH_*(\Lambda) = HH^k_*(\Lambda)$ is
defined to be the homology of this chain complex.  It is bigraded,
first by the Hochschild degree~$q$ and second by the internal grading
from $\Lambda$.  When $\Lambda$ is commutative (always in the graded
sense) the shuffle product of chains defines a product
$$
\phi \: HH_*(\Lambda) \otimes_\Lambda HH_*(\Lambda)
\to HH_*(\Lambda)
$$
that makes $HH_*(\Lambda)$ a commutative $\Lambda$-algebra, with unit
corresponding to the inclusion of $0$-chains $\Lambda \to
HH_*(\Lambda)$.

When $\Lambda$ is commutative there is also a chain level coproduct
$\psi \: C_*(\Lambda) \to C_*(\Lambda) \otimes_\Lambda C_*(\Lambda)$
given in degree~$q$ by
\begin{multline}
\psi(\lambda_0 \otimes \lambda_1 \otimes \dots \otimes \lambda_q) = \\
\sum_{i=0}^q (\lambda_0 \otimes \lambda_1 \otimes \dots \otimes \lambda_i)
\otimes_\Lambda
(1 \otimes \lambda_{i+1} \otimes \dots \otimes \lambda_q) \,.
\label{e2.1}
\end{multline}
It is essential to tensor over~$\Lambda$ in the target of this chain map.
When $HH_*(\Lambda)$ is flat as a $\Lambda$-module, the chain level
coproduct $\psi$ induces a coproduct
$$
\psi \: HH_*(\Lambda) \to HH_*(\Lambda) \otimes_\Lambda HH_*(\Lambda)
$$
on Hochschild homology.  Here the right hand side is identified with
the homology of $C_*(\Lambda) \otimes_\Lambda C_*(\Lambda)$ by the
K{\"u}nneth theorem.  We shall now compare this chain level definition
of the coproduct on $HH_*(\Lambda)$ with an equivalent definition given
in more simplicial terms.

Let $B_*(\Lambda) = B_*(\Lambda, \Lambda, \Lambda)$ be the two-sided
bar construction \cite[X.2]{Mac75} for the $k$-algebra $\Lambda$.
It has $B_q(\Lambda) = \Lambda \otimes \Lambda^{\otimes q} \otimes
\Lambda$ in degree $q$, and is a free resolution of $\Lambda$ in the
category of $\Lambda$-bimodules.  We use the bar notation $\lambda_0
[\lambda_1 | \dots | \lambda_q] \lambda_{q+1}$ for a typical generator
of $B_q(\Lambda)$.  The Hochschild complex is obtained from the
two-sided bar construction by tensoring it with $\Lambda$ viewed as a
$\Lambda$-bimodule: $C_*(\Lambda) = \Lambda \otimes_{\Lambda-\Lambda}
B_*(\Lambda)$.

When $\Lambda$ is commutative, $B_*(\Lambda)$ is the chain complex
$\Ch(\Lambda \otimes \Delta^1)$ associated to the simplicial
$\Lambda$-bimodule $[q] \mapsto \Lambda \otimes \Delta^1_q$.
Here $\Delta^1$ is the simplicial $1$-simplex, and $\Lambda \otimes
\Delta^1_q$ denotes the tensor product of one copy of $\Lambda$ for each
element of $\Delta^1_q$.  See Section~3 below for more on this notation.
The $\Lambda$-bimodule structure on $B_*(\Lambda)$ is derived from the
inclusion of the two boundary points $\partial \Delta^1 \to \Delta^1$,
and $C_*(\Lambda)$ equals the chain complex $\Ch(\Lambda \otimes S^1)$
associated to the simplicial $\Lambda$-module $[q] \mapsto \Lambda \otimes
S^1_q$, where $S^1 = \Delta^1/\partial \Delta^1$ is the simplicial circle.

We now discuss three maps of $\Lambda$-bimodule chain complexes.
First, there is a canonical chain level coproduct $\psi \: B_*(\Lambda)
\to B_*(\Lambda) \otimes_\Lambda B_*(\Lambda)$ of $\Lambda$-bimodules,
given in degree~$q$ by
$$
\psi(\lambda_0 [\lambda_1 | \dots | \lambda_q] \lambda_{q+1})
= \sum_{i=0}^q \lambda_0 [\lambda_1 | \dots | \lambda_i] 1
\otimes_\Lambda 1 [\lambda_{i+1} | \dots | \lambda_q] \lambda_{q+1}
\,.
$$
When $\Lambda$ is commutative the chain level coproduct $\psi$ on
$C_*(\Lambda)$ is derived from this, as the obvious composite map
$$
\xymatrix{
C_*(\Lambda) = \Lambda \otimes_{\Lambda-\Lambda} B_*(\Lambda)
\ar[r]^-{1\otimes\psi} & \Lambda \otimes_{\Lambda-\Lambda} (B_*(\Lambda)
\otimes_\Lambda B_*(\Lambda)) \ar[d]^-{\psi'} \\
& C_*(\Lambda) \otimes_\Lambda C_*(\Lambda) \,.
}
$$

Second, there is a shuffle equivalence $sh \: B_*(\Lambda) \otimes_\Lambda
B_*(\Lambda) \to dB_*(\Lambda)$ of $\Lambda$-bimodules, by the
Eilenberg--Zilber theorem \cite[VIII.8.8]{Mac75} applied to two copies
of the simplicial $\Lambda$-bimodule $\Lambda \otimes \Delta^1$.
Here
$$
dB_*(\Lambda) = \Ch(\Lambda \otimes d\Delta^1) = \Ch((\Lambda \otimes
\Delta^1) \otimes_\Lambda (\Lambda \otimes \Delta^1))
$$
is the chain complex associated to the simplicial tensor product of
two copies of $\Lambda \otimes \Delta^1$, considered as simplicial
$\Lambda$-modules by way of the right and left actions, respectively.
This simplicial tensor product equals $\Lambda \otimes d\Delta^1$,
where the ``double $1$-simplex'' $d\Delta^1 = \Delta^1 \cup_{\Delta^0}
\Delta^1$ is the union of two $1$-simplices that are compatibly oriented.
(So $d\Delta^1$ is the $2$-fold edgewise subdivision of $\Delta^1$
\cite[\S1]{BHM93}.)  More explicitly,
$$
sh(x \otimes_\Lambda y) = \sum_{(\mu,\nu)} \sgn(\mu, \nu)
(s_\nu(x) \otimes_\Lambda s_\mu(y)) \,,
$$
where $x \in B_i(\Lambda)$, $y \in B_{q-i}(\Lambda)$, the sum is taken
over all $(i, q-i)$-shuffles $(\mu, \nu)$, $\sgn(\mu, \nu)$ is the sign
of the associated permutation and $s_\nu(x)$ and $s_\mu(y)$ are the
appropriate iterated degeneracy operations on $x$ and $y$,
respectively.

Third, there is a chain equivalence $\pi \: dB_*(\Lambda) \to
B_*(\Lambda)$ of $\Lambda$-bimodules, induced by the simplicial map
$\pi \: d\Delta^1 = \Delta^1 \cup_{\Delta^0} \Delta^1 \to \Delta^1$
that collapses the second $\Delta^1$ in $d\Delta^1$ to a point.
It is given by
$$
\pi(x \otimes_\Lambda y) = x \cdot \epsilon(y)
$$
for $x, y \in B_q(\Lambda)$, where $\epsilon(\lambda_0[\lambda_1|\dots
|\lambda_q]\lambda_{q+1}) = \lambda_0\lambda_1 \dots
\lambda_q\lambda_{q+1}$ is the augmentation.

\begin{lem}
\label{l2.2}
Let $\Lambda$ be a commutative $k$-algebra.  The maps $sh \circ \psi \:
B_*(\Lambda) \to dB_*(\Lambda)$ and $\pi \:  dB_*(\Lambda) \to
B_*(\Lambda)$ of $\Lambda$-bimodule chain complexes are mutual chain
inverses.  Hence the induced composite
$$
\xymatrix@C+11pt{
dC_*(\Lambda) = \Lambda \otimes_{\Lambda-\Lambda}
	dB_*(\Lambda) \ar[r]^-{1\otimes\pi} &
C_*(\Lambda) \ar[r]^-{1\otimes(sh \circ \psi)} & dC_*(\Lambda)
}
$$
is chain homotopic to the identity.
\end{lem}

\begin{proof}
All three chain complexes $B_*(\Lambda)$, $B_*(\Lambda) \otimes_\Lambda
B_*(\Lambda)$ and $dB_*(\Lambda)$ are free $\Lambda$-bimodule resolutions
of $\Lambda$, and the maps $\psi$, $sh$ and $\pi$ are $\Lambda$-bimodule
chain maps, so it suffices to verify that the composite $\pi \circ sh
\circ \psi \: B_*(\Lambda) \to B_*(\Lambda)$ covers the identity on
$\Lambda$.  In degree zero, $\psi(\lambda_0[]\lambda_1) = \lambda_0[]1
\otimes_\Lambda 1[]\lambda_1$, $sh(\lambda_0[]1 \otimes_\Lambda
1[]\lambda_1) = \lambda_0[]1 \otimes_\Lambda 1[]\lambda_1$ and
$\pi(\lambda_0[]1 \otimes_\Lambda 1[]\lambda_1) = \lambda_0[]\lambda_1$,
as required.  If desired, the explicit formulas can be composed also
in higher degrees, to show that $(\pi \circ sh \circ \psi)(x) \equiv
x$ modulo simplicially degenerate terms, for all $x \in B_q(\Lambda)$.
In either case, we see that $\pi \circ sh \circ \psi$ is chain homotopic
to the identity.  The remaining conclusions follow by uniqueness of inverses.
\end{proof}

Note that $dC_*(\Lambda)$ defined above is the chain complex associated
to the simplicial $\Lambda$-module $\Lambda \otimes d'S^1$, where the
``double circle'' $d'S^1 = d\Delta^1/\partial d\Delta^1$ is the
quotient of the double $1$-simplex $d\Delta^1 = \Delta^1
\cup_{\Delta^0} \Delta^1$ by its two end-points $\partial d\Delta^1$.
The chain equivalence $1 \otimes \pi \: dC_*(\Lambda) = \Ch(\Lambda
\otimes d'S^1) \to \Ch(\Lambda \otimes S^1) = C_*(\Lambda)$ obtained by
tensoring down the $\Lambda$-bimodule chain equivalence induced from
the collapse map $\pi \: d\Delta^1 \to \Delta^1$, is then more directly
obtained from the collapse map $\pi \: d'S^1 \to S^1$ that collapses
the second of the two non-degenerate $1$-simplices in $d'S^1$ to a
point.

There is also a simplicial pinch map $\psi \: d'S^1 \to S^1 \vee S^1$
to the one-point union (wedge sum) of two circles, that collapses the
$0$-skeleton of $d'S^1$ to a point.  It induces a map
$$
\psi' \: dC_*(\Lambda) = \Ch(\Lambda \otimes d'S^1) \to \Ch(\Lambda
\otimes (S^1 \vee S^1)) \,.
$$
The target is the chain complex associated to the simplicial tensor
product of two copies of $\Lambda \otimes S^1$, considered as a
simplicial $\Lambda$-module.  It is therefore also the target of a
shuffle equivalence, namely $sh \: C_*(\Lambda) \otimes_\Lambda
C_*(\Lambda) \to \Ch(\Lambda \otimes (S^1 \vee S^1))$.

\begin{prop}
\label{p2.3}
Let $\Lambda$ be a commutative $k$-algebra.
The composite map
$$
\xymatrix{
dC_*(\Lambda) \ar[r]^-{1\otimes\pi}_-{\simeq}
	& C_*(\Lambda) \ar[r]^-{\psi}
	& C_*(\Lambda) \otimes_\Lambda C_*(\Lambda) \ar[r]^-{sh}_-{\simeq}
	& \Ch(\Lambda \otimes (S^1 \vee S^1))
}
$$
of the chain level coproduct $\psi$ in $C_*(\Lambda)$ (see
formula~(\ref{e2.1})), with the chain equivalence $1 \otimes \pi$
induced by the simplicial collapse map $\pi \: d'S^1 \to S^1$ and the
shuffle equivalence $sh$, is chain homotopic to the map
$$
\xymatrix{
dC_*(\Lambda) = \Ch(\Lambda \otimes d'S^1) \ar[r]^-{\psi'} & \Ch(\Lambda
\otimes (S^1 \vee S^1))
}
$$
induced by the simplicial pinch map $\psi \: d'S^1 \to S^1 \vee S^1$.

Hence, when $HH_*(\Lambda)$ is flat over~$\Lambda$, the coproduct $\psi$
on Hochschild homology agrees, via the identifications induced by $1
\otimes \pi$ and $sh$, with the map $\psi'$ induced by the simplicial
pinch map.
\end{prop}

\begin{proof}
Consider the following diagram.
$$
\xymatrix{
C_*(\Lambda) \ar[r]^-{1\otimes\psi} & \Lambda \otimes_{\Lambda-\Lambda}
(B_*(\Lambda) \otimes_\Lambda B_*(\Lambda)) \ar[r]^-{\psi'}
\ar[d]^-{1\otimes sh} &
C_*(\Lambda) \otimes_\Lambda C_*(\Lambda) \ar[d]^-{sh} \\
& dC_*(\Lambda) \ar[r]^-{\psi'} \ar[ul]^{1\otimes\pi}
& {}\Ch(\Lambda \otimes (S^1 \vee S^1))
}
$$
The composite along the upper row is the coproduct $\psi$, the
composite around the triangle is chain homotopic to the identity by
Lemma~\ref{l2.2}, and the square commutes by naturality of the shuffle
map with respect to the pinch map $\psi'$.  A diagram chase then
provides the claimed chain homotopy.
\end{proof}

We shall make use of the following standard calculations of Hochschild
homology.  The formulas for the coproduct $\psi$ follow directly from
the chain level formula~(\ref{e2.1}) above.  Let $P(x) = k[x]$ and $E(x) =
k[x]/(x^2)$ be the polynomial and exterior algebras over~$k$ in one
variable $x$, and let $\Gamma(x) = k\{ \gamma_i(x) \mid i\ge0 \}$
be the {\it divided power algebra\/} with multiplication
$$
\gamma_i(x) \cdot \gamma_j(x) = (i, j) \ \gamma_{i+j}(x) \,,
$$
where $(i, j) = (i+j)!/i!j!$ is the binomial coefficient.

\begin{prop}
\label{p2.4}
For $x \in \Lambda$ let $\sigma x \in HH_1(\Lambda)$ be the homology class
of the cycle $1 \otimes x \in C_1(\Lambda)$ in the Hochschild
complex.  For $\Lambda = P(x)$ there is a $P(x)$-algebra isomorphism
$$
HH_*(P(x)) = P(x) \otimes E(\sigma x) \,.
$$
The class $\sigma x$ is coalgebra primitive, i.e., $\psi(\sigma x) =
\sigma x \otimes 1 + 1 \otimes \sigma x$.
For $\Lambda = E(x)$ there is an $E(x)$-algebra isomorphism
$$
HH_*(E(x)) = E(x) \otimes \Gamma(\sigma x) \,.
$$
The $i$th divided power $\gamma_i(\sigma x)$ is the homology class
of the cycle $1 \otimes x \otimes \dots \otimes x \in C_i(\Lambda)$.
The coproduct is given by
$$
\psi(\gamma_k(\sigma x)) = \sum_{i+j=k} \gamma_i(\sigma x) \otimes
\gamma_j(\sigma x) \,.
$$
There is a K{\"u}nneth formula
$$
HH_*(\Lambda_1 \otimes \Lambda_2) \cong HH_*(\Lambda_1) \otimes
HH_*(\Lambda_2)
\,.
$$
\end{prop}

Let $P_h(x) = k[x]/(x^h)$ be the truncated polynomial algebra of height
$h$.  We write $P(x_i \mid i\ge0) = P(x_0, x_1, \dots) = P(x_0) \otimes
P(x_1) \otimes \dots$, and so on.  When $k$ is of prime characteristic
$p$ it is a standard calculation with binomial coefficients that
\begin{equation}
\Gamma(x) = P_p(\gamma_{p^i}(x) \mid i\ge0)
\label{e2.5}
\end{equation}
as a $k$-algebra.

We have already noted that the Hochschild complex $C_*(\Lambda)$ is the
chain complex associated to a simplicial graded $k$-module $[q] \mapsto
C_q(\Lambda) = \Lambda^{\otimes(q+1)}$, with face maps $d_i$
corresponding to the individual terms in the alternating sum defining
the Hochschild boundary~$\partial$.  In fact, this is a {\it cyclic\/}
graded $k$-module in the sense of Connes, with cyclic structure maps
$t_q$ that cyclically permute the $(q+1)$ tensor factors in
$C_q(\Lambda)$, up to sign.  It follows that the geometric realization
$HH(\Lambda) = |[q] \mapsto C_q(\Lambda)|$ admits a natural
$S^1$-action $\alpha \:  HH(\Lambda) \wedge S^1_+ \to HH(\Lambda)$, and
that the Hochschild homology groups are the homotopy groups of this
space: $HH_*(\Lambda) = \pi_* HH(\Lambda)$.

\medskip

The basic idea in the definition of topological Hochschild homology is
to replace the ground ring $k$ by the sphere spectrum $S$, and the
symmetric monoidal category of graded $k$-modules under the tensor
product $\otimes = \otimes_k$ by the symmetric monoidal category of
spectra, interpreted as $S$-modules, under the smash product $\wedge =
\wedge_S$.  A monoid in the first category is a graded $k$-algebra
$\Lambda$, which then gets replaced by a monoid in the second category,
i.e., an $S$-algebra $R$.  To make sense of this we will work in the
framework of \cite{EKMM97}, but we could also use \cite{HSS00} or any
other reasonable setting that gives a symmetric monoidal category of
spectra.

The original definition of topological Hochschild homology was given by
B{\"o}k\-stedt in the mid 1980's \cite{Bo1}, inspired by work and conjectures
of Goodwillie and Waldhausen.  The following definition is not the one
originally used by B{\"o}kstedt, since he did not have the symmetric
monoidal smash product from \cite{EKMM97} or \cite{HSS00} available, but
it agrees with the heuristic definition that his more complicated definition
managed to make sense of, with the more elementary technology that he
had at hand.

\begin{defn}
\label{d2.6}
Let $R$ be an $S$-algebra, with multiplication $\mu \: R \wedge R \to
R$ and unit $\eta \: S \to R$.  The topological Hochschild homology of
$R$ is the geometric realization $\THH(R)$ of the simplicial $S$-module
$\THH_\bullet(R)$ with
$$
\THH_q(R) = R^{\wedge(q+1)}
$$
in simplicial degree~$q$.  The simplicial structure is like that on the
simplicial $k$-module underlying the Hochschild complex.  More
precisely, the $i$-th face map $d_i \: R^{\wedge(q+1)} \to R^{\wedge
q}$ equals $id_R^i \wedge \mu \wedge id_R^{q-i-1}$ for $0 \le i < q$,
while $d_q = (\mu \wedge id_R^{q-1}) t_q$, where $t_q$ cyclically
permutes the $(q+1)$ smash factors $R$ by moving the last factor to the
front.  The $j$-th degeneracy map $s_j \: R^{\wedge(q+1)} \to
R^{\wedge(q+2)}$ equals $id_R^{j+1} \wedge \eta \wedge id_R^{q-j}$, for
$0 \le j \le q$.

Furthermore, the cyclic operators $t_q \: \THH_q(R) \to \THH_q(R)$ make
$\THH_\bullet(R)$ a cyclic $S$-module, so that its geometric realization
has a natural $S^1$-action
$$
\alpha \: \THH(R) \wedge S^1_+ \to \THH(R) \,.
$$
The topological Hochschild homology groups of $R$ are defined
to be the homotopy groups $\pi_* \THH(R)$.
\end{defn}

When $R$ is a commutative $S$-algebra there is a product on $\THH(R)$
that makes it an augmented commutative $R$-algebra, with unit
corresponding to the inclusion of $0$-simplices $R \to \THH(R)$.
See~(\ref{e3.3}) below.  There is a discussion of alternative definitions
of $\THH(R)$ in \cite[IX]{EKMM97}.

\section{ Hopf algebra structure on $\THH$ }

Already B{\"o}kstedt noted that the simplicial structure on $\THH(R)$,
as defined above, is derived from the simplicial structure on the
standard simplicial circle $S^1 = \Delta^1/\partial \Delta^1$.  This
can be made most precise in the case when $R$ is a commutative
$S$-algebra, in which case there is a formula $\THH(R) \cong R \otimes
S^1$ in terms of the simplicial tensor structure on the category of
commutative $S$-algebras.  A corresponding formula $\THH(R) \simeq R
\otimes |S^1|$ in terms of the topological tensor structure was
discussed by McClure, Schw{\"a}nzl and Vogt in \cite{MSV97}.  We shall
stick to the simplicial context, since we will make use of the
resulting skeletal filtrations to form spectral sequences.

We now make this ``tensored structure'' explicit.  Let $R$ be a
commutative $S$-algebra and $X$ a finite set, and let
$$
R \otimes X = \bigwedge_{x \in X} R
$$
be the smash product of one copy of $R$ for each element of $X$.  It is
again a commutative $S$-algebra.  Now let $f \: X \to Y$ be a function
between finite sets, and let $R \otimes f \: R \otimes X \to R \otimes Y$
be the smash product over all $y \in Y$ of the maps
$$
R \otimes f^{-1}(y) = \bigwedge_{x \in f^{-1}(y)} R \to R = R \otimes \{y\}
$$
that are given by the iterated multiplication from the $\# f^{-1}(y)$
copies of $R$ on the left to the single copy of $R$ on the right.
If $f^{-1}(y)$ is empty, this is by definition the unit map $\eta \:
S \to R$.  Since $R$ is commutative, there is no ambiguity in how
these iterated multiplications are to be formed.

Note that the construction $R \otimes X$ is functorial in both $R$
and $X$.  Given an injection $X \to Y$ and any function $X \to Z$, there
is a natural isomorphism
$$
R \otimes (Y \cup_X Z)
\cong 
(R \otimes Y) \wedge_{(R \otimes X)} (R \otimes Z) \,.
$$
There is also a natural map
$$
R \wedge X_+ = \bigvee_{x \in X} R \to \bigwedge_{x \in X} R = R \otimes X
$$
whose restriction to the wedge summand indexed by $x \in X$ is $R \otimes
i_x$, where $i_x \: \{x\} \to X$ is the inclusion.

By naturality, these constructions all extend degreewise to simplicial
finite sets $X \: [q] \mapsto X_q$, simplicial maps $f \: X \to Y$, etc.
In particular, we can define the simplicial commutative $S$-algebra
$$
R \otimes X = \Bigl( [q] \mapsto R \otimes X_q \Bigr)
$$
with structure maps $R \otimes d_i$, $R \otimes s_j$, etc.
There is then a useful natural map
\begin{equation}
\omega \: R \wedge X_+ \to R \otimes X \,.
\label{e3.1}
\end{equation}

Consider the special case $X = S^1 = \Delta^1/\partial \Delta^1$.  Here
$\Delta^1_q$ has $(q+2)$ elements $\{x_0, \dots, x_{q+1}\}$ where $x_t
\: [q] \to [1]$ has $\# x_t^{-1}(0) = t$.  The quotient $S^1_q$ has
$(q+1)$ elements, obtained by identifying $x_0 \sim x_{q+1}$.  Then
$d_i(x_t) = x_t$ for $t \le i$ and $ d_i(x_t) = x_{t-1}$ for $t > i$,
while $s_j(x_t) = x_t$ for $t \le j$ and $s_j(x_t) = x_{t+1}$ for $t >
j$.  A direct check shows that there is a natural isomorphism
\begin{equation}
\THH(R) \cong R \otimes S^1
\label{e3.2}
\end{equation}
of (simplicial) commutative $S$-algebras.  In degree~$q$ it is the
obvious identification $R^{\wedge(q+1)} \cong R \otimes S^1_q$.

There are natural maps $\eta \: * \to S^1$, $\epsilon \: S^1 \to *$ and
$\phi \: S^1 \vee S^1 \to S^1$ that map to the base point of $S^1$,
retract to $*$ and fold two copies of $S^1$ to one, respectively.  By
naturality, these induce the following maps of commutative
$S$-algebras:
\begin{equation}
\begin{aligned}
\eta &\: R \to \THH(R) \\
\epsilon &\: \THH(R) \to R \\
\phi &\: \THH(R) \wedge_R \THH(R) \to \THH(R) \,.
\end{aligned}
\label{e3.3}
\end{equation}
In the last case, the product map involves the identification $$\THH(R)
\wedge_R \THH(R) \cong R \otimes (S^1 \cup_* S^1),$$ where $S^1 \cup_*
S^1 = S^1 \vee S^1$.  Taken together, these maps naturally make
$\THH(R)$ an augmented commutative $R$-algebra.

There is also a natural map
\begin{equation}
\omega \: R \wedge S^1_+ \to \THH(R) \,,
\label{e3.4}
\end{equation}
derived from~(\ref{e3.1}), which captures part of the circle action upon
$\THH(R)$.  More precisely, the map $\omega$ admits the following
factorization:
\begin{equation}
\omega = \alpha \circ (\eta \wedge id) \: 
R \wedge S^1_+ \to \THH(R) \wedge S^1_+ \to \THH(R) \,.
\label{e3.5}
\end{equation}
This is clear by inspection of the definition of the circle action
$\alpha$ on the $0$-simplices of $\THH(R)$.

\medskip

We would like to have a coproduct on $\THH(R)$, coming from a pinch map
$S^1 \to S^1 \vee S^1$, but there is no such simplicial map with our
basic model for $S^1$.  To fix this we again consider a ``double model''
for $S^1$, denoted $dS^1$, with
$$
dS^1 = (\Delta^1 \sqcup \Delta^1) \cup_{(\partial \Delta^1 \sqcup
\partial \Delta^1)} \partial \Delta^1 \,.
$$
Here $\partial \Delta^1 \sqcup \partial \Delta^1 \to \partial \Delta^1$ is
the identity map on each summand, so the two non-degenerate $1$-simplices
of $dS^1$ have opposing orientations in the geometrically realized circle.
It is the quotient of the barycentric subdivision of $\Delta^1$ by
its boundary.

Then we have a simplicial pinch map $\psi \: dS^1 \to S^1 \vee S^1$
that collapses $\partial \Delta^1 \subset dS^1$ to $*$, as well as a
simplicial flip map $\chi \: dS^1 \to dS^1$ that interchanges the two
copies of $\Delta^1$.

\begin{rem}
\label{r3.6}
The simplicial set $dS^1$ introduced here differs from the double
circle $d'S^1$ considered in Section~2, in that the orientation of the
second $1$-simplex has been reversed.  The switch is necessary here to
make the flip map $\chi$ simplicial.  In principle, we could have used
the same $dS^1$ in Section~2 as here, but this would have entailed the
cost of discussing the anti-simplicial involution
$\lambda_0[\lambda_1|\dots|\lambda_q]\lambda_{q+1} \mapsto \pm
\lambda_{q+1}[\lambda_q|\dots|\lambda_1]\lambda_0$ of $B_*(\Lambda)$,
and complicating the formula~(\ref{e2.1}) for the chain level coproduct
$\psi$.  We choose instead to suppress this point.
\end{rem}

We define a corresponding ``double model'' for $\THH(R)$, denoted
$d\THH(R)$, by
$$
d\THH(R) = R \otimes dS^1 \,.
$$
The pinch and flip maps now induce the
following natural maps of commutative $S$-algebras:
\begin{align}
\begin{aligned}
\psi' &\: d\THH(R) \to \THH(R) \wedge_R \THH(R) \\
\chi' &\: d\THH(R) \to d\THH(R) \,.
\end{aligned}
\label{e3.7}
\end{align}

\begin{lem}
\label{l3.8}
Let $R$ be cofibrant as an $S$-module.  Then the collapse map
$\pi \: dS^1 \to S^1$ that takes the second $\Delta^1$ to $*$
induces a weak equivalence
$$
\xymatrix{
\pi \: d\THH(R) \ar[r]^-{\simeq} & \THH(R) \,.
}
$$
\end{lem}

\begin{proof}
Consider the commutative diagram
$$
\xymatrix{
B(R) \ar@{=}[d] & R \wedge R \ar[l] \ar@{=}[d] \ar[r] & B(R) \ar[d] \\
B(R) & R \wedge R \ar[l] \ar[r] & R
}
$$
of commutative $S$-algebras.  Here $B(R) = B(R, R, R) = R \otimes
\Delta^1$ is the two-sided bar construction, its augmentation $B(R)
\to R$ is a weak equivalence, and the inclusion $R \wedge R \to B(R)$
is a cofibration of $S$-modules.  From \cite[III.3.8]{EKMM97} we know
that the categorical pushout (balanced smash product) in this
case preserves weak equivalences.  Pushout along the upper row gives
$d\THH(R)$ and pushout along the lower row gives $\THH(R)$, so
the induced map $\pi \: d\THH(R) \to \THH(R)$ is indeed a weak equivalence.
\end{proof}

\begin{thm}
\label{t3.9}
Let $R$ be a commutative $S$-algebra.  Its topological Hochschild
homology $\THH(R)$ is naturally an augmented commutative $R$-algebra,
with unit, counit and product maps $\eta$, $\epsilon$ and $\phi$
defined as in~(\ref{e3.3}) above.  In the stable homotopy category,
these maps, the coproduct map
$$
\psi = \psi' \circ \pi^{-1} \: \THH(R) \to \THH(R) \wedge_R \THH(R)
$$
and the conjugation map
$$
\chi = \pi \circ \chi' \circ \pi^{-1} \: \THH(R) \to \THH(R)
$$
naturally make $\THH(R)$ an $R$-Hopf algebra.
\end{thm}

\begin{proof}
To check that $\THH(R)$ is indeed a Hopf algebra over~$R$ in
the stable homotopy category we must verify that a number of
diagrams commute.  We will do one case that illustrates the
technique, and leave the rest to the reader.

Let $T = \THH(R)$.  In order to show that the diagram
\begin{equation}
\xymatrix{
T \ar[r]^-{\psi} \ar[dr]^{\epsilon} \ar[d]_-{\psi} &
T \wedge_R T \ar[dr]^{id \wedge \chi} \\
T \wedge_R T \ar[dr]_{\chi \wedge id} &
R \ar[dr]^{\eta} &
T \wedge_R T \ar[d]^-{\phi} \\
& T \wedge_R T \ar[r]_-{\phi} & T
}
\label{e3.10}
\end{equation}
commutes in the stable homotopy category it suffices to check
that the diagram of simplicial sets
$$
\xymatrix{
tS^1 \ar[r]^-{\psi_1} \ar[dr]^{\epsilon} \ar[d]_-{\psi_2} &
S^1 \vee dS^1 \ar[dr]^{id \vee \chi'} \\
dS^1 \vee S^1 \ar[dr]_{\chi' \vee id} &
{*} \ar[dr]^{\eta} &
S^1 \vee dS^1 \ar[d]^-{\phi(id \vee \pi)} \\
& dS^1 \vee S^1 \ar[r]_-{\phi(\pi \vee id)} & S^1
}
$$
homotopy commutes (simplicially).

Here $tS^1 = \partial \Delta^2$ is a ``triple model'' for $S^1$, with
three non-degenerate $1$-simplices.  The pinch map $\psi_1$ identifies
the vertices $0$ and $1$ in $\partial \Delta^2$ and takes the face
$\delta_0$ to the first $\Delta^1$ in $dS^1$.  Then the composite
$\phi(id \vee \pi)(id \vee \chi')\psi_1$ factors as
$$
\xymatrix{
\partial \Delta^2 \subset \Delta^2 \ar[r]^-{s_1} & \Delta^1 \to S^1 \,,
}
$$
and $\Delta^1$ is simplicially contractible.
Similarly, $\psi_2$ identifies the vertices $1$ and $2$ and
takes $\delta_2$ to the first $\Delta^1$ in $dS^1$.
Then $\phi(\pi \vee id)(\chi' \vee id)\psi_2$ factors as
$$
\xymatrix{
\partial \Delta^2 \subset \Delta^2 \ar[r]^-{s_0} & \Delta^1 \to S^1 \,,
}
$$
and again this map is simplicially contractible.

Finally we use a weak equivalence $t\THH(R) = R \otimes tS^1 \to
\THH(R)$, as in Lemma~\ref{l3.8}, to deduce that the diagram~(\ref{e3.10})
indeed homotopy commutes.
\end{proof}

\begin{rem}
\label{r3.11}
As noted above, $\THH(R)$ is commutative as an $R$-algebra.  The pinch
map $\psi \: dS^1 \to S^1 \vee S^1$ is not homotopy cocommutative, so
we do not expect that the coproduct $\psi \: \THH(R) \to \THH(R) \wedge_R
\THH(R)$ will be cocommutative in any great generality.  The Hopf
algebras arising in this fashion will therefore be commutative, but not
cocommutative.
\end{rem}

The inclusion of the base point $\eta \: * \to S^1$ induces a cofiber
sequence of $R$-modules
$$
\xymatrix{
R = R \wedge *_+ \ar[r]^-{1\wedge\eta_+} & R \wedge S^1_+ \ar[r]^-{j}
	& R \wedge S^1 = \Sigma R
}
$$
which is canonically split by the retraction map
$$
\xymatrix{
R \wedge S^1_+ \ar[r]^-{1\wedge\epsilon_+} & R \wedge *_+ = R \,.
}
$$
Hence in the stable homotopy category there is a canonical section $\kappa
\: \Sigma R \to R \wedge S^1_+$ to the map labeled $j$ above.  We let
\begin{equation}
\sigma = \omega \circ \kappa \: \Sigma R \to R \wedge S^1_+ \to \THH(R)
\label{e3.12}
\end{equation}
be the composite stable map.  It induces an operator
$$
\sigma \: H_*(\Sigma R; \F_p) \to H_*(\THH(R); \F_p) \,,
$$
which we in Proposition~\ref{p4.9} shall see is compatible with that of
Proposition~\ref{p2.4}.

\section{ The B{\"o}kstedt spectral sequence }

Let $R$ be an $S$-algebra.  To calculate the mod~$p$ homology $H_*(\THH(R);
\F_p)$ of its topological Hochschild homology, B{\"o}kstedt constructed
a strongly convergent spectral sequence
\begin{equation}
E^2_{s,*}(R) = HH_s(H_*(R; \F_p))
\Longrightarrow H_{s+*}(\THH(R); \F_p) \,,
\label{e4.1}
\end{equation}
using the skeleton filtration on $\THH(R)$.  In fact,
$$
E^1_{s,*}(R) = H_*(R^{\wedge(s+1)}; \F_p) \cong H_*(R;
\F_p)^{\otimes(s+1)} = C_s(H_*(R; \F_p))
$$
equals the Hochschild $s$-chains of the algebra $\Lambda = H_*(R; \F_p)$
over $k = \F_p$, and the $d^1$-differential can as usual be identified
with the Hochschild boundary operator~$\partial$.  To be quite precise, the
$E^1$-term is really the associated normalized complex $\Lambda \otimes
\bar\Lambda^{\otimes s}$, with $\bar\Lambda = \Lambda/k$, but this change
does not affect the $E^2$-term.

This is naturally a spectral sequence of $A_*$-comodules, where $A_* =
H_*(H\F_p; \F_p)$ is the dual of the mod~$p$ Steenrod algebra, since it
is obtained by applying mod~$p$ homology to a filtered spectrum.  If
$R$ is a commutative $S$-algebra, the spectral sequence admits more
structure.

\begin{prop}
\label{p4.2}
Let $R$ be a commutative $S$-algebra.  Then $H_*(\THH(R); \F_p)$ is an
augmented commutative $A_*$-comodule $H_*(R; \F_p)$-algebra, and the
B{\"o}k\-stedt spectral sequence $E^r_{**}(R)$ is an augmented
commutative $A_*$-comodule $H_*(R; \F_p)$-algebra spectral sequence,
converging to $H_*(\THH(R); \F_p)$.
\end{prop}

\begin{proof}
We know that $\THH(R)$ is an augmented commutative $R$-algebra by
Theorem~\ref{t3.9}, and the relevant structure maps~(\ref{e3.3}) are
all maps of simplicial spectra.  Hence they respect the skeleton
filtration on $\THH(R)$, and we have in particular a composite map of
spectral sequences
$$
\xymatrix{
E^r_{**}(R) \otimes_\Lambda E^r_{**}(R) \to {}'E^r_{**} \ar[r]^-{\phi} &
E^r_{**}(R)
}
$$
with ${}'E^r_{**}$ the spectral sequence associated to the skeleton
filtration on the smash product $\THH(R) \wedge_R \THH(R)$.  The left
hand map is induced by the usual homology cross product from the
$E^1$-term and onwards.  This defines the algebra structure on
$E^r_{**}(R)$, and the remaining claims are straightforward.
\end{proof}

In fact, we can settle for less than strict commutativity to have an
algebra structure in the B{\"o}kstedt spectral sequence.  We shall make
use of this for some examples related to the Brown--Peterson spectrum
in Section~5.  An {\it $E_n$ ring spectrum} \cite[I.4]{BMMS86} is an
$S$-module with an action by the little $n$-cubes operad $\CalC_n$, or
one that is related to it by a chain of $\Sigma$-equivariant weak
equivalences.  Fiedorowicz and Vogt \cite[3.4]{FV} show that any
$E_{n+1}$ ring spectrum is weakly equivalent to an associative
$S$-algebra with an action by the little $n$-cubes operads, i.e., to a
$\CalC_n$-algebra in the category of associative $S$-algebras.
Basterra and Mandell (unpublished) have obtained the same result by a
different method.

Thus any continuous functor from associative $S$-algebras takes an
$E_{n+1}$ ring spectrum to a $\CalC_n$-algebra in its target
category.  This applies, for example, to the functor $\THH_{\bullet}(-)$
from associative $S$-algebras to cyclic $S$-modules, and to its
geometric realization $\THH(-)$ from associative $S$-algebras to
$S^1$-equivariant $S$-modules.  Therefore $\THH$ of an $E_{n+1}$ ring
spectrum $R$ is an $S^1$-equivariant $E_n$ ring spectrum, and the unit
$\eta \: R \to \THH(R)$ is a map of $E_n$ ring spectra.  We do not
discuss whether an $E_{n+1}$ ring spectrum structure on $R$, for some
finite $n$, suffices to define a more-or-less strictly coassociative
coproduct or Hopf algebra structure on $\THH(R)$.

\begin{prop}
\label{p4.3}
Let $R$ be an $E_2$ ring spectrum.  Then the B{\"o}kstedt spectral
sequence $E^r_{**}(R)$ converging to $H_*(\THH(R); \F_p)$ is an
$A_*$-comodule $\F_p$-algebra spectral sequence.  If $R$ is an
$E_3$ ring spectrum, then $E^r_{**}(R)$ is a commutative $H_*(R;
\F_p)$-algebra spectral sequence in $A_*$-comodules.
\end{prop}

\begin{proof}
To make sense of $\THH(R)$, we have implicitly replaced $R$ with a
weakly equivalent associative $S$-algebra.  By the discussion above,
for $R$ an $E_2$ ring spectrum there is a $\CalC_1$-algebra action on
the cyclic $S$-module $\THH_\bullet(R)$, on the $A_*$-comodule
$H_*(\THH(R); \F_p)$ and on the B{\"o}kstedt spectral sequence
$E^r_{**}(R)$.  This $\CalC_1$-action induces the desired
$A_*$-comodule $\F_p$-algebra structures.  For example, the spectral
sequence pairing is the composite
$$
\xymatrix{
E^r_{**}(R) \otimes_{\F_p} E^r_{**}(R)
\to {}' E^r_{**} \ar[r]^-{\mu} & E^r_{**}(R) \,,
}
$$
with ${}'E^r_{**}$ the spectral sequence associated to $\THH(R) \wedge
\THH(R)$ and $\mu$ induced by any point in $\CalC_1(2)$.  When $R$ is an
$E_3$ ring spectrum, $\THH(R)$ is an $E_2$ ring spectrum and homotopy
commutative, so the algebra structures are commutative.  The unit
$S$-algebra map $\eta \: R \to \THH(R)$ is then homotopy central, and
supplies the $H_*(R; \F_p)$-algebra structure.
\end{proof}

Here is a typical application of this algebra structure.

\begin{cor}
\label{c4.4}
If $R$ is an $E_2$ ring spectrum and $H_*(R; \F_p)$ is a polynomial
algebra over~$\F_p$, then $E^r_{**}(R)$ collapses at the $E^2$-term, so
$E^2_{**}(R) = E^\infty_{**}(R)$.  Furthermore, there are no nontrivial
(left or right) $H_*(R; \F_p)$-module extensions.
\end{cor}

\begin{proof}
If $H_*(R; \F_p) = P(x_i)$ is a polynomial algebra, where $i$ ranges
through some indexing set $I$, then $E^2_{**}(R) = HH_*(P(x_i)) =
P(x_i) \otimes E(\sigma x_i)$ by Proposition~\ref{p2.4} (and passage to
colimits).  All of the $\F_p$-algebra generators are in
filtration~$\le1$, so all differentials on these classes are zero,
since the B{\"o}kstedt spectral sequence is a right half plane
homological spectral sequence.  Thus there are no further differentials
in this spectral sequence, and $E^2_{**}(R) = E^\infty_{**}(R)$.
Furthermore, the $E^\infty$-term is a free $H_*(R; \F_p)$-module, so
also $H_*(\THH(R); \F_p) \cong H_*(R; \F_p) \otimes E(\sigma x_i)$ is a
free $H_*(R; \F_p)$-module.
\end{proof}

There may, as we shall see in Section~5, be multiplicative extensions
between $E^\infty_{**}(R)$ and $H_*(\THH(R); \F_p)$, as well as
$A_*$-comodule extensions.

A flatness hypothesis is required for the spectral sequence to carry
the coproduct and full Hopf algebra structure.  Our Sections~5, 6 and~7
will show many examples of B{\"o}kstedt spectral sequences with this
structure.

\begin{thm}
\label{t4.5}
Let $R$ be a commutative $S$-algebra and write $\Lambda = H_*(R; \F_p)$.

{\rm (a)}\qua
If $H_*(\THH(R); \F_p)$ is flat over~$\Lambda$, then there is a
coproduct
$$
\psi \: H_*(\THH(R); \F_p)
\to H_*(\THH(R); \F_p) \otimes_\Lambda H_*(\THH(R); \F_p)
$$
and $H_*(\THH(R); \F_p)$ is an $A_*$-comodule $\Lambda$-Hopf algebra.

{\rm (b)}\qua
If each term $E^r_{**}(R)$ for $r \ge 2$
is flat over~$\Lambda$, then there is a coproduct
$$
\psi \: E^r_{**}(R) \to E^r_{**}(R) \otimes_\Lambda E^r_{**}(R)
$$
and $E^r_{**}(R)$ is an $A_*$-comodule $\Lambda$-Hopf algebra spectral
sequence.  In particular, the differentials $d^r$ respect the coproduct
$\psi$.
\end{thm}
 
\begin{proof}
Write $T = \THH(R)$.
There is a K{\"u}nneth spectral sequence with
$$
E^2_{**} = \Tor^{\Lambda}_{**}(H_*(T; \F_p), H_*(T; \F_p))
\Longrightarrow H_*(T \wedge_R T; \F_p) \,.
$$
When $H_*(T; \F_p)$ is flat over~$\Lambda$ the higher
$\Tor$-groups vanish, the spectral sequence collapses, and
the map
$$
\psi \: H_*(T; \F_p) \to H_*(T \wedge_R T; \F_p)
\cong H_*(T; \F_p) \otimes_\Lambda H_*(T; \F_p)
$$
induces the coproduct in part~(a).

For part~(b) let $dE^r_{**}$ be the spectral sequence associated to the
skeleton filtration on $d\THH(R) = R \otimes dS^1$, and let
${}'E^r_{**}$ be the spectral sequence associated to $T \wedge_R T = R
\otimes (S^1 \vee S^1)$, as in the proof of Proposition~\ref{p4.2}.
Then there are natural maps of spectral sequences
\begin{equation}
\xymatrix{
E^r_{**}(R) & \ar[l]_-{\pi} dE^r_{**} \ar[r]^-{\psi'} & {}'E^r_{**} &
\ar[l]_-{sh} E^r_{**}(R) \otimes_\Lambda E^r_{**}(R) \,.
}
\label{e4.6}
\end{equation}
Here $\pi \: dE^r_{**} \to E^r_{**}(R)$ is an isomorphism for $r \ge
2$, by the algebraic analogue of Lemma~\ref{l3.8}.  The map $\psi' \:
dE^r_{**} \to {}'E^r_{**}$ is induced by the simplicial pinch map
$\psi'$ from~(\ref{e3.7}).  As regards the final (shuffle) map $sh$, we have
$$
{}'E^1_{s,*} \cong \Lambda^{\otimes(s+1)} \otimes_\Lambda
\Lambda^{\otimes(s+1)} \cong E^1_{s,*}(R) \otimes_\Lambda E^1_{s,*}(R)
$$
by the collapsing K{\"u}nneth spectral sequence for
$$
H_*(\THH_s(R) \wedge_R \THH_s(R); \F_p) \,,
$$
and the map $sh$ for $r=1$ is the shuffle equivalence from the bigraded
tensor product $[E^1_{**}(R) \otimes_\Lambda E^1_{**}(R)]_{s,*}$.  By
assumption $E^2_{**}(R) = HH_*(\Lambda)$ is flat over~$\Lambda$, so by
the algebraic K{\"u}nneth spectral sequence and the Eilenberg--Zilber
theorem
$$
\xymatrix{
sh \: E^2_{**}(R) \otimes_\Lambda E^2_{**}(R)
\cong H_*(E^1_{**}(R) \otimes_\Lambda E^1_{**}(R))
\ar[r]^-{\cong} & H_*({}'E^1_{**}) = {}'E^2_{**}
}
$$
is an isomorphism.  Inductively, suppose that $sh$ is an isomorphism
for a fixed $r\ge2$, and that $E^{r+1}_{**}(R)$ is flat over~$\Lambda$.
Then by the algebraic K{\"u}nneth spectral sequence again
$$
\xymatrix{
sh \: E^{r+1}_{**}(R) \otimes_\Lambda E^{r+1}_{**}(R)
\cong H_*(E^r_{**}(R) \otimes_\Lambda E^r_{**}(R))
\ar[r]^-{\cong} & H_*({}'E^r_{**}) = {}'E^{r+1}_{**}
}
$$
is also an isomorphism, as desired.  The coproduct $\psi$ on $E^r_{**}(R)$
is then the composite map $(sh)^{-1} \circ \psi' \circ \pi^{-1}$, for
$r \ge 2$.  The conjugation $\chi$ on $E^r_{**}(R)$ is more
simply defined, as the composite map $\pi \circ \chi' \circ \pi^{-1}$.
\end{proof}

By Proposition~\ref{p2.3}, the coproduct (and Hopf algebra structure)
on the $E^2$-term $E^2_{**}(R) = HH_*(H_*(R; \F_p))$ that is derived
from the $R$-Hopf algebra structure on $\THH(R)$ agrees with the
algebraically defined structure on the Hochschild homology
$HH_*(\Lambda)$ of the commutative algebra $\Lambda = H_*(R; \F_p)$.

\begin{rem}
\label{r4.7}
The same proof shows that if only the initial terms
$$
E^2_{**}(R), \dots, E^{r_0}_{**}(R)
$$
of the B{\"o}kstedt spectral sequence are flat over~$\Lambda$, for some
integer $r_0$, then these are all $A_*$-comodule $\Lambda$-Hopf
algebras and the corresponding differentials $d^2, \dots, d^{r_0}$ all
respect that structure.
\end{rem}

There are natural examples that show that the flatness hypothesis is
not always realistic.  For example, Ausoni \cite{Au} has studied the
case of $R = ku$ at an odd prime $p$, where $H_*(ku; \F_p) = H_*(\ell;
\F_p) \otimes P_{p-1}(x)$ for a class $x$ in degree~$2$, and already
the $E^2$-term $E^2_{**}(ku) = HH_*(H_*(ku; \F_p))$ is not flat over
$H_*(ku; \F_p)$.
We shall also see in Proposition~\ref{p7.13}(c) that for $R = j$, the
connective, real image-of-$J$ spectrum at $p=2$, the $E^2$-term
$E^2_{**}(j) = HH_*(H_*(j;\F_2))$ is not flat over $H_*(j; \F_2)$.

We say that the graded $k$-algebra $\Lambda$ is {\it connected\/} when
it is trivial in negative degrees and the unit map $\eta \: k \to
\Lambda$ is an isomorphism in degree~$0$.  We write $\nu \: H_*(X;
\F_p) \to A_* \otimes H_*(X; \F_p)$ for the $A_*$-comodule coaction
map, and say that a class $x \in H_*(X; \F_p)$ is {\it $A_*$-comodule
primitive\/} if $\nu(x) = 1 \otimes x$.

\begin{prop}
\label{p4.8}
Let $R$ be a commutative $S$-algebra with $\Lambda = H_*(R; \F_p)$
connected and such that $HH_*(\Lambda)$ is flat over~$\Lambda$.  Then the
$E^2$-term of the first quadrant B{\"o}kstedt spectral sequence
$$
E^2_{**}(R) = HH_*(\Lambda)
$$
is an $A_*$-comodule $\Lambda$-Hopf algebra, and a shortest non-zero
differential $d^r_{s,t}$ in lowest total degree $s+t$, if one exists,
must map from an algebra indecomposable to a coalgebra primitive and
$A_*$-comodule primitive, in $HH_*(\Lambda)$.
\end{prop}

\begin{proof}
If $d^2, \dots, d^{r-1}$ are all zero, then $E^2_{**}(R) = E^r_{**}(R)$
is still an $A_*$-comodule $\Lambda$-Hopf algebra.  If $d^r(xy) \ne 0$,
where $xy$ is decomposable (a product of classes of positive degree),
then the Leibniz formula
$$
d^r(xy) = d^r(x) y \pm x d^r(y)
$$
implies that $d^r(x) \ne 0$ or $d^r(y) \ne 0$, so $xy$ cannot be in the
lowest possible total degree for the source of a differential.  Dually,
if $d^r(z)$ is not coalgebra primitive, with $\psi(z) = z \otimes 1 +
1 \otimes z + \sum_i z'_i \otimes z''_i$, then the co-Leibniz formula
$$
\psi \circ d^r = (d^r \otimes 1 \pm 1 \otimes d^r) \psi
$$
(tensor products over~$\Lambda$) implies that some term $d^r(z'_i)
\ne 0$ or $d^r(z''_i) \ne 0$, so $z$ cannot be in the lowest possible
total degree.  Finally, if $d^r(z)$ is not $A_*$-comodule primitive,
with $\nu(z) = 1 \otimes z + \sum_i a_i \otimes z_i$, then the
co-linearity condition
$$
\nu \circ d^r = (1 \otimes d^r) \nu
$$
implies that some term $d^r(z_i) \ne 0$, so $z$ cannot be in the lowest
possible total degree.
\end{proof}

The last two arguments are perhaps easier to visualize in the
$\F_p$-vector space dual spectral sequence.  The $\F_p$-dual of the
coproduct $\psi$ maps from the cotensor product
$$
HH_*(\Lambda)^* \cotensor_{\Lambda^*} HH_*(\Lambda)^*
\subset
HH_*(\Lambda)^* \otimes_{\F_p} HH_*(\Lambda)^*
$$
to $HH_*(\Lambda)^*$.

\begin{prop}
\label{p4.9}
Let $R$ be a commutative $S$-algebra.
For each element $x \in H_t(R; \F_p)$ the image $\sigma x \in
H_{t+1}(\THH(R); \F_p)$ is the coalgebra primitive
$$
\psi(\sigma x) = \sigma x \otimes 1 + 1 \otimes \sigma x
$$
that is represented in $E^\infty_{**}(R)$ by the class
$\sigma x = [1 \otimes x] \in E^2_{1,t}(R) = HH_1(\Lambda)_t$.
\end{prop}

\begin{proof}
Note that the coproduct on $\THH(R)$ is compatible under
$\omega \: R \wedge S^1_+ \to \THH(R)$ with the pinch map $R \wedge
dS^1_+ \to R \wedge (S^1 \vee S^1)_+$, and thus under $\sigma \: \Sigma
R \to \THH(R)$ with the pinch map $\Sigma R \to \Sigma R \vee \Sigma R$.
The claims then follow by inspection of the definitions in Section~3.
\end{proof}

See also \cite[3.2]{MS93} for the last claim above.	

\section{ Differentials and algebra extensions }

We now apply the B{\"o}kstedt spectral sequence~(\ref{e4.1}) to compute
the mod~$p$ homology of $\THH(R)$ for the Brown--Peterson and
Johnson--Wilson $S$-algebras $R = BP\langle m{-}1\rangle$ for $0 \le m
\le \infty$.  In each case we can replace $R$ by its $p$-localization
$R_{(p)}$ or $p$-completion $R_p$ without changing the mod~$p$ homology
of $\THH(R)$, so we will sometimes do so without further comment.  The
cases $R = H\F_p$ and $R = H\Z$ were first treated by B{\"o}kstedt
\cite{Bo2}, and the case $R = \ell_p \subset ku_p$ (the Adams summand
of $p$-complete connective $K$-theory) is due to McClure and Staffeldt
\cite[4.2]{MS93}.

\begin{sshead}
\label{s5.1}
The dual Steenrod algebra
\end{sshead}
Let $A = H^*(H\F_p; \F_p)$ be the Steenrod algebra, with generators $Sq^i$
for $p=2$ and $\beta$ and $P^i$ for $p$ odd.  We recall the structure of
its dual $A_* = H_*(H\F_p; \F_p)$ from \cite[Thm.~2]{Mi60}.  When $p=2$
we have
$$
A_* = P(\xi_k \mid k \ge 1) = P(\bar\xi_k \mid k \ge 1)
$$
where $\xi_k$ has degree $2^k-1$ and $\bar\xi_k = \chi(\xi_k)$ is the
conjugate class.  Most of the time it will be more convenient for us to
use the conjugate classes.  The coproduct is given by
$$
\psi(\bar\xi_k) = \sum_{i+j=k} \bar\xi_i \otimes \bar\xi_j^{2^i} \,,
$$
where as usual we read $\bar\xi_0$ to mean $1$.
When $p$ is odd we have
$$
A_* = P(\bar\xi_k \mid k \ge 1) \otimes E(\bar\tau_k \mid k \ge 0)
$$
with $\bar\xi_k = \chi(\xi_k)$ in degree $2(p^k-1)$ and $\bar\tau_k
= \chi(\tau_k)$ in degree $2p^k-1$.  The coproduct is given
by
$$
\psi(\bar\xi_k) = \sum_{i+j=k} \bar\xi_i \otimes \bar\xi_j^{p^i}
\qquad\text{and}\qquad
\psi(\bar\tau_k) = 1\otimes\bar\tau_k + \sum_{i+j=k} \bar\tau_i \otimes
\bar\xi_j^{p^i} \,.
$$
The mod~$p$ homology Bockstein
satisfies $\beta(\bar\tau_k) = \bar\xi_k$.

Any commutative $S$-algebra $R$ has a canonical structure as an
$E_\infty$ ring spectrum \cite[II.3.4]{EKMM97}.  In particular, its
mod~$p$ homology $H_*(R; \F_p)$ admits natural Dyer--Lashof operations
$$
Q^k \: H_*(R; \F_2) \to H_{*+k}(R; \F_2)
$$
for $p=2$ and
$$
Q^k \: H_*(R; \F_p) \to H_{* + 2k(p-1)}(R; \F_p)
$$
for $p$ odd.  Their formal properties are summarized in
\cite[III.1.1]{BMMS86}, and include Cartan formulas, Adem relations and
Nishida relations.  For $p=2$, $Q^k(x) = 0$ when $k < |x|$ and $Q^k(x)
= x^2$ when $k = |x|$.  For $p$ odd, $Q^k(x) = 0$ when $k < 2|x|$ and
$Q^k(x) = x^p$ when $k = 2|x|$.  In the special case of $R = H\F_p$,
the Dyer--Lashof operations in $A_* = H_*(H\F_p; \F_p)$ satisfy
$$
Q^{p^k}(\bar\xi_k) = \bar\xi_{k+1}
$$
for all primes $p$, and
$$
Q^{p^k}(\bar\tau_k) = \bar\tau_{k+1}
$$
for $p$ odd.  These formulas were first obtained by Leif Kristensen
(unpublished), and appeared in print in \cite[III.2.2 and
III.2.3]{BMMS86}.

When $R$ is just an $E_n$ ring spectrum, for $1 \le n < \infty$, only
some of these operations are naturally defined.  For $p=2$ the
operations $Q^k(x)$ are defined for $k \le |x| + n - 1$, and for $p$
odd they are defined for $2k \le |x| + n - 1$.  In the case of
equality, the ``top'' operation $Q^k(x)$ (denoted $\xi_{n-1}(x)$ in
{\it op.~cit.}) is not additive, and special care is called for.  We
refer to \cite[III.3.1--3]{BMMS86} for detailed information about
the formal properties that are satisfied in the homology of an $E_n$
ring spectrum.

\begin{sshead}
\label{s5.2}
The Johnson--Wilson spectra $BP\langle m{-}1\rangle$
\end{sshead}
For any prime $p$ let $BP = BP\langle\infty\rangle$ be the
Brown--Peterson spectrum and let $BP\langle m{-}1\rangle$ for $0 \le m
< \infty$ be the spectrum introduced by Johnson and Wilson in
\cite{JW73}, with mod~$p$ cohomology
$$
H^*(BP; \F_p) \cong A/\!/E = A \otimes_E \F_p
$$
and
$$
H^*(BP\langle m{-}1\rangle; \F_p) \cong A/\!/E_{m-1} =
	A \otimes_{E_{m-1}} \F_p \,.
$$
Here $E \subset A$ is the exterior subalgebra generated by the Milnor
primitives $Q_k$ for $k \ge 0$, while $E_{m-1}$ is the exterior
subalgebra generated by $Q_0, \dots, Q_{m-1}$.  These elements are
inductively defined by $Q_0 = Sq^1$ and $Q_k = [Sq^{2^k}, Q_{k-1}]$ for
$p=2$, and by $Q_0 = \beta$ and $Q_{k+1} = [P^{p^k}, Q_k]$ for $p$ odd,
see \cite[\S1]{Mi60}.  We note the special cases $BP\langle{-}1\rangle
= H\F_p$, $BP\langle0\rangle = H\Z_{(p)}$ and $BP\langle1\rangle = \ell
\subset ku_{(p)}$.  For $p=2$ we have $\ell = ku_{(2)}$.

The Brown--Peterson and Johnson--Wilson spectra have homotopy groups
$$
\pi_* BP = \Z_{(p)}[v_k \mid k \ge 1]
$$
for $m = \infty$ and
$$
\pi_* BP\langle m{-}1\rangle = \Z_{(p)}[v_1, \dots, v_{m-1}]
= \pi_* BP / (v_k \mid k \ge m)
$$
for $0 \le m < \infty$, with $v_0 = p$.  The class $v_k$ is detected
in the Adams spectral sequence $E_2^{**} = \Ext_{A_*}^{**}(\F_p,
H_*(BP; \F_p))$ for $\pi_* BP$ by the normalized cobar cocycle
$$
\sum_{i+j=k+1} [\bar\xi_i] \bar\xi_j^{2^i}
$$
for $p=2$ (the term for $i=0$ is zero) and
$$
- \sum_{i+j=k} [\bar\tau_i] \bar\xi_j^{p^i}
$$
for $p$ odd \cite[p.~63]{Ra04}.  Under the change-of-rings isomorphism
to $E_2^{**} \cong \Ext_{E_*}^{**}(\F_p, \F_p)$, where $E_* =
E(\bar\xi_k \mid k\ge1)$ for $p=2$ and $E_* = E(\bar\tau_k \mid k\ge0)$
for $p$ odd, these cobar cocycles correspond to $[\bar\xi_{k+1}]$ and
$-[\bar\tau_k]$, respectively.  Modulo decomposables, we have
$\bar\xi_{k+1} \equiv \xi_{k+1}$ for $p=2$ and $-\bar\tau_k \equiv
\tau_k$ for $p$ odd.

In each case $0 \le m \le \infty$, the spectrum $BP\langle m{-}1\rangle$
admits the structure of an associative $S$-algebra.  More precisely,
there is a tower of associative $MU$-algebras
$$
BP \to \dots \to BP\langle m\rangle \to BP\langle m{-}1\rangle \to
\dots \to H\F_p \,,
$$
which by omission of structure is also a tower of $S$-algebras.
See~\cite[3.5]{BJ02}.  It is well known that $H\F_p$ and $H\Z_{(p)}$
admit unique structures as commutative $S$-algebras, and that the
$p$-complete Adams summand $\ell_p \subset ku_p$ admits at least one
such structure \cite[\S9]{MS93}.  In fact each of $\ell$, $ku$,
$\ell_p$ and $ku_p$ admits a unique commutative $S$-algebra structure
\cite{BaRi}.

It remains a well-known open problem whether $BP$ admits a commutative
$S$-algebra structure.  Basterra and Mandell (unpublished) have shown
that $BP$ is at least an $E_4$ ring spectrum.  We shall make use of the
existence of an $E_3$ ring spectrum structure on $BP$ for the
calculations in this section.

For $p=2$ and $2 < m < \infty$, Strickland \cite[6.5]{Str99} has shown
that $BP\langle m{-}1\rangle$ is not a homotopy commutative $MU$-ring
spectrum with the most common (e.g.~Araki or Hazewinkel) choices of
maps $MU \to BP\langle m{-}1\rangle$, but he also shows that there is a
less familiar replacement $BP\langle m{-}1\rangle'$ that may serve as a
substitute \cite[2.10]{Str99}.  We do not know for which $p$ and $m$
the spectrum $BP\langle m{-}1\rangle$ admits an $E_2$- or $E_3$ ring
spectrum structure.

\begin{prop}
\label{p5.3}
For each $p$ and $0 \le m \le \infty$ the unique $S$-algebra map
$BP\langle m{-}1\rangle \to H\F_p$ induces the identifications
$$
H_*(BP\langle m{-}1\rangle; \F_2) =
P(\bar\xi_1^2, \dots, \bar\xi_m^2, \bar\xi_k \mid k \ge m+1)
$$
when $p=2$, and
$$
H_*(BP\langle m{-}1\rangle; \F_p) =
P(\bar\xi_k \mid k\ge1) \otimes E(\bar\tau_k \mid k \ge m)
$$
when $p$ is odd, as $\F_p$-subalgebras of $A_*$.
\end{prop}

\begin{proof}
See \cite[1.7]{Wi75} and dualize.  In the case $m = \infty$, the
notation is meant to imply that there are no classes $\bar\xi_k$ with
$k \ge m+1$, and no classes $\bar\tau_k$ with $k \ge m$.
\end{proof}

In particular,
\begin{align*}
H_*(H\Z; \F_2) &= P(\bar\xi_1^2, \bar\xi_k \mid k \ge 2) \\
H_*(ku; \F_2) &= P(\bar\xi_1^2, \bar\xi_2^2, \bar\xi_k \mid k \ge 3) \\
H_*(BP; \F_2) &= P(\bar\xi_k^2 \mid k \ge 1)
\end{align*}
for $p=2$, and
\begin{align*}
H_*(H\Z; \F_p) &= P(\bar\xi_k \mid k\ge1) \otimes E(\bar\tau_k \mid k \ge 1) \\
H_*(\ell; \F_p) &= P(\bar\xi_k \mid k\ge1) \otimes E(\bar\tau_k \mid k \ge 2) \\
H_*(BP; \F_p) &= P(\bar\xi_k \mid k \ge 1)
\end{align*}
for $p$ odd.

The $E^2$-term of the B{\"o}kstedt spectral sequence for $BP\langle
m{-}1\rangle$ can now be computed from Proposition~\ref{p2.4}.  It is
$$
E^2_{**}(BP\langle m{-}1\rangle) =
H_*(BP\langle m{-}1\rangle; \F_2) \otimes E(\sigma\bar\xi_1^2, \dots,
\sigma\bar\xi_m^2, \sigma\bar\xi_k \mid k \ge m+1)
$$
when $p=2$, and
$$
E^2_{**}(BP\langle m{-}1\rangle) =
H_*(BP\langle m{-}1\rangle; \F_p) \otimes E(\sigma\bar\xi_k \mid k \ge 1)
\otimes \Gamma(\sigma\bar\tau_k \mid k \ge m)
$$
when $p$ is odd.  Note that in all cases for $p=2$, as well as the case
$m = \infty$ for $p$ odd, the $E^2$-term is generated as an
$\F_p$-algebra by classes in filtration~$\le 1$.

\begin{sshead}
\label{s5.4}
Odd-primary differentials
\end{sshead}
In the remaining cases, with $p$ odd and $0 \le m < \infty$, there are
non-trivial $d^{p-1}$-differentials in the B{\"o}kstedt spectral
sequence.  These can all be determined by naturality with respect to
the map of $S$-algebras $BP\langle m{-}1\rangle \to H\F_p$, since the
map of $E^2$-terms
$$
E^2_{**}(BP\langle m{-}1\rangle) \to E^2_{**}(H\F_p) = A_* \otimes
E(\sigma\bar\xi_k \mid k \ge 1) \otimes \Gamma(\sigma\bar\tau_k \mid k
\ge 0)
$$
is injective.  Now $H\F_p$ is a connective, commutative $S$-algebra,
and the $E^2$-term is free, hence flat, over $A_* = H_*(H\F_p; \F_p)$,
so we are in the situation of Proposition~\ref{p4.8}.  The
$A_*$-algebra generators (indecomposables) of the $E^2$-term are in
filtrations~$p^i$ for $i\ge0$, by formula~(\ref{e2.5}), and the
$A_*$-coalgebra primitives are all in filtration~$1$, by
Proposition~\ref{p2.4}, so the shortest non-zero differential in lowest
total degree, if any, must go from some filtration~$p^i$ for $i\ge1$ to
filtration~$1$.

Indeed, B{\"o}kstedt \cite[1.3]{Bo2} found that in the case $R =
H\F_p$ there are non-trivial $d^{p-1}$-differentials in his spectral
sequence computing $H_*(\THH(R); \F_p)$.  In our notation, they are
given by the formula
\begin{equation}
d^{p-1}(\gamma_j(\sigma\bar\tau_k)) = \sigma\bar\xi_{k+1} \cdot
\gamma_{j-p}(\sigma\bar\tau_k)
\label{e5.5}
\end{equation}
for $j \ge p$, up to a unit in $\F_p$.  This way of writing
B{\"o}kstedt's formula first appears in \cite[p.~21]{MS93}.  A proof of
a more general result that implies~(\ref{e5.5}) was published by Hunter
\cite[Thm.~1]{Hu96}, as recalled below.

\begin{prop}
\label{p5.6}
Let $R$ be a commutative $S$-algebra, and let $x \in H_{2i-1}(R; \F_p)$
for some $i\ge1$.  Then in the B{\"o}kstedt spectral sequence
$E^r_{**}(R)$ the differentials $d^r$ all vanish for $2 \le r \le p-2$,
and there is a differential
$$
d^{p-1}(\gamma_p(\sigma x)) = \sigma(\beta Q^i(x))
$$
up to a unit in $\F_p$.
\end{prop}

A short proof of B{\"o}kstedt's formula in particular cases was later
found by Ausoni \cite[4.3]{Au}. More recently, the first author has
given a construction of $\THH$ of any $A_\infty$ ring spectrum as a kind
of geometric realization using cyclohedra \cite{Sta97}, similar to
Stasheff's bar construction of an $A_\infty$ $H$-space in \cite{Sta63}.
This construction exhibits the $d^{p-1}$-differential on
$\gamma_p(\sigma x)$ as being represented by $\sigma(\langle x,\dots, x
\rangle)$ in the B{\"o}kstedt spectral sequence, where $\langle x,
\dots, x \rangle \in H_{2pi-2}(R; \F_p)$ denotes the $p$-fold Massey
product.  In this formulation only the strict associativity of $R$ is
needed.

The proposition above applies with $R = H\F_p$, $x =\bar\tau_k$, $i =
p^k$ and $\beta Q^i(x) = \beta \bar\tau_{k+1} = \bar\xi_{k+1}$, to
establish B{\"o}kstedt's formula~(\ref{e5.5}) in the special case of $j
= p$.  The general case follows from this by induction on $j \ge p$ and
the coalgebra structure on the B{\"o}kstedt spectral sequence.  In more
detail, a comparison of $\psi d^{p-1}(\gamma_j(\sigma\bar\tau_k)) =
(d^{p-1} \otimes 1 + 1 \otimes d^{p-1})(\psi
\gamma_j(\sigma\bar\tau_k))$ and $\psi(\sigma\bar\xi_{k+1} \cdot
\gamma_{j-p}(\sigma\bar\tau_k))$ shows that the difference
$d^{p-1}(\gamma_j(\sigma\bar\tau_k)) - \sigma\bar\xi_{k+1} \cdot
\gamma_{j-p}(\sigma\bar\tau_k)$ must be a coalgebra primitive, and
there are none such other than zero in its bidegree when $j > p$.

By naturality, the formula~(\ref{e5.5}) for $d^{p-1}$ holds also in the
B{\"o}kstedt spectral sequence for $BP\langle m{-}1\rangle$ at odd
primes~$p$, whether this $S$-algebra is commutative or not.  We view
its $E^2 = E^{p-1}$-term as the tensor product of the complexes
$E(\sigma\bar\xi_{k+1}) \otimes \Gamma(\sigma\bar\tau_k)$ for $k \ge m$
and the remaining terms $H_*(BP\langle m{-}1\rangle; \F_p) \otimes
E(\sigma\bar\xi_1, \dots, \sigma\bar\xi_m)$.  Applying the K{\"u}nneth
formula, we compute its homology to be
$$
E^p_{**}(BP\langle m{-}1\rangle) =
H_*(BP\langle m{-}1\rangle; \F_p) \otimes
E(\sigma\bar\xi_1, \dots, \sigma\bar\xi_m)
\otimes P_p(\sigma\bar\tau_k \mid k \ge m) \,.
$$
At this point we note that also for $0 \le m < \infty$ and $p$ odd, the
$E^p$-term of the B{\"o}kstedt spectral sequence is generated as an
$\F_p$-algebra by classes in filtration~$\le 1$.

To proceed, we shall need to assume that the B{\"o}kstedt spectral
sequence is an algebra spectral sequence.  We know that the hypothesis
of the following proposition is satisfied in the commutative cases $0
\le m \le 2$ and in the $E_4$ ring spectrum case $m = \infty$.

\begin{prop}
\label{p5.7}
Let $p$ and $0 \le m \le \infty$ be such that $BP\langle m{-}1\rangle$
admits the structure of an $E_2$ ring spectrum.  Then its B{\"o}kstedt
spectral sequence collapses at the $E^p = E^\infty$-term, which equals
$$
E^\infty_{**}(BP\langle m{-}1\rangle) =
H_*(BP\langle m{-}1\rangle; \F_2) \otimes E(\sigma\bar\xi_1^2, \dots,
\sigma\bar\xi_m^2, \sigma\bar\xi_k \mid k \ge m+1)
$$
when $p=2$, and
$$
E^\infty_{**}(BP\langle m{-}1\rangle) =
H_*(BP\langle m{-}1\rangle; \F_p) \otimes
E(\sigma\bar\xi_1, \dots, \sigma\bar\xi_m)
\otimes P_p(\sigma\bar\tau_k \mid k \ge m)
$$
when $p$ is odd.  In each case the $E^\infty$-term exhibits
$H_*(\THH(BP\langle m{-}1\rangle); \F_p)$ as a free $H_*(BP\langle
m{-}1\rangle; \F_p)$-module.
\end{prop}

\begin{proof}
By Proposition~\ref{p4.3} we know that the B{\"o}kstedt spectral
sequence is an $\F_p$-algebra spectral sequence, and by the
calculations above we know that its $E^p$-term is generated as an
$\F_p$-algebra by classes in filtrations~$\le1$, which must be infinite
cycles since the spectral sequence is concentrated in the first
quadrant.  Thus there is no room for any further differentials and $E^p
= E^\infty$.  In each case the $E^\infty$-term is free as an
$H_*(BP\langle m{-}1\rangle; \F_p)$-module, so there is no room for
additive (module-)extensions, either.
\end{proof}

\begin{sshead}
\label{s5.8}
Algebra extensions
\end{sshead}
To determine $H_*(\THH(BP\langle m{-}1\rangle); \F_p)$ as an
$H_*(BP\langle m{-}1\rangle; \F_p)$-algebra, we need to resolve the
possible multiplicative extensions.  For this we will use the first
nontrivial Dyer--Lashof operation after the Pontryagin power, the
existence of which presupposes that $\THH(BP\langle m{-}1\rangle)$ is an
$E_2$ ring spectrum.  We shall therefore assume that $p$ and $m$ are
such that $BP\langle m{-}1\rangle$ is an $E_3$ ring spectrum.

For a commutative $S$-algebra $R$, the map $\sigma \: \Sigma R \to
\THH(R)$ relates the Dyer--Lashof operations on $R$ to those of $\THH(R)$
by the following formula of B{\"o}kstedt \cite[2.9]{Bo2}:
$$
Q^k(\sigma x) = \sigma Q^k(x) \,.
$$
See also \cite[p.~22]{MS93}.  There appears to be no published proof of
this key relation, so we offer the following slightly more general
result.

Recall from Definition~\ref{d2.6} the $S^1$-action map $\alpha \:
\THH(R) \wedge S^1_+ \to \THH(R)$, inducing the homomorphism
\begin{multline*}
\alpha \: H_*(\THH(R); \F_p) \otimes H_*(S^1_+; \F_p) \cong \\
H_*(\THH(R) \wedge S^1_+; \F_p) \to H_*(\THH(R); \F_p)
\end{multline*}
in homology.  Let $s_1 \in H_1(S^1_+; \F_p)$ be the canonical
generator.  Also recall that any commutative $S$-algebra $R$ is an
$E_\infty$ ring spectrum, so that the following proposition applies to
such an $R$ with $n = \infty$, and thus for all integers $k$.

\begin{prop}
\label{p5.9}
Let $R$ be an $E_{n+1}$ ring spectrum, so that $\THH(R)$ is an $E_n$
ring spectrum.  Then we have
$$
Q^k(\alpha(x \otimes s_1)) = \alpha(Q^k(x) \otimes s_1)
$$
for all classes $x \in H_*(\THH(R); \F_p)$ and integers $k$ such that
the operation $Q^k(x)$ is naturally defined, i.e., for $k \le |x| + n -
1$ when $p=2$ and for $2k \le |x| + n - 1$ when $p$ is odd.  In
particular, we have
$$
Q^k(\sigma x) = \sigma Q^k(x)
$$
for all classes $x \in H_*(R; \F_p)$ and the same integers $k$.
\end{prop}

\begin{proof}
By \cite[3.4]{FV}, or the definition $\THH(R) \cong R \otimes S^1$ in
the strictly commutative case, the circle acts on $\THH(R)$ by
$\CalC_n$-algebra maps.  Thus the right adjoint map
$$
\widetilde\alpha \: \THH(R) \to F(S^1_+, \THH(R))
$$
is a map of $\CalC_n$-algebras, where the product structure on the
right is given by pointwise multiplication, using the strictly
cocommutative diagonal map $S^1_+ \to S^1_+ \wedge S^1_+$.

Let $DS^1_+ = F(S^1_+, S)$ be the functional dual of $S^1_+$, also with
the pointwise multiplication.  It has mod~$p$ homology $H_*(DS^1_+;
\F_p) \cong H^{-*}(S^1_+; \F_p) = E(\iota_1)$ for a canonical class
$\iota_1 \in H^1(S^1_+; \F_p)$ dual to the class $s_1 \in H_1(S^1_+;
\F_p)$.  The Dyer--Lashof operations $Q^k$ on $H_*(DS^1_+; \F_p)$
correspond \cite[III.1.2 and VIII.3]{BMMS86} to the Steenrod operations
$P^{-k}$ on $H^{-*}(S^1_+; \F_p)$, hence are trivial for $k \ne 0$.
The same conclusions hold for $p=2$, with somewhat different notation.

There is a canonical map of $\CalC_n$-algebras
$$
\nu \: \THH(R) \wedge DS^1_+ \to F(S^1_+, \THH(R)) \,,
$$
given by the composition of functions
$$
F(S, \THH(R)) \wedge F(S^1_+, S) \to F(S^1_+, \THH(R)) \,.
$$
Compare \cite[III.1]{LMS86}.  The map $\nu$ is an equivalence since
$S^1_+$ is a finite CW complex (i.e., by Spanier--Whitehead duality).
Hence there are homomorphisms
$$
\xymatrix{
H_*(\THH(R); \F_p) \ar[r]^-{\widetilde\alpha} & 
H_*(F(S^1_+, \THH(R)); \F_p) \\
& \ar[u]_-{\nu}^-{\cong} H_*(\THH(R); \F_p) \otimes H_*(DS^1_+; \F_p)
}
$$
that take $x \in H_*(\THH(R); \F_p)$ to
$$
\nu^{-1}\widetilde\alpha(x) = x \otimes 1 + \alpha(x \otimes s_1)
\otimes \iota_1\,.
$$
Since $\widetilde\alpha$ and $\nu$ are maps of $\CalC_n$-algebras,
we have
$$
Q^k(\nu^{-1}\widetilde\alpha(x)) = \nu^{-1}\widetilde\alpha(Q^k(x))
$$
when $k$ is such that $Q^k$ is naturally defined on $x \in H_*(\THH(R);
\F_p)$.  The external Cartan formula for Dyer--Lashof operations
\cite[III.1.1(6)]{BMMS86} then gives us
$$
Q^k(x \otimes 1 + \alpha(x \otimes s_1) \otimes \iota_1)
= Q^k(x) \otimes 1 + Q^k(\alpha(x \otimes s_1)) \otimes \iota_1
$$
since $Q^i(1) = 0$ and $Q^i(\iota_1) = 0$ for $i \ne 0$.  Stringing
these formulas together yields
$$
Q^k(x) \otimes 1 + Q^k(\alpha(x \otimes s_1)) \otimes \iota_1 =
Q^k(x) \otimes 1 + \alpha(Q^k(x) \otimes s_1) \otimes \iota_1
\,,
$$
and we can read off $Q^k(\alpha(x \otimes s_1)) = \alpha(Q^k(x) \otimes
s_1)$, as desired.

Specializing to classes $\eta(x) \in H_*(\THH(R); \F_p)$ that are images
under the $E_n$ ring spectrum map $\eta \:  R \to \THH(R)$ of elements
$x \in H_*(R; \F_p)$, we have $\sigma x = \alpha(\eta(x) \otimes s_1)$
and $\sigma Q^k(x) = \alpha(\eta(Q^k(x)) \otimes s_1) =
\alpha(Q^k(\eta(x)) \otimes s_1)$.  Thus we obtain B{\"o}kstedt's
formula $Q^k(\sigma x) = \sigma Q^k(x)$, valid for the same integers
$k$ as above.
\end{proof}

The same ideas can be used to prove that $\sigma \: H_*(R; \F_p) \to
H_{*+1}(\THH(R); \F_p)$ is a graded derivation, when $R$ is sufficiently
commutative for this to make sense.

\begin{prop}
\label{p5.10}
Let $R$ be an $E_2$ ring spectrum, so that $\THH(R)$ is an $S$-algebra.
Then we have
$$
\alpha(xy \otimes s_1) = x \cdot \alpha(y \otimes s_1) +
(-1)^{|y|} \alpha(x \otimes s_1) \cdot y
$$ 
for $x, y \in H_*(\THH(R); \F_p)$.  In particular, we have
the Leibniz rule
$$
\sigma(x \cdot y) = x \cdot \sigma(y) + (-1)^{|y|} \sigma(x) \cdot y
$$
for $x, y \in H_*(R; \F_p)$.
\end{prop}

\begin{proof}
We keep the notation of the proof of Proposition~\ref{p5.9}.  Since
$\tilde\alpha$ and $\nu$ are maps of $\CalC_1$-algebras,
$$
\nu^{-1}\tilde\alpha \: H_*(\THH(R); \F_p) \to H_*(\THH(R); \F_p)
\otimes E(\iota_1)
$$
is an $\F_p$-algebra homomorphism.  Thus
$$
xy \otimes 1 + \alpha(xy \otimes s_1) \otimes \iota_1
= (x \otimes 1 + \alpha(x \otimes s_1) \otimes \iota_1)
\cdot (y \otimes 1 + \alpha(y \otimes s_1) \otimes \iota_1)
$$
for $x, y \in H_*(\THH(R); \F_p)$.  Multiplying out and comparing
$\iota_1$-coefficients gives the claimed formulas.
\end{proof}

We also wish to describe the $A_*$-comodule structure on $H_*(\THH(R);
\F_p)$.  In the cases $R = BP\langle m{-}1\rangle$ the following observations
will suffice.  The $A_*$-comodule coaction map
$$
\nu \: H_*(R; \F_p) \to A_* \otimes H_*(R; \F_p)
$$
is in each case given by restricting the coproduct
$$
\psi \: A_* \to A_* \otimes A_*
$$
given in Subsection~5.1, to the subalgebra $H_*(R; \F_p) \subset A_*$,
since the latter inclusion is induced by a spectrum map $R \to H\F_p$.
The operator
$$
\sigma \: H_*(\Sigma R; \F_p) \to H_*(\THH(R); \F_p)
$$
is also induced by a spectrum map, hence is an $A_*$-comodule
homomorphism.  Hence the coaction map
$$
\nu \: H_*(\THH(R); \F_p) \to A_* \otimes H_*(\THH(R); \F_p)
$$
satisfies
\begin{equation}
\nu \circ \sigma = (1 \otimes \sigma) \nu \,.
\label{e5.11}
\end{equation}

As before, we know that the hypotheses of the following proposition are
all satisfied in the commutative cases $0 \le m \le 2$, and for
parts~(a) and~(c) in the $E_4$ ring spectrum case $m = \infty$.

\begin{thm}
\label{t5.12}
{\rm (a)}\qua
Let $p$ and $0\le m\le\infty$ be such that $BP\langle m{-}1\rangle$
admits the structure of an $E_3$ ring spectrum.  Then for $m < \infty$
we have
\begin{multline*}
H_*(\THH(BP\langle m{-}1\rangle); \F_2) = \\
H_*(BP\langle m{-}1\rangle; \F_2)
\otimes E(\sigma\bar\xi_1^2, \dots, \sigma\bar\xi_m^2)
\otimes P(\sigma\bar\xi_{m+1})
\end{multline*}
when $p=2$, and
\begin{multline*}
H_*(\THH(BP\langle m{-}1\rangle); \F_p) = \\
H_*(BP\langle m{-}1\rangle; \F_p) \otimes
E(\sigma\bar\xi_1, \dots, \sigma\bar\xi_m)
\otimes P(\sigma\bar\tau_m)
\end{multline*}
when $p$ is odd, as $H_*(BP\langle m{-}1\rangle; \F_p)$-algebras.
For $m=\infty$ we have
$$
H_*(\THH(BP); \F_2) = H_*(BP; \F_2) \otimes
E(\sigma\bar\xi_k^2 \mid k \ge 1)
$$
when $p=2$, and
$$
H_*(\THH(BP); \F_p) =
H_*(BP; \F_p) \otimes
E(\sigma\bar\xi_k \mid k \ge 1)
$$
when $p$ is odd, as $H_*(BP; \F_p)$-algebras.

{\rm (b)}\qua
When $p$ and $m$ are such that $BP\langle m{-}1\rangle$ is a
commutative $S$-algebra, then these are isomorphisms of primitively
generated $H_*(BP\langle m{-}1\rangle; \F_p)$-Hopf algebras.

{\rm(c)}\qua
The $A_*$-comodule coaction on $H_*(\THH(BP\langle m{-}1\rangle); \F_p)$
is given on the tensor factor $H_*(BP\langle m{-}1\rangle; \F_p)$ by
restricting the coproduct on $A_*$.  For $p=2$ the algebra generators
$\sigma\bar\xi_k^2$ for $1 \le k \le m$ are $A_*$-comodule primitives,
while
$$
\nu(\sigma\bar\xi_{m+1}) = 1 \otimes \sigma\bar\xi_{m+1} + \bar\xi_1
\otimes \sigma\bar\xi_m^2 \,.
$$
For $p$ odd the algebra generators $\sigma\bar\xi_k$ for $1 \le k
\le m$ are $A_*$-comodule primitives, while
$$
\nu(\sigma\bar\tau_m) = 1 \otimes \sigma\bar\tau_m + \bar\tau_0
\otimes \sigma\bar\xi_m \,.
$$
\end{thm}

As usual, $\bar\xi_0$ is read as $1$ in such formulas.  Thus for $m=0$,
$\sigma\bar\xi_1$ and $\sigma\bar\tau_0$ are also primitive.

\begin{proof}
{\rm (a)}\qua
We first resolve the algebra extensions in the $E^\infty$-terms from
Proposition~\ref{p5.7}.  In each case the B{\"o}kstedt spectral
sequence is one of commutative $H_*(BP\langle m{-}1\rangle;
\F_p)$-algebras, by the $E_3$ ring spectrum hypothesis and
Proposition~\ref{p4.3}.  By the functoriality of the zeroth Postnikov
section, the $S$-algebra maps $BP\langle m{-}1\rangle \to H\Z_{(p)} \to
H\F_p$ and their composite are $E_3$ ring spectrum maps.  (We only need
this for $m \ne 0$.)  Since the induced map in homology is injective,
we can read off the $E_3$ ring spectrum Dyer--Lashof operations in
$H_*(BP\langle m{-}1\rangle; \F_p)$ from the formulas in $A_* =
H_*(H\F_p; \F_p)$.  This tells us that $Q^{2^{k+1}-1}(\bar\xi_k^2) = 0$
(by the Cartan formula) and $Q^{2^k}(\bar\xi_k) = \bar\xi_{k+1}$ for
$p=2$, and that $Q^{p^k}(\bar\tau_k) = \bar\tau_{k+1}$ for $p$ odd.

We will also use the formulas $x^2 = Q^k(x)$ for $|x| = k$ and $p=2$,
and $x^p = Q^k(x)$ for $|x| = 2k$ and $p$ odd, from
\cite[III.1.1(4)]{BMMS86}.   Likewise we will use the formulas
$Q^k(\sigma x) = \sigma Q^k(x)$ for $|x| = k-1$ and $p=2$, and
$Q^k(\sigma x) = \sigma Q^k(x)$ for $|x| = 2k-1$ and $p$ odd, from
Proposition~\ref{p5.9}.  These are all are valid in the homology of the
$E_2$ ring spectrum $\THH(BP\langle m{-}1\rangle)$.

With these preliminaries, we are ready to compute.  For $p=2$ the
squares in $H_*(\THH(BP\langle m{-}1\rangle); \F_2)$ of the algebra
generators in $E^\infty_{**}(BP\langle m{-}1\rangle)$ are
$$
(\sigma\bar\xi_k^2)^2 = Q^{2^{k+1}-1}(\sigma\bar\xi_k^2)
= \sigma Q^{2^{k+1}-1}(\bar\xi_k^2) = 0
$$
for $k=1, \dots, m$ and
$$
(\sigma\bar\xi_k)^2 = Q^{2^k}(\sigma\bar\xi_k) = \sigma Q^{2^k}(\bar\xi_k)
= \sigma\bar\xi_{k+1}
$$
for $k \ge m+1$.

For $p$ odd the classes $\sigma\bar\xi_k$ remain exterior, since they
are of odd degree in the graded commutative algebra $H_*(\THH(BP\langle
m{-}1\rangle); \F_p)$.  The $p$-th powers of the truncated polynomial
generators in $E^\infty_{**}(BP\langle m{-}1\rangle)$ are
$$
(\sigma\bar\tau_k)^p = Q^{p^k}(\sigma\bar\tau_k)
= \sigma Q^{p^k}(\bar\tau_k) = \sigma\bar\tau_{k+1}
$$
for $k \ge m$.  Hence these truncated polynomial algebras assemble to a
polynomial algebra on the single generator $\sigma\bar\tau_m$.

{\rm (b)}\qua
In the strictly commutative case, Proposition~\ref{p4.9} tells us that
all the algebra generators are coalgebra primitive, since they are of
the form form $\sigma x$ in Hochschild filtration~$1$.

{\rm(c)}\qua
To compute the $A_*$-comodule structure we use that $\sigma$ is a
graded derivation, in view of Proposition~\ref{p5.10}, so $\sigma(y^p)
= 0$ and $\sigma(1) = 0$.  Then for $p=2$ and $1 \le k \le m$ we have
$$
\nu(\sigma\bar\xi_k^2) = \sum_{i+j=k} \bar\xi_i^2 \otimes
\sigma\bar\xi_j^{2^{i+1}} = 1 \otimes \sigma\bar\xi_k^2
$$
while
$$
\nu(\sigma\bar\xi_{m+1}) = \sum_{i+j=m+1} \bar\xi_i \otimes
\sigma\bar\xi_j^{2^i} = 1 \otimes \sigma\bar\xi_{m+1} + \bar\xi_1 \otimes
\sigma\bar\xi_m^2 \,.
$$
For $p$ odd and $1 \le k \le m$ we get
$$
\nu(\sigma\bar\xi_k) = \sum_{i+j=k} \bar\xi_i \otimes
\sigma\bar\xi_j^{p^i} = 1 \otimes \sigma\bar\xi_k
$$
while
$$
\nu(\sigma\bar\tau_m) = 1 \otimes \sigma\bar\tau_m +
\sum_{i+j=m} \bar\tau_i \otimes \sigma\bar\xi_j^{p^i} = 1 \otimes
\sigma\bar\tau_m + \bar\tau_0 \otimes \sigma\bar\xi_m \,.
$$
\end{proof}

We can get similar conclusions under different hypotheses, by use of
naturality with respect to the maps $BP \to BP\langle m{-}1\rangle \to
H\F_p$.  Strickland's result \cite[2.10]{Str99} suggests that it may be
doable, but nontrivial, to relate a given $E_2$ ring spectrum structure
on $BP\langle m{-}1\rangle$ to the one on $BP$ by such a structured
map.

\begin{cor}
\label{c5.13}
Let $p$ and $2 < m < \infty$ be such that there is an $E_2$ ring
spectrum map $BP \to BP\langle m{-}1\rangle$.  Then there is a split
extension of associative $\F_p$-algebras
\begin{multline*}
E(\sigma\bar\xi_1^2, \dots, \sigma\bar\xi_m^2) \to \\
H_*(\THH(BP\langle m{-}1\rangle); \F_2)
\to H_*(BP\langle m{-}1\rangle; \F_2) \otimes P(\sigma\bar\xi_{m+1})
\end{multline*}
when $p=2$, and
\begin{multline*}
E(\sigma\bar\xi_1, \dots, \sigma\bar\xi_m) \to \\
H_*(\THH(BP\langle m{-}1\rangle); \F_p)
\to H_*(BP\langle m{-}1\rangle; \F_p) \otimes P(\sigma\bar\tau_m)
\end{multline*}
when $p$ is odd.

If furthermore $H_*(\THH(BP\langle m{-}1\rangle); \F_p)$ is commutative,
for which it suffices that $\THH(BP\langle m{-}1\rangle)$ is homotopy
commutative, then these extensions are trivial and
\begin{multline*}
H_*(\THH(BP\langle m{-}1\rangle); \F_2) = \\
H_*(BP\langle m{-}1\rangle; \F_2)
\otimes E(\sigma\bar\xi_1^2, \dots, \sigma\bar\xi_m^2)
\otimes P(\sigma\bar\xi_{m+1})
\end{multline*}
for $p=2$, and
\begin{multline*}
H_*(\THH(BP\langle m{-}1\rangle); \F_p) = \\
H_*(BP\langle m{-}1\rangle; \F_p)
\otimes E(\sigma\bar\xi_1, \dots, \sigma\bar\xi_m)
\otimes P(\sigma\bar\tau_m)
\end{multline*}
for $p$ odd, as commutative $H_*(BP\langle m{-}1\rangle;
\F_p)$-algebras.
\end{cor}

\begin{proof}
For brevity we write $R = BP\langle m{-}1\rangle$.
The zeroth Postnikov section provides an $E_2$ ring spectrum map $R \to
H\Z_{(p)}$, which we can continue to $H\F_p$.  The $E_2$ ring spectrum
maps $BP \to R \to H\F_p$ then induce $S$-algebra maps
$$
\THH(BP) \to \THH(R) \to \THH(H\F_p) \,.
$$
The B{\"o}kstedt spectral sequences at the two ends were determined in
Theorem~\ref{t5.12}, as the cases $m=\infty$ and $m=0$, respectively:
$H_*(\THH(BP); \F_p) = H_*(BP; \F_p) \otimes E(\sigma\bar\xi_k \mid
k\ge1)$ and $H_*(\THH(H\F_p); \F_p) = A_* \otimes P(\sigma\bar\tau_0)$,
for $p$ odd.  The case $p=2$ is similar, and will be omitted.

The left hand map embeds the exterior algebra $E(\sigma\bar\xi_1,
\dots, \sigma\bar\xi_m)$ as a graded commutative subalgebra of
$H_*(\THH(R); \F_p)$.  It maps trivially (by the augmentation) under the
second map.

The right hand map embeds $H_*(R; \F_p)$ in $A_*$.  It also embeds the
remaining factor $P_p(\sigma\bar\tau_k \mid k \ge m)$ in the
B{\"o}kstedt $E^\infty$-term as the polynomial subalgebra of
$P(\sigma\bar\tau_0)$ that is generated by the single class
$\sigma\bar\tau_m = (\sigma\bar\tau_0)^{p^m}$.  This uses the
multiplicative relations $(\sigma\bar\tau_k)^p = \sigma\bar\tau_{k+1}$
for $k\ge0$ in $H_*(\THH(H\F_p); \F_p)$.

The product of the unit map $R \to \THH(R)$ and the suspension operator
$\sigma$ provide a splitting $H_*(R; \F_p) \otimes P(\sigma\bar\tau_m)
\to H_*(\THH(R); \F_p)$.  In the commutative case, this splitting
amounts to a trivialization.
\end{proof}

From here on we only consider the topological Hochschild homology of
strictly commutative $S$-algebras, not of $E_n$ ring spectra for any
finite $n$.

For later reference we extract from Theorem~\ref{t5.12} the following
special cases, which correspond to $m=2$.  Recall that $H_*(ku; \F_2) =
(A/\!/E_1)_* \subset A_*$ and $H_*(\ell; \F_p) = (A/\!/E_1)_* \subset
A_*$.

\begin{cor}
\label{c5.14}
{\rm (a)}\qua
There is an isomorphism
$$
H_*(\THH(ku); \F_2) \cong H_*(ku; \F_2) \otimes E(\sigma\bar\xi_1^2,
\sigma\bar\xi_2^2) \otimes P(\sigma\bar\xi_3)
$$
of primitively generated $H_*(ku; \F_2)$-Hopf algebras.

The $A_*$-comodule coaction $\nu \: H_*(\THH(ku); \F_2) \to A_* \otimes
H_*(\THH(ku); \F_2)$ is given on $H_*(ku; \F_2)$ by restricting
the coproduct $\psi \: A_* \to A_* \otimes A_*$, and on the algebra
generators by $\nu(\sigma\bar\xi_1^2) = 1 \otimes \sigma\bar\xi_1^2$,
$\nu(\sigma\bar\xi_2^2) = 1 \otimes \sigma\bar\xi_2^2$ and
$$
\nu(\sigma\bar\xi_3) = 1 \otimes \sigma\bar\xi_3 + \bar\xi_1 \otimes
\sigma\bar\xi_2^2 \,.
$$

{\rm (b)}\qua
There is an isomorphism
$$
H_*(\THH(\ell); \F_p) \cong H_*(\ell; \F_p) \otimes E(\sigma\bar\xi_1,
\sigma\bar\xi_2) \otimes P(\sigma\bar\tau_2)
$$
of primitively generated $H_*(\ell; \F_p)$-Hopf algebras.

The $A_*$-comodule coaction $\nu \: H_*(\THH(\ell); \F_p) \to A_* \otimes
H_*(\THH(\ell); \F_p)$ is given on $H_*(\ell; \F_p)$ by restricting
the coproduct $\psi \: A_* \to A_* \otimes A_*$, and on the algebra
generators by $\nu(\sigma\bar\xi_1) = 1 \otimes \sigma\bar\xi_1$,
$\nu(\sigma\bar\xi_2) = 1 \otimes \sigma\bar\xi_2$ and
$$
\nu(\sigma\bar\tau_2) = 1 \otimes \sigma\bar\tau_2 + \bar\tau_0 \otimes
\sigma\bar\xi_2 \,.
$$
\end{cor}

\section{ The higher real cases }

For $p=2$ there are a few more known examples of commutative
$S$-algebras such that $H^*(R; \F_p)$ is a cyclic $A$-module.  Let $ko$
be the connective real $K$-theory spectrum, with $H^*(ko; \F_2) \cong
A/\!/A_1 = A \otimes_{A_1} \F_2$, and let $\tmf$ be the
Hopkins--Mahowald topological modular forms spectrum, with $H^*(\tmf;
\F_2) \cong A/\!/A_2 = A \otimes_{A_2} \F_2$.  See
e.g.~\cite[21.5]{Re01}.  Here $A_n \subset A$ is the subalgebra
generated by $Sq^1, \dots, Sq^{2^n}$, so $A_1$ has rank~$8$ and $A_2$
has rank~$64$.  It is well known that $ko$ is a commutative
$S$-algebra, and in the case of $\tmf$ this is a consequence of the
Hopkins--Miller theory, as being presented in the commutative case by
Goerss and Hopkins.

\begin{prop}
\label{p6.1}
There are maps of commutative $S$-algebras $ko \to H\F_2$ and $\tmf \to
H\F_2$ that induce the following identifications:
$$
H_*(ko; \F_2) = P(\bar\xi_1^4, \bar\xi_2^2, \bar\xi_k \mid k \ge 3)
$$
and
$$
H_*(\tmf; \F_2) = P(\bar\xi_1^8, \bar\xi_2^4, \bar\xi_3^2, \bar\xi_k
\mid k \ge 4) \,.
$$
\end{prop}

\begin{proof}
The maps can be obtained from the zeroth Postnikov sections $ko \to
H\Z$ and $\tmf \to H\Z$.  The homology computations are then immediate by
dualization from $H^*(ko; \F_2) \cong A/\!/A_1$, cf.~\cite{Sto63}, and
$H^*(\tmf; \F_2) \cong A/\!/A_2$.
\end{proof}

We now follow the outline of Section~5.  By Proposition~\ref{p2.4} the
$E^2$-terms of the respective B{\"o}kstedt spectral sequences are
$$
E^2_{**}(ko) = H_*(ko; \F_2) \otimes E(\sigma\bar\xi_1^4,
\sigma\bar\xi_2^2, \sigma\bar\xi_k \mid k \ge 3)
$$
and
$$
E^2_{**}(\tmf) = H_*(\tmf; \F_2) \otimes E(\sigma\bar\xi_1^8,
\sigma\bar\xi_2^4, \sigma\bar\xi_3^2, \sigma\bar\xi_k \mid k \ge 4) \,.
$$
By Corollary~\ref{c4.4} both spectral sequences collapse at the
$E^2$-term, so $E^2_{**}(R) = E^\infty_{**}(R)$.  To resolve the
algebra extensions we use the Dyer--Lashof operations and
Proposition~\ref{p5.9}.  The squares in $H_*(\THH(ko); \F_2)$ of the
algebra generators in $E^\infty_{**}(ko)$ are
\begin{align*}
(\sigma\bar\xi_1^4)^2 &= Q^5(\sigma\bar\xi_1^4) = \sigma Q^5(\bar\xi_1^4)
= 0 \\
(\sigma\bar\xi_2^2)^2 &= Q^7(\sigma\bar\xi_2^2) = \sigma Q^7(\bar\xi_2^2)
= 0
\end{align*}
by the formula $Q^k(y^2) = 0$ for $p=2$ and $k$ odd, and
$$
(\sigma\bar\xi_k)^2 = Q^{2^k}(\sigma\bar\xi_k) =
	\sigma Q^{2^k}(\bar\xi_k) = \sigma\bar\xi_{k+1}
$$
for all $k \ge 3$.  Similar calculations show that $(\sigma\bar\xi_1^8)^2
= 0$, $(\sigma\bar\xi_2^4)^2 = 0$ and $(\sigma\bar\xi_3^2)^2 = 0$ in
$H_*(\THH(\tmf); \F_2)$, while $(\sigma\bar\xi_k)^2 = \sigma\bar\xi_{k+1}$
for all $k \ge 4$.  The $A_*$-comodule coaction map $\nu$ on the
resulting algebra generators is obtained from the coproduct on $A_*$
and formula~(\ref{e5.11}), as in the proof of Theorem~\ref{t5.12}.
The result is as follows:

\begin{thm}
\label{t6.2}
{\rm (a)}\qua
There is an isomorphism
$$
H_*(\THH(ko); \F_2) \cong H_*(ko; \F_2) \otimes E(\sigma\bar\xi_1^4,
\sigma\bar\xi_2^2) \otimes P(\sigma\bar\xi_3)
$$
of primitively generated $H_*(ko; \F_2)$-Hopf algebras.
The $A_*$-comodule structure is given on $H_*(ko; \F_2)$ by
restricting the coproduct on $A_*$, and on the algebra generators by
$\nu(\sigma\bar\xi_1^4) = 1 \otimes \sigma\bar\xi_1^4$ and
\begin{align*}
\nu(\sigma\bar\xi_2^2) &= 1 \otimes \sigma\bar\xi_2^2 + \bar\xi_1^2 \otimes
\sigma\bar\xi_1^4 \\
\nu(\sigma\bar\xi_3) &= 1 \otimes \sigma\bar\xi_3 + \bar\xi_1 \otimes
\sigma\bar\xi_2^2 + \bar\xi_2 \otimes \sigma\bar\xi_1^4 \,.
\end{align*}

{\rm (b)}\qua
There is an isomorphism
$$
H_*(\THH(\tmf); \F_2) \cong H_*(\tmf; \F_2) \otimes
E(\sigma\bar\xi_1^8, \sigma\bar\xi_2^4, \sigma\bar\xi_3^2) \otimes
P(\sigma\bar\xi_4)
$$
of primitively generated $H_*(\tmf; \F_2)$-Hopf algebras.
The $A_*$-comodule structure is given on $H_*(\tmf; \F_2)$ by
restricting the coproduct on $A_*$, and on the algebra generators by
$\nu(\sigma\bar\xi_1^8) = 1 \otimes \sigma\bar\xi_1^8$ and
\begin{align*}
\nu(\sigma\bar\xi_2^4) &= 1 \otimes \sigma\bar\xi_2^4 + \bar\xi_1^4 \otimes
\sigma\bar\xi_1^8 \\
\nu(\sigma\bar\xi_3^2) &= 1 \otimes \sigma\bar\xi_3^2 + \bar\xi_1^2 \otimes
\sigma\bar\xi_2^4 + \bar\xi_2^2 \otimes \sigma\bar\xi_1^8 \\
\nu(\sigma\bar\xi_4) &= 1 \otimes \sigma\bar\xi_4 + \bar\xi_1 \otimes
\sigma\bar\xi_3^2 + \bar\xi_2 \otimes \sigma\bar\xi_2^4 + \bar\xi_3 \otimes
\sigma\bar\xi_1^8 \,.
\end{align*}
\end{thm}

\begin{proof}
Assemble the computations above.
\end{proof}

\section{ The real and complex image-of-$J$ }

We now turn to the various image-of-$J$ spectra that are commutative
$S$-alge\-bras.  Their mod~$p$ cohomology is no longer cyclic as a module
over the Steenrod algebra, but of rank~$2$, so extra work is needed to
describe their homology as an $A_*$-comodule algebra.

\newpage 

\begin{sshead}
\label{s7.1}
The image-of-$J$ spectra
\end{sshead}
For any prime $p$ let the $p$-local, connective {\it complex image-of-$J$
spectrum\/} be $ju = K(\F_r)_{(p)}$, where $r = 3$ for $p=2$ and $r$
is a prime power that topologically generates the $p$-adic units for
$p$ odd.  Being the localized algebraic $K$-theory of a field, $ju$ is a
commutative $S$-algebra \cite[VIII.3.1]{May77}.
For $p=2$ there is a cofiber sequence of spectra
\begin{equation}
\xymatrix{
ju \ar[r]^-{\kappa} & ku_{(2)} \ar[r]^-{\psi^3-1} & bu_{(2)} \,,
}
\label{e7.2}
\end{equation}
where $bu \simeq \Sigma^2 ku$ is the $1$-connected cover of $ku$
and $\psi^3$ is the Adams operation.
For odd $p$ the cofiber sequence appears as
\begin{equation}
\xymatrix{
ju \ar[r]^-{\kappa} & \ell \ar[r]^-{\psi^r-1} & \Sigma^q \ell
}
\label{e7.3}
\end{equation}
with $q = 2p-2$, where $\ell \subset ku_{(p)}$ is the $p$-local,
connective Adams summand and $\psi^r$ is the $r$-th Adams operation
\cite[V.5.16]{May77}.

Let the $2$-local, connective {\it real image-of-$J$ spectrum\/} be $j
= K\CalN(\F_3)_{(2)}$, as defined in \cite[VIII.3.1]{May77}.  Being the
localized algebraic $K$-theory of a symmetric bimonoidal category, $j$
is a commutative $S$-algebra.  There is a cofiber sequence of spectra
\begin{equation}
\xymatrix{
j \ar[r]^-{\kappa} & ko_{(2)} \ar[r]^-{\psi^3-1} & bspin_{(2)} \,,
}
\label{e7.4}
\end{equation}
where $bspin \simeq \Sigma^4 ksp$ is the $3$-connected cover of $ko$
\cite[V.5.16]{May77}.

The fiber map $\kappa \: ju \to \ell$ is a map of commutative
$S$-algebras, at least after $p$-adic completion, because there is a
discrete model $K(k')_p$ for $\ell_p$ with $k'$ a suitable subfield
of the algebraic closure of $\F_r$, and applying the functor $K(-)_p$
to the field inclusions $\F_r \subset k' \subset \overline\F_r$
produces the commutative $S$-algebra maps
$$
\xymatrix{
ju_p \ar[r]^-{\kappa} & \ell_p \subset ku_p \,.
}
$$
See \cite[VIII.3.2]{May77} and \cite[\S9]{MS93}.  Similarly, $\kappa \:
j \to ko$ becomes a map of commutative $S$-algebras after $2$-adic
completion, since there is a discrete model $K\CalO(\overline\F_3)_2$
for $ko_2$, and $\kappa$ can be identified with the natural map
$$
K\CalN(\F_3)_2 \to K\CalO(\overline\F_3)_2 \,.
$$
See \cite[VIII.2.6 and 3.2]{May77}.

\newpage 

\begin{sshead}
\label{s7.5}
Cohomology modules
\end{sshead}
Recall that
$$
H^*(\ell; \F_p) = A/\!/E_1
\cong \F_p\{1, P^1, \dots, P^p, \dots\} \,,
$$
where $E_1 = E(Q_0, Q_1) \subset A$ is the exterior algebra generated
by $Q_0 = \beta$ and $Q_1 = [P^1, \beta]$, and that
$$
H^*(ko; \F_2) = A/\!/A_1 \cong \F_2\{1, Sq^4, Sq^2 Sq^4 \equiv Sq^6,
Sq^1 Sq^2 Sq^4 \equiv Sq^7, \dots\} \,,
$$
where $A_1 = \langle Sq^1, Sq^2 \rangle \subset A$ is the subalgebra
generated by $Sq^1$ and $Sq^2$.  There are also $A$-module isomorphisms
\begin{align*}
H^*(bo; \F_2) &\cong \Sigma A/A Sq^2 \\
H^*(bso; \F_2) &\cong \Sigma^2 A/A Sq^3 \\
H^*(bspin; \F_2) &\cong \Sigma^4 A/A\{Sq^1, Sq^2Sq^3\} \,.
\end{align*}
See \cite[2.5, 2.4 and p.~501]{AP76}.  Here $Sq^2 Sq^3 = Sq^5 + Sq^4 Sq^1$
in admissible form, but the shorter expression is perhaps more memorable.
For $p$ odd we let $A_n \subset A$ be the subalgebra generated by $\beta,
P^1, \dots, P^{p^{n-1}}$.  In particular, $A_1 = \langle \beta, P^1
\rangle$ contains $E_1$, and
$$
A/\!/A_1 \cong \F_p\{1, P^p, P^1 P^p \equiv - P^{p+1}, \beta P^1 P^p \equiv
Q_2, \dots\} \,.
$$

\begin{lem}
\label{l7.6}
{\rm (a)}\qua
For $p=2$ the map $\psi^3-1 \: ku_{(2)} \to \Sigma^2 ku_{(2)}$ induces
right multiplication by $Sq^2$ on mod~$2$ cohomology:
$$
(\psi^3-1)^* = Sq^2 \: \Sigma^2 A/\!/E_1 \to A/\!/E_1 \,.
$$

{\rm (b)}\qua
For $p$ odd the map $\psi^r-1 \: \ell \to \Sigma^q \ell$ induces right
multiplication by $P^1$ on mod~$p$ cohomology:
$$
(\psi^r-1)^* = P^1 \: \Sigma^q A/\!/E_1 \to A/\!/E_1 \,.
$$

{\rm(c)}\qua
For $p=2$ the map $\psi^3-1 \: ko_{(2)} \to bspin_{(2)}$ induces
right multiplication by $Sq^4$ on mod~$2$ cohomology:
$$
(\psi^3-1)^* = Sq^4 \: \Sigma^4 A/A\{Sq^1, Sq^2Sq^3\} \to A/\!/A_1 \,.
$$
\end{lem}

Case~(c) is due to Mahowald and Milgram \cite[3.4]{MaMi76}.

\begin{proof}
The $S$-algebra unit map $e \: S_{(p)} \to ju$ is well-known to be
$2$-connected for $p=2$ and $(pq-2)$-connected for $p$ odd, since this
is the degree of the first element $\beta_1$ in the $p$-primary
cokernel of $J$ \cite[1.1.14]{Ra04}.  Hence $e^* \: H^*(ju; \F_p) \to
H^*(S; \F_p) = \F_p$ is cohomologically $2$-connected for $p=2$ and
cohomologically $(pq-2)$-connected for $p$ odd, meaning that the
homomorphism is injective in the stated degree and an isomorphism in
lower degrees.  In particular, $H^2(ju; \F_2) = 0$ for $p=2$ and
$H^q(ju; \F_p) = 0$ for $p$ odd.

So in the long exact cohomology sequence
\begin{equation}
\xymatrix@C+7pt{
\Sigma^2 A/\!/E_1 \ar[r]^-{(\psi^3-1)^*} & A/\!/E_1 \ar[r]^-{\kappa^*}
	& H^*(ju; \F_2)
}
\label{e7.7}
\end{equation}
associated to the cofiber sequence~(\ref{e7.2}), the non-zero class of
$Sq^2$ in $A/\!/E_1$ maps to zero under $\kappa^*$, hence is in the
image of $(\psi^3-1)^*$.  The latter is a left $A$-module homomorphism,
and can only take $\Sigma^2(1)$ to $Sq^2$, hence is given by right
multiplication by $Sq^2$.  This proves~(a).  For part~(b) we use the
same argument for the exact sequence
\begin{equation}
\xymatrix@C+7pt{
\Sigma^q A/\!/E_1 \ar[r]^-{(\psi^r-1)^*} & A/\!/E_1 \ar[r]^-{\kappa^*}
	& H^*(ju; \F_p)
}
\label{e7.8}
\end{equation}
associated to~(\ref{e7.3}), in cohomological degree~$q$.  The non-zero
class of $P^1$ in $A/\!/E_1$ maps to zero under $\kappa^*$, hence must
equal the image of $\Sigma^q(1)$ under $(\psi^r-1)^*$.  This
proves~(b).

The unit map $e \: S_{(2)} \to j$ is likewise well-known to be
$6$-connected, since this is the degree of the first element $\nu^2$ in
the $2$-primary cokernel of $J$, so $e^* \: H^*(j; \F_2) \to H^*(S;
\F_2) = \F_2$ is cohomologically $6$-connected.  In particular, $H^4(j;
\F_2) = 0$.  So in the long exact cohomology sequence associated to the
cofiber sequence~(\ref{e7.4})
\begin{equation}
\xymatrix@C+7pt{
\Sigma^4 A/A\{Sq^1, Sq^2Sq^3\} \ar[r]^-{(\psi^3-1)^*}
	& A/\!/A_1 \ar[r]^-{\kappa^*} & H^*(j; \F_2)
}
\label{e7.9}
\end{equation}
the non-zero class of $Sq^4$ in $A/\!/A_1$ maps to zero under $\kappa^*$,
hence is in the image of $(\psi^3-1)^*$.  The only class that can hit
it is $\Sigma^4(1)$, which proves~(c).
\end{proof}

\begin{lem}
\label{l7.10}
{\rm (a)}\qua
For $p=2$ there is a uniquely split extension of $A$-modules
$$
0 \to A/\!/A_1 \to H^*(ju; \F_2) \to \Sigma^3 A/\!/A_1 \to 0 \,.
$$
Hence there is a canonical $A$-module isomorphism $H^*(ju; \F_2) \cong
A/\!/A_1\{1, x\}$, with $x$ a class in degree~$3$.

{\rm (b)}\qua
For $p$ odd there is a non-split extension of $A$-modules
$$
0 \to A/\!/A_1 \to H^*(ju; \F_p) \to \Sigma^{pq-1} A/\!/A_1 \to 0 \,.
$$
As an $A$-module, $H^*(ju; \F_p)$ is generated by two classes
$1$ and $x$ in degrees~$0$ and $(pq-1)$, respectively, with
$\beta(x) = P^p(1)$.

{\rm(c)}\qua
There is a unique non-split extension of $A$-modules
$$
0 \to A/\!/A_2 \to H^*(j; \F_2) \to A \otimes_{A_2} \Sigma^7 K \to 0
\,.
$$
The cyclic $A_2$-module $K = A_2/A_2\{Sq^1, Sq^7, Sq^4Sq^6+Sq^6Sq^4\}$
has rank~$17$ over~$\F_2$.  As an $A$-module, $H^*(j; \F_2)$ is generated
by two classes~$1$ and~$x$ in degrees~$0$ and~$7$, respectively, with
$Sq^1(x) = Sq^8(1)$.
\end{lem}

For case~(b), see also \cite[5.1(b)]{Rog03}.  Case~(c) is due to Davis
\cite[Thm.~1]{Da75}, who also shows that $H^*(j; \F_2)$ is a free
$A/\!/A_3$-module.

\begin{proof}
{\rm (a)}\qua
In the long exact sequence extending~(\ref{e7.7}) the $A$-module
homomorphism $(\psi^3-1)^*$ is induced up from the $A_1$-module
homomorphism
$$
Sq^2 \: \Sigma^2 A_1/\!/E_1 \to A_1/\!/E_1 = \F_2\{1, Sq^2\}
$$
with kernel $\Sigma^2 \F_2\{Sq^2\}$ and cokernel $\F_2\{1\}$.  Since $A$
is flat (in fact free) over $A_1$, it follows that $\ker(\psi^3-1)^*
\cong \Sigma^4 A/\!/A_1$ and $\cok(\psi^3-1)^* \cong A/\!/A_1$.  Hence there
is an extension of $A$-modules
$$
0 \to A/\!/A_1 \to H^*(ju; \F_2) \to \Sigma^3 A/\!/A_1 \to 0 \,,
$$
as asserted.  The group of such extensions is trivial,
by the change-of-rings isomorphism
$$
\Ext^1_A(\Sigma^3 A/\!/A_1, A/\!/A_1) \cong \Ext^1_{A_1}(\Sigma^3 \F_2,
A/\!/A_1) \,.
$$
For in any $A_1$-module extension
$$
0 \to A/\!/A_1 \to E \to \Sigma^3 \F_2 \to 0
$$
let $x \in E$ be the unique class in degree~$3$ that maps to
$\Sigma^3(1)$.  Then $Sq^2 x = 0$ since $A/\!/A_1 = \F_2\{1, Sq^4, Sq^6,
Sq^7, \dots\}$ is trivial in degree~$5$.  Furthermore $Sq^1 x = Sq^4(1)$
would contradict the Adem relation $Sq^2 Sq^2 = Sq^3 Sq^1$, since $Sq^3
Sq^4 \equiv Sq^7$ in $A/\!/A_1$.  So $Sq^1 x = 0$ and the extension $E$
is trivial.

Two choices of splitting maps for the trivial extension describing
$H^*(ju; \F_2)$ differ by an $A$-module homomorphism
$\Sigma^3 A/\!/A_1 \to A/\!/A_1$, which must be zero since $A/\!/A_1$ is
trivial in degree~$3$.  Therefore the splitting is unique, as
claimed.

{\rm (b)}\qua
Similarly, in~(\ref{e7.8}) the $A$-module homomorphism $(\psi^r-1)^*$ is
induced up from the $A_1$-module homomorphism
$$
P^1 \: \Sigma^q A_1/\!/E_1 \to A_1/\!/E_1 = \F_p\{1, P^1, \dots, P^{p-1}\}
$$
with kernel $\Sigma^q \F_p\{P^{p-1}\}$ and cokernel $\F_p\{1\}$.
As above it follows that $\ker(\psi^r-1)^* \cong
\Sigma^{pq} A/\!/A_1$ and $\cok(\psi^r-1)^* \cong A/\!/A_1$.
Hence there is an extension of $A$-modules
$$
0 \to A/\!/A_1 \to H^*(ju; \F_p) \to \Sigma^{pq-1} A/\!/A_1 \to 0 \,,
$$
as asserted.  This time the group of extensions is non-trivial;
in fact it is isomorphic to $\Z/p$ and generated by the extension
above.  To see this, we again use the change-of-rings isomorphism
$$
\Ext^1_A(\Sigma^{pq-1} A/\!/A_1, A/\!/A_1) \cong
\Ext^1_{A_1}(\Sigma^{pq-1} \F_p, A/\!/A_1)
$$
and consider $A_1$-module extensions
$$
0 \to A/\!/A_1 \to E \to \Sigma^{pq-1} \F_p \to 0 \,.
$$
Let $x \in E$ be the unique class that maps to $\Sigma^{pq-1}(1)$.
Then $P^1 x = 0$ since $A/\!/A_1 = \F_p\{1, P^p, P^1 P^p, \dots\}$ is
trivial in degree $(p+1)q-1$.  But $\beta x$ is a multiple of $P^p$,
and this multiple in $\Z/p$ classifies the extension.

To see that $\beta x$ is non-zero in the case of $H^*(ju; \F_p)$,
recall again that the first class $\beta_1$ in the cokernel of $J$ is
in degree $(pq-2)$ and has order~$p$.  Let $c \to S_{(p)} \to ju$ be the
usual cofiber sequence.  Then by the Hurewicz and universal coefficient
theorems, the lowest class $x$ in $H^*(\Sigma c; \F_p)$ sits in degree
$(pq-1)$ and supports a non-trivial mod~$p$ Bockstein $\beta x \ne 0$.
Furthermore, $H^*(\Sigma c; \F_p) \cong H^*(ju; \F_p)$ in positive
degrees $* > 0$, so also in $H^*(ju; \F_p)$ we have $\beta x \ne 0$.

{\rm(c)}\qua
In~(\ref{e7.9}), the $A$-module homomorphism $(\psi^3-1)^*$ is induced up
from the $A_2$-module homomorphism
$$
Sq^4 \: \Sigma^4 A_2/A_2\{Sq^1, Sq^2 Sq^3\} \to A_2/\!/A_1 \,.
$$
A direct calculation shows that $A_2/A_2\{Sq^1, Sq^2Sq^3\}$ has
rank~$24$ and $A_2/\!/A_1$ has rank~$8$, as $\F_2$-vector spaces.  The
cokernel of the homomorphism $Sq^4$ is $A_2/\!/A_2 = \F_2\{1\}$, of
rank~$1$, so its kernel $\Sigma^8 K$ has rank~$17$.  Here
\begin{multline*}
\Sigma^4 K = \F_2\{Sq^4,\ Sq^6,\ Sq^7,\ Sq^6Sq^2,\ Sq^9,\
	Sq^{10}+Sq^8Sq^2,\ Sq^7Sq^3, \\
Sq^{11}+Sq^9Sq^2,\ Sq^{10}Sq^2,\ Sq^{13}+Sq^{10}Sq^3,\ Sq^{11}Sq^2,\
	Sq^{11}Sq^3, \\
Sq^{13}Sq^2+Sq^{12}Sq^3,\ Sq^{13}Sq^3,\ Sq^{17}+Sq^{15}Sq^2,\
	Sq^{17}Sq^2+Sq^{16}Sq^3,\ Sq^{17}Sq^3\}
\end{multline*}
as a submodule of $A_2/A_2\{Sq^1, Sq^2 Sq^3\}$.  By another direct
calculation, $\Sigma^4 K$ is in fact the cyclic $A_2$-submodule
generated by $Sq^4$.  The annihilator ideal turns out to be generated
by $Sq^1$, $Sq^7$ and $Sq^4Sq^6 + Sq^6Sq^4 = Sq^{10}+Sq^8Sq^2+Sq^7Sq^3$
(in admissible form), so
$$
K \cong A_2/A_2\{Sq^1, Sq^7, Sq^4Sq^6+Sq^6Sq^4\} \,.
$$
Hence there is an extension of $A$-modules
$$
0 \to A/\!/A_2 \to H^*(j; \F_2) \to A \otimes_{A_2} \Sigma^7 K
\to 0
$$
with $A \otimes_{A_2} \Sigma^7 K \cong \Sigma^7 A/A\{Sq^1, Sq^7,
Sq^4Sq^6+Sq^6Sq^4\}$.  The group of such $A$-module extensions is
$$
\Ext^1_A(A \otimes_{A_2} \Sigma^7 K, A/\!/A_2) \cong \Z/2 \,,
$$
and the extension is determined by the action of $Sq^1$ on the generator
$x$ in degree~$7$ that maps to $\Sigma^7(1)$.  ($Sq^7 x = 0$ by the Adem
relation $Sq^1 Sq^7 = 0$ and the fact that $Sq^1$ acts injectively from
degree~$14$ of $A/\!/A_2$.  $(Sq^4Sq^6+Sq^6Sq^4)x = 0$ since $A/\!/A_2$
is trivial in degree~$17$.)

To see that $Sq^1(x) = Sq^8(1) \ne 0$ in $H^*(j; \Z/2)$, we once again
use the cofiber sequence $c \to S_{(2)} \to j$ and the fact that $c$ is
$5$-connected with $\pi_6(c) = \Z/2\{\nu^2\}$.  Hence the lowest class
$x$ in $H^*(\Sigma c; \F_2)$ sits in degree~$7$ and supports a
non-trivial $Sq^1 x \ne 0$.  Again, $H^*(\Sigma c; \F_2) \cong H^*(j;
\F_2)$ in positive degrees, so also in $H^*(j; \F_2)$ we have $Sq^1 x
\ne 0$.  The only possible nonzero value is the $Sq^8$ from
$A/\!/A_2$.
\end{proof}

We display the $A_2$-module $\Sigma^4 K \subset A_2/A_2\{Sq^1,
Sq^2Sq^3\}$ below.  Here $(i)$ or $(i,j)$ denotes an admissible class
with lexicographically leading term $Sq^i$ or $Sq^i Sq^j$,
respectively.  The arrows indicate the $Sq^1$- and $Sq^2$-operations.
The $Sq^4$-operations can be deduced from the relations $Sq^4 (6) =
(10)$ and $Sq^4 (13) = (17)$, but are omitted to avoid cluttering the
diagram.

$$
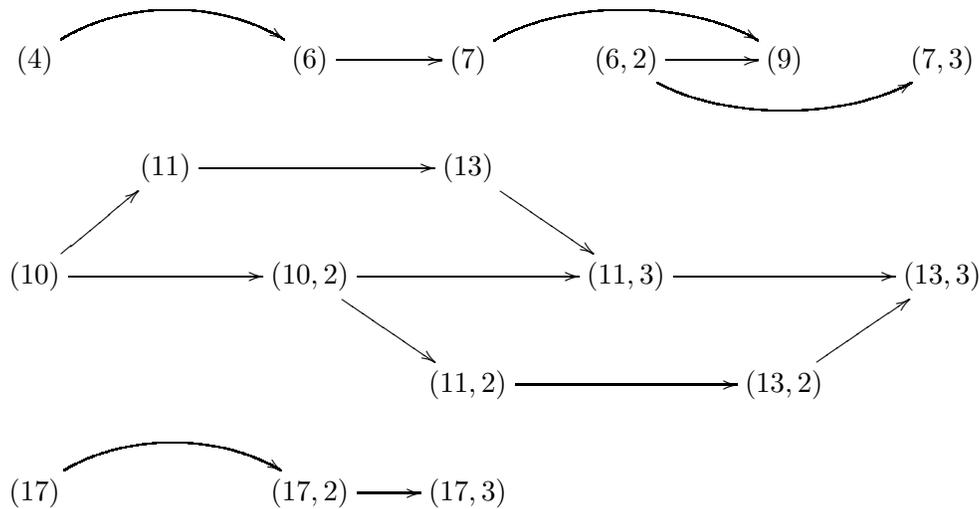
\xymatrix{
(4) \ar@(ur,ul)[rr] & & (6) \ar[r] &
	(7) \ar@(ur,ul)[rr] & (6,2) \ar[r] \ar@(dr,dl)[rr] &
	(9) & (7,3) \\
& (11) \ar[rr] & & (13) \ar[dr] \\
(10) \ar[ur] \ar[rr] & & (10,2) \ar[dr] \ar[rr]
	& & (11,3) \ar[rr] & & (13,3) \\
& & & (11,2) \ar[rr] & & (13,2) \ar[ur] \\
(17) \ar@(ur,ul)[rr] & &
	(17,2) \ar[r] & (17,3)
}
$$

\begin{sshead}
\label{s7.11}
Homology algebras
\end{sshead}
Let us write $(A/\!/A_1)_* \subset A_*$ for the $A_*$-comodule subalgebra
dual to the quotient $A$-module coalgebra $A/\!/A_1$ of $A$.  For $p=2$
we recall from Proposition~\ref{p6.1} that
$$
(A/\!/A_1)_* = P(\bar\xi_1^4, \bar\xi_2^2, \bar\xi_k \mid k\ge3) \cong
H_*(ko; \F_2) \,.
$$
For $p$ odd
$$
(A/\!/A_1)_* = P(\bar\xi_1^p, \bar\xi_k \mid k \ge 2) \otimes E(\bar\tau_k
\mid k \ge 2) \,,
$$
but there is no spectrum with mod~$p$ homology realizing
$(A/\!/A_1)_*$.  For $p=2$ and for $p$ odd there are extensions of
$A_*$-comodules
$$
\xymatrix{
0 \to \Sigma^3 (A/\!/A_1)_* \to H_*(ju; \F_2) \ar[r]^-{\kappa}
        & (A/\!/A_1)_* \to 0
}
$$
and
$$
\xymatrix{
0 \to \Sigma^{pq-1}(A/\!/A_1)_* \to H_*(ju; \F_p) \ar[r]^-{\kappa}
	& (A/\!/A_1)_* \to 0 \,,
}
$$
dual to the $A$-module extensions of Lemma~\ref{l7.10}(a) and~(b),
respectively, where in both cases $\kappa$ is an $A_*$-comodule algebra
homomorphism.

Likewise, we write $(A/\!/A_2)_* \subset A_*$ for the $A_*$-comodule
subalgebra dual to the quotient $A$-module coalgebra $A/\!/A_2$ of $A$.
For $p=2$ we recall from Proposition~\ref{p6.1} that
$$
(A/\!/A_2)_* = P(\bar\xi_1^8, \bar\xi_2^4, \bar\xi_3^2, \bar\xi_k
\mid k\ge4) \cong H_*(\tmf; \F_2) \,.
$$
There is an extension of $A_*$-comodules
$$
\xymatrix{
0 \to A_* \cotensor_{A_{2*}} \Sigma^7 K_* \to H_*(j; \F_2) \ar[r]^-{\kappa} &
(A/\!/A_2)_* \to 0 \,,
}
$$
dual to the $A$-module extension of Lemma~\ref{l7.10}(c), where
$\kappa$ is an $A_*$-comodule algebra homomorphism.  Here $K_* \subset
A_{2*}$ is the $A_{2*}$-comodule dual to the cyclic $A_2$-module $K =
A_2/A_2\{Sq^1, Sq^7, Sq^4Sq^6+Sq^6Sq^4\}$.

\begin{prop}
\label{p7.12}
{\rm (a)}\qua
For $p=2$, let $b \in H_3(ju; \F_2)$ be the image of $\Sigma^3(1)$
in $\Sigma^3 (A/\!/A_1)_*$.  Then there is an $A_*$-comodule algebra
isomorphism
$$
H_*(ju; \F_2) \cong (A/\!/A_1)_* \otimes E(b) \,,
$$
where $(A/\!/A_1)_*$ has the subalgebra structure from $A_*$ and $b$
is $A_*$-comodule primitive.

The Dyer--Lashof operations in $H_*(ju; \F_2)$ satisfy $Q^4(b) = 0$,
$Q^5(\bar\xi_1^4) = 0$, $Q^7(\bar\xi_2^2) = 0$ and $Q^{2^k}(\bar\xi_k)
= \bar\xi_{k+1}$ for all $k\ge3$, so $H_*(ju; \F_2)$ is generated by
$\bar\xi_1^4$, $\bar\xi_2^2$, $\bar\xi_3$ and $b$ as an algebra over
the Dyer--Lashof algebra.

{\rm (b)}\qua
For $p$ odd, let $b \in H_{pq-1}(ju; \F_p)$ be the image of
$\Sigma^{pq-1}\!(1)$ in $\Sigma^{pq-1}\! (A/\!/\!A_1)_*$.  There is an algebra
isomorphism
$$
(A/\!/A_1)_* \otimes E(b) \cong H_*(ju; \F_p)
$$
that takes the algebra generators $\bar\xi_1^p$, $\bar\xi_k$, $\bar\tau_k$
(for $k\ge 2$) and~$b$ to classes $\tilde\xi_1^p$, $\tilde\xi_k$,
$\tilde\tau_k$ and~$b$ in $H_*(ju; \F_p)$, respectively.  These map
under $\kappa$ to $\bar\xi_1^p$, $\bar\xi_k$, $\bar\tau_k$ and $0$,
respectively.

The Dyer--Lashof operations on $H_*(ju; \F_p)$ satisfy $Q^{pq/2}(b)
= 0$ and $Q^{p^k}(\tilde\tau_k) = \tilde\tau_{k+1}$ for all $k\ge2$,
and $\beta(\tilde\tau_k) = \tilde\xi_k$ for all $k\ge2$.  Thus $H_*(ju;
\F_p)$ is generated as an algebra over the Dyer--Lashof algebra by
$\tilde\xi_1^p$, $\tilde\xi_2$, $\tilde\tau_2$ and $b$.

The $A_*$-comodule structure is determined by
\begin{align*}
\nu(b) &= 1 \otimes b \\
\nu(\tilde\xi_1^p) &= 1 \otimes \tilde\xi_1^p - \tau_0 \otimes b +
\bar\xi_1^p \otimes 1 \\
\nu(\tilde\xi_2) &= 1 \otimes \tilde\xi_2 + \bar\xi_1 \otimes \tilde\xi_1^p
+ \tau_1 \otimes b + \bar\xi_2 \otimes 1 \\
\nu(\tilde\tau_2) &= 1 \otimes \tilde\tau_2 + \bar\tau_0 \otimes
\tilde\xi_2 + \bar\tau_1 \otimes \tilde\xi_1^p - \tau_0 \tau_1 \otimes
b + \bar\tau_2 \otimes 1 \,.
\end{align*}
(The class $b \in H_{pq-1}(ju; \F_p)$ maps under the connecting map for
the cofiber sequence $c \to S_{(p)} \to ju$ to the mod~$p$ Hurewicz
image of $\beta_1 \in \pi_{pq-2}(c)$, so the letter $b$ is chosen to
correspond to $\beta$.)

{\rm(c)}\qua
There is a square-zero extension of $A_*$-comodule algebras
$$
\xymatrix{
0 \to A_* \cotensor_{A_{2*}} \Sigma^7 K_* \to H_*(j; \F_2)
\ar[r]^-{\kappa} & (A/\!/A_2)_* \to 0 \,,
}
$$
where $\kappa$ is split as an algebra homomorphism.  As an
$(A/\!/A_2)_*$-module,
$$
\ker(\kappa) = A_* \cotensor_{A_{2*}} \Sigma^7 K_* \cong (A/\!/A_2)_*
\otimes \Sigma^7 K_*
$$
is free of rank~$17$.  There is an algebra isomorphism
$$
H_*(j; \F_2) \cong (A/\!/A_2)_* \otimes (\F_2 \oplus \Sigma^7 K_*)
$$
where $\F_2 \oplus \Sigma^7 K_*$ is the split square-zero extension of
$\F_2$ with kernel $\Sigma^7 K_*$ of rank~$17$.
\end{prop}

\begin{proof}
{\rm (a)}\qua
Let $x \in H^3(ju; \F_2)$ be the class that maps to $\Sigma^3(1)$ in the
uniquely split $A$-module extension of Lemma~\ref{l7.10}(a).  Then $\psi(x)
= x \otimes 1 + 1 \otimes x$ since $H^*(ju; \F_2) = 0$ for $0 < *
< 3$, so $E(x) = \F_2\{1, x\}$ (no algebra structure is implied) is
a sub-coalgebra of $H^*(ju; \F_2)$ and there a surjective composite
$A$-module coalgebra homomorphism
$$
A \otimes E(x) \to A \otimes H^*(ju; \F_2) \to H^*(ju; \F_2) \,.
$$
Since the $A$-module extension is split, the generators $Sq^1$ and
$Sq^2$ of $A_1$ act trivially on $1$ and $x$ in $H^*(ju; \F_2)$, so the
surjection factors through a surjection
$$
A/\!/A_1 \otimes E(x) \to H^*(ju; \F_2) \,,
$$
which by a dimension count must be an $A$-module coalgebra isomorphism.

Dually, let $b \in H_3(ju; \F_2)$ be the image of $\Sigma^3(1)$ in the
split $A_*$-comodule extension
$$
\xymatrix{
0 \to \Sigma^3 (A/\!/A_1)_* \to H_*(ju; \F_2) \ar[r]^-{\kappa}
        & (A/\!/A_1)_* \to 0 \,.
}
$$
Then $b$ is dual to $x$, and the dual of the above isomorphism is
an $A_*$-comodule algebra isomorphism
$$
H_*(ju; \F_2) \cong (A/\!/A_1)_* \otimes E(b) \,.
$$
In particular, the unique $A_*$-comodule splitting $(A/\!/A_1)_* \to H_*(ju;
\F_2)$ is an algebra map, $b$ is an $A_*$-comodule primitive, and $b^2=0$.

To determine the Dyer--Lashof operations that we will need in $H_*(ju;
\F_2)$ we make use of the Nishida relations and the known
$A_*$-comodule structure.  Some of the Nishida relations that we shall
use are
$$
Sq^1_* Q^s = \begin{cases} Q^{s-1} & \text{for $s$ even} \\
0 & \text{for $s$ odd,} \end{cases}
$$
and
$$
Sq^2_* Q^s = \begin{cases} Q^{s-2} + Q^{s-1} Sq^1_*
	& \text{for $s \equiv 0,1 \mod 4$} \\
Q^{s-1} Sq^1_* & \text{for $s \equiv 2, 3 \mod 4$.} \end{cases}
$$
See \cite[III.1.1(8)]{BMMS86}.

The Dyer--Lashof operation $Q^4(b)$ lands in $H_7(ju; \F_2) \cong
\F_2\{\bar\xi_3, \bar\xi_1^4 b\}$.  From the $A_*$-comodule structure
we can read off the dual Steenrod operations $Sq^1_*(\bar\xi_3) =
\bar\xi_2^2$ and $Sq^4_*(\bar\xi_1^4 b) = b$, since $Sq^i$ is dual to
$\xi_1^i$.  These are linearly independent, so $Q^4(b)$ is determined
by its images under $Sq^1_*$ and $Sq^4_*$.  By a Nishida relation we
get that $Sq^1_* Q^4(b) = Q^3(b)$, and $Q^3(b) = b^2 = 0$ since $|b| =
3$ and $b$ is an exterior class.  By another Nishida relation $Sq^4_*
Q^4(b) = Q^2 Sq^2_*(b)$, and $Sq^2_*(b) = 0$ since $b$ is $A_*$-comodule
primitive.  Thus $Sq^1_* Q^4(b) = 0$ and $Sq^4_* Q^4(b) = 0$, and the
only possibility is that $Q^4(b) = 0$.

The operation $Q^5(\bar\xi_1^4)$ lands in $H_9(ju; \F_2) \cong
\F_2\{\bar\xi_2^2 b\}$.  By a Nishida relation $Sq^2_* Q^5(\bar\xi_1^4)
= (Q^3 + Q^4 Sq^1_*)(\bar\xi_1^4) = 0$, while $Sq^2_*(\bar\xi_2^2 b) =
\bar\xi_1^4 b \ne 0$.  So $Q^5(\bar\xi_1^4)$ must be zero.

The operation $Q^7(\bar\xi_2^2)$ lands in $H_{13}(ju; \F_2) \cong
\F_2\{\bar\xi_2^2 \bar\xi_3, \bar\xi_1^4 \bar\xi_2^2 b\}$.  Its image
under $\kappa$ in $H_*(ku; \F_2) \subset A_*$ must vanish, by the
Cartan formula, so in fact $Q^7(\bar\xi_2^2) \in \F_2\{\bar\xi_1^4
\bar\xi_2^2 b\}$.  By a Nishida relation $Sq^2_* Q^7(\bar\xi_2^2) =
Q^6 Sq^1_*(\bar\xi_2^2) = 0$, while $Sq^2_*(\bar\xi_1^4 \bar\xi_2^2 b) =
(\bar\xi_1^4)^2 b \ne 0$, so $Q^7(\bar\xi_2^2)$ must be zero.

The claim that $Q^{2^k}(\bar\xi_k) = \bar\xi_{k+1}$ follows in the same
way as in $A_* = H_*(H\F_2; \F_2)$, see \cite[\S III.6]{BMMS86}.  It
suffices to show that $Sq^{2^m}_* Q^{2^k}(\bar\xi_k) =
Sq^{2^m}_*(\bar\xi_{k+1})$ for $0 \le m \le k$, because the only
$A_*$-comodule primitives in $H_*(ju; \F_2)$ are $1$ and~$b$.  The
right hand side equals $\bar\xi_k^2$ for $m=0$, and is zero otherwise,
by the formula for $\psi(\bar\xi_{k+1})$.  The Nishida relations imply
that $Sq^1_* Q^{2^k}(\bar\xi_k) = Q^{2^k-1}(\bar\xi_k) = \bar\xi_k^2$,
while $Sq^2_* Q^{2^k}(\bar\xi_k) = Q^{2^k-1} Sq^1_*(\bar\xi_k) =
Q^{2^k-1}(\bar\xi_{k-1}^2) = 0$ by the Cartan formula.  For $2 \le m
\le k$ we have $Sq^{2^m}_* Q^{2^k}(\bar\xi_k) = Q^{2^k-2^{m-1}}
Sq^{2^{m-1}}_*(\bar\xi_k) = Q^{2^k-2^{m-1}}(0) = 0$ by the formula for
$\psi(\bar\xi_k)$ in $A_*$.  Hence $Q^{2^k}(\bar\xi_k) = \bar\xi_{k+1}$
in $H_*(ju; \F_2)$ is the only possibility.

{\rm (b)}\qua
Let $x \in H^{pq-1}(ju; \F_p)$ be the class that maps to
$\Sigma^{pq-1}(1)$ in the $A$-module extension of Lemma~\ref{l7.10}(b).
Then $\psi(x) = x \otimes 1 + 1 \otimes x$ since $H^*(ju; \F_p) = 0$ for
$0 < * < pq-1$, so $E(x) = \F_p\{1, x\}$ is a sub-coalgebra of $H^*(ju;
\F_p)$ and there a surjective composite $A$-module coalgebra homomorphism
$$
A \otimes E(x) \to A \otimes H^*(ju; \F_p) \to H^*(ju; \F_p) \,.
$$
Dually, let $b \in H_{pq-1}(ju; \F_p)$ be the image of $\Sigma^{pq-1}(1)$
in the $A_*$-comodule extension
$$
\xymatrix{
0 \to \Sigma^{pq-1}(A/\!/A_1)_* \to H_*(ju; \F_p) \ar[r]^-{\kappa}
        & (A/\!/A_1)_* \to 0 \,.
}
$$
Then $E(b)$ is a quotient algebra of $H_*(ju; \F_p)$ and the dual of
the surjection above is an injective $A_*$-comodule algebra homomorphism
$$
H_*(ju; \F_p) \to A_* \otimes E(b) \,.
$$
We shall describe $H_*(ju; \F_p)$ in terms of its image under this injection.

Since $(A/\!/A_1)_*$ is a free graded commutative algebra, we can
choose an algebra section $s \: (A/\!/A_1)_* \to H_*(ju; \F_p)$ to the
surjection $\kappa \: H_*(ju; \F_p) \to (A/\!/A_1)_*$.  We write
$\tilde\xi_1^p = s(\bar\xi_1^p)$, $\tilde\xi_k = s(\bar\xi_k)$ and
$\tilde\tau_k = s(\bar\tau_k)$ for the lifted classes in $H_*(ju;
\F_p)$.  Since $\Sigma^{pq-1}(A/\!/A_1)$ vanishes in the degrees of
$\bar\xi_1^p$ and $\bar\tau_2$, the respective lifts $\tilde\xi_1^p$
and $\tilde\tau_2$ are unique, and we can use the commutative
$S$-algebra structure on $ju$ and the resulting Dyer--Lashof operations
on $H_*(ju; \F_p)$ to fix the other lifts by the formulas $\tilde\xi_k
= \beta(\tilde\tau_k)$ and $\tilde\tau_{k+1} = Q^{p^k}(\tilde\tau_k)$
for $k\ge2$.  This specifies the algebra section $s$ uniquely.  (But
beware, $s$ is not an $A_*$-comodule homomorphism!)

Since $|b| = pq-1$ is odd we have $b^2=0$ and the exterior algebra
$E(b)$ is a subalgebra of $H_*(ju; \F_p)$.  Writing $i \: E(b)
\to H_*(ju; \F_p)$ for the inclusion, we obtain an algebra map
$$
s \otimes i \: (A/\!/A_1)_* \otimes E(b) \to H_*(ju; \F_p) \,,
$$
which we claim is an isomorphism.  By a dimension count it suffices
to show that its composite with the injection $H_*(ju; \F_p)
\to A_* \otimes E(b)$ is injective.  We have a diagram of algebra maps
$$
\xymatrix{
(A/\!/A_1)_* \otimes E(b) \ar[r]^-{s \otimes i} \ar[d] &
H_*(ju; \F_p) \ar[r] \ar[d]^-\kappa & A_* \otimes E(b) \ar[d] \\
(A/\!/A_1)_* \ar[r]^-{=} & (A/\!/A_1)_* \ar[r] & A_*
}
$$
where the vertical maps take $b$ to zero.  The lower map takes $\xi \in
(A/\!/A_1)_*$ to $\xi \in A_*$, so the upper composite takes $\xi \otimes
1$ to $\xi \otimes 1 \pmod b$.  Hence the latter takes $\xi \otimes b$
to $\xi \otimes b \pmod {b^2 = 0}$, and it follows that the upper composite
indeed is injective.

In low degrees, $H^*(ju; \F_p) \cong \F_p\{1, x, P^p(1), P^1 P^p(1),
\beta P^1 P^p(1), \dots\}$ is dual to $H_*(ju; \F_p) \cong \F_p\{1, b,
- \tilde\xi_1^p, \tilde\xi_2, \tilde\tau_2, \dots\}$.  By
Lemma~\ref{l7.10}(b) we have $\beta(x) = P^p(1)$, so dually the
$A_*$-comodule coactions are given by
\begin{align*}
\nu(b) &= 1 \otimes b \\
\nu(\tilde\xi_1^p) &= 1 \otimes \tilde\xi_1^p - \tau_0 \otimes b +
\bar\xi_1^p \otimes 1 \\
\nu(\tilde\xi_2) &= 1 \otimes \tilde\xi_2 + \bar\xi_1 \otimes \tilde\xi_1^p
+ \tau_1 \otimes b + \bar\xi_2 \otimes 1 \\
\nu(\tilde\tau_2) &= 1 \otimes \tilde\tau_2 + \bar\tau_0 \otimes
\tilde\xi_2 + \bar\tau_1 \otimes \tilde\xi_1^p - \tau_0 \tau_1 \otimes
b + \bar\tau_2 \otimes 1 \,.
\end{align*}
In particular, the images in $A_* \otimes E(b)$ of these classes are $1
\otimes b$, $-\tau_0 \otimes b + \bar\xi_1^p \otimes 1$, $\tau_1 \otimes b
+ \bar\xi_2 \otimes 1$ and $-\tau_0\tau_1 \otimes b + \bar\tau_2 \otimes
1$, respectively.

The Dyer--Lashof operation $Q^{pq/2}(b)$ lands in $H_{p^2q-1}(ju; \F_p)
= \F_p\{\tilde\xi_1^{pq/2} b\}$.  Here $P^{pq/2}_*(\tilde\xi_1^{pq/2} b)
= b$ is non-zero, in view of the formula above for $\nu(\tilde\xi_1^p)$.
By a Nishida relation $P^{pq/2}_* Q^{pq/2}(b) = Q^{q/2} P^{q/2}_*(b) =
0$, so $Q^{pq/2}(b) = 0$.

{\rm(c)}\qua
Let $x \in H^7(j; \F_2)$ be the class that maps to $\Sigma^7(1)$ in the
$A$-module extension of Lemma~\ref{l7.10}(c).  Then $\psi(x) = x
\otimes 1 + 1 \otimes x$ since $H^*(j; \F_2) = 0$ for $0 < * < 7$, so
$E(x) = \F_2\{1, x\}$ is a sub-coalgebra of $H^*(j; \F_2)$ and there is
a surjective composite $A$-module coalgebra homomorphism
$$
A \otimes E(x) \to A \otimes H^*(j; \F_2) \to H^*(j; \F_2) \,.
$$
Dually, let $b \in H_7(j; \F_2)$ be the image of $\Sigma^7(1)$
in the $A_*$-comodule extension
$$
\xymatrix{
0 \to A_* \cotensor_{A_{2*}} \Sigma^7 K_* \to H_*(j; \F_2) \ar[r]^-{\kappa}
	& (A/\!/A_2)_* \to 0 \,.
}
$$
Here $K_* \subset A_{2*}$ has rank~$17$, and contains $1 \in K_*
\subset A_{2*}$ in degree~$0$.  The dual of the surjection above is
an injective $A_*$-comodule algebra homomorphism
$$
H_*(j; \F_2) \to A_* \otimes E(b)
$$
that may otherwise be described as the composite of the $A_*$-comodule
coaction $H_*(j; \F_2) \to A_* \otimes H_*(j; \F_2)$ and the algebra
surjection $H_*(j; \F_2) \to E(b)$.

We obtain a vertical map of $A_*$-comodule extensions
$$
\xymatrix{
A_* \cotensor_{A_{2*}} \Sigma^7 K_* \ar[r] \ar[d] &
H_*(j; \F_2) \ar[r]^-\kappa \ar[d] &
(A/\!/A_2)_* \ar[d] \\
A_*\{b\} \ar[r] & A_* \otimes E(b) \ar[r] & A_*
}
$$
where the right hand square consists of $A_*$-comodule algebra
homomorphisms, and the vertical maps are injective.  At the left hand
side we find the composite map
$$
A_* \cotensor_{A_{2*}} \Sigma^7 K_* \to A_* \cotensor_{A_{2*}} \Sigma^7 A_{2*}
\cong \Sigma^7 A_* \cong A_*\{b\}
$$
that is dual to the surjection $A\{x\} \cong \Sigma^7 A \to A
\otimes_{A_2} \Sigma^7 K$.

In the lower row the ideal $A_*\{b\}$ has square zero, since $b^2=0$,
so also the ideal $\ker(\kappa) = A_* \cotensor_{A_{2*}} \Sigma^7 K_*$ is a
square-zero ideal \cite[XIV.2]{CE56} in $H_*(j; \F_2)$.  Its module action
by $H_*(j; \F_2)$ therefore descends to one by $(A/\!/A_2)_*$.  It can be
described in terms of the algebra product $\phi$ on $A_*$ by the following
commutative diagram with injective vertical maps:
$$
\xymatrix{
(A/\!/A_2)_* \otimes (A_* \cotensor_{A_{2*}} \Sigma^7 K_*) \ar[r] \ar[d]
& (A_* \cotensor_{A_{2*}} \Sigma^7 K_*) \ar[d] \\
A_* \otimes \Sigma^7 A_* \ar[r]^-{\Sigma^7\phi} & \Sigma^7 A_*
}
$$
The proof is easy given the algebra embedding of $H_*(j; \F_2)$
into $A_* \otimes E(b)$.

In fact, the square-zero ideal $A_* \cotensor_{A_{2*}} \Sigma^7 K_*$
is a free $(A/\!/A_2)_*$-module of rank~$17$.  For $A_2 \subset A$ is a
direct summand as an $A_2$-module, so dually $A_* \to A_{2*}$ admits an
$A_{2*}$-comodule section $s \: A_{2*} \to A_*$.  For example, the image
of $s$ may be $\F_2\{\bar\xi_1^i \bar\xi_2^j \bar\xi_3^k \mid i<8, j<4,
k<2\} \subset A_*$.  We then have a map
$$
s \cotensor id \: \Sigma^7 K_* \cong A_{2*} \cotensor_{A_{2*}} \Sigma^7 K_*
\to A_* \cotensor_{A_{2*}} \Sigma^7 K_* \,.
$$
Its composite with the inclusion $A_* \cotensor_{A_{2*}} \Sigma^7 K_*
\to A_* \cotensor_{A_{2*}} \Sigma^7 A_{2*} \cong \Sigma^7 A_*$ factors as
the two inclusions $\Sigma^7 K_* \to \Sigma^7 A_{2*} \to \Sigma^7 A_*$.

Combining $s \cotensor id$ with the $(A/\!/A_2)_*$-module action on $A_*
\cotensor_{A_{2*}} \Sigma^7 K_*$ we obtain the left hand map $f$
in a commuting diagram
$$
\xymatrix{
(A/\!/A_2)_* \otimes \Sigma^7 K_* \ar[r] \ar[d]^-f &
(A/\!/A_2)_* \otimes \Sigma^7 A_{2*} \ar[d]^-\cong \\
A_* \cotensor_{A_{2*}} \Sigma^7 K_* \ar[r] &
A_* \cotensor_{A_{2*}} \Sigma^7 A_{2*}  \,.
}
$$
Here $A_* \cotensor_{A_{2*}} \Sigma^7 A_{2*} \cong \Sigma^7 A_*$
and the right hand isomorphism exhibits $\Sigma^7 A_*$ as a free
$(A/\!/A_2)_*$-module on the generators given by the section $s \: \Sigma^7
A_{2*} \to \Sigma^7 A_*$.  (It is a case of the Milnor--Moore comodule
algebra theorem \cite[4.7]{MiMo65}.)  The upper map is injective, hence
so is
$$
f \: (A/\!/A_2)_* \otimes \Sigma^7 K_*
	\to A_* \cotensor_{A_{2*}} \Sigma^7 K_* \,.
$$
But both sides have the same, finite dimension over~$\F_2$ in each degree,
so in fact $f$ is an isomorphism of $(A/\!/A_2)_*$-modules.

The fact that $\kappa$ splits as an algebra homomorphism is clear since
$(A/\!/A_2)_*$ is a free graded commutative algebra over~$\F_2$.
However, the splitting is not an $A_*$-comodule homomorphism, and
the $(A/\!/A_2)_*$-module isomorphism $f$ is not an $A_*$-comodule
isomorphism.  The $A_*$-comodule algebra structure on $H_*(j; \F_2)$
may, if desired, be obtained by describing the image of the algebra
generators under the $A_*$-comodule algebra embedding $H_*(j; \F_2)
\to A_* \otimes E(b)$.
\end{proof}

\begin{prop}
\label{p7.13}
{\rm (a)}\qua
For $p=2$ the B{\"o}kstedt spectral sequence $E^r_{**}(ju)$
collapses at the $E^2$-term, with
$$
E^\infty_{**}(ju) \cong H_*(ju; \F_2) \otimes
E(\sigma\bar\xi_1^4, \sigma\bar\xi_2^2, \sigma\bar\xi_k \mid k\ge3)
\otimes \Gamma(\sigma b) \,.
$$

{\rm (b)}\qua
For $p$ odd the B{\"o}kstedt spectral sequence $E^r_{**}(ju)$
collapses at the $E^p$-term, with
$$
E^\infty_{**}(ju) \cong H_*(ju; \F_p)
	\otimes E(\sigma\tilde\xi_1^p, \sigma\tilde\xi_2)
	\otimes P_p(\sigma\tilde\tau_k \mid k\ge2)
	\otimes \Gamma(\sigma b) \,.
$$

{\rm(c)}\qua
The B{\"o}kstedt spectral sequence for $j$ at $p=2$ has
$E^2$-term
$$
E^2_{**}(j) \cong HH_*((A/\!/A_2)_*) \otimes HH_*(\F_2 \oplus \Sigma^7 K_*)
$$
where
$$
HH_*((A/\!/A_2)_*) \cong (A/\!/A_2)_* \otimes E(\sigma\bar\xi_1^8,
\sigma\bar\xi_2^4, \sigma\bar\xi_3^2, \sigma\bar\xi_k \mid k\ge4)
$$
and
$$
HH_q(\F_2 \oplus \Sigma^7 K_*) \cong [(\Sigma^7 K_*)^{\otimes q}]^{C_q}
\oplus [(\Sigma^7 K_*)^{\otimes {q+1}}]_{C_{q+1}} \,.
$$
In particular, $E^2_{**}(j)$ is not flat as a module over $H_*(j; \F_2)$.
\end{prop}

\begin{proof}
{\rm (a)}\qua
The B{\"o}kstedt spectral sequence for $ju$ at $p=2$ begins
$$
E^2_{**}(ju) = H_*(ju; \F_2) \otimes E(\sigma\bar\xi_1^4,
\sigma\bar\xi_2^2, \sigma\bar\xi_k \mid k\ge3) \otimes \Gamma(\sigma b)
\,.
$$
Proposition~\ref{p4.8} applies, so a shortest non-zero differential
must map from an algebra indecomposable to a coalgebra primitive and
$A_*$-comodule primitive.  Here we are referring to the $A_*$-comodule
$H_*(ju; \F_2)$-Hopf algebra structure on $E^2_{**}(ju)$.  The only
possible algebra indecomposables are the $\gamma_{2^k}(\sigma b)$ in
degrees $2^{k+2}$, for $k\ge2$.  The coalgebra primitives are $H_*(ju;
\F_2) \{\sigma b, \sigma\bar\xi_1^4, \sigma\bar\xi_2^2, \sigma\bar\xi_k
\mid k\ge 3\}$, all in filtration $s=1$.  A calculation with
$\nu(\sigma\bar\xi_2^2)$ and $\nu(\sigma\bar\xi_3)$ shows that among
these, the $A_*$-comodule primitives are equal to
$$
E(b) \otimes
\F_2\{\sigma b, \sigma\bar\xi_1^4, \sigma\bar\xi_k \mid k\ge4\} \,.
$$
These live in degrees $4$, $5$, $7$, $8$, $2^k$ and $2^k+3$, for
$k\ge4$.  The image of a differential on $\gamma_{2^k}(\sigma b)$ must
be in total degree $2^{k+2}-1$, for $k\ge2$, but these degrees do not
contain any simultaneous coalgebra- and comodule primitives.  Therefore
there are no non-zero differentials, and the spectral sequence
collapses at the $E^2$-term.

{\rm (b)}\qua
For $p$ odd the spectral sequence begins
$$
E^2_{**}(ju) = H_*(ju; \F_p) \otimes
E(\sigma\tilde\xi_1^p, \sigma\tilde\xi_k \mid k\ge2) \otimes
\Gamma(\sigma b, \sigma\tilde\tau_k \mid k\ge2) \,.
$$
By Proposition~\ref{p5.6} we have $E^2 = E^{p-1}$ and there are differentials
$$
d^{p-1}(\gamma_p(\sigma\tilde\tau_k)) = \sigma\tilde\xi_{k+1}
$$
for $k\ge2$.  This uses the relation $\beta Q^{p^k}(\tilde\tau_k)
= \beta(\tilde\tau_{k+1}) = \tilde\xi_{k+1}$.  There is also a
potential differential
$$
d^{p-1}(\gamma_p(\sigma b)) = \sigma(\beta Q^{pq/2}(b)) \,,
$$
but $\beta Q^{pq/2}(b)$ is in degree $p^2q-2$ of $H_*(ju; \F_p)$,
which is a trivial group, so this differential is zero.
Hence
$$
E^p_{**}(ju) = H_*(ju; \F_p) \otimes
E(\sigma\tilde\xi_1^p, \sigma\tilde\xi_2) \otimes
P_p(\sigma\tilde\tau_k \mid k\ge2) \otimes \Gamma(\sigma b) \,.
$$
Proposition~\ref{p4.8} applies again, so a shortest differential must
map from one of the algebra indecomposables $\gamma_{p^k}(\sigma b)$ in
degrees $p^{k+1}q$, for $k\ge2$.  Its target must be among the
coalgebra primitives, which are $H_*(ju; \F_p)\{\sigma b,
\sigma\tilde\xi_1^p, \sigma\tilde\xi_2, \sigma\tilde\tau_k \mid
k\ge2\}$, all in filtration $s=1$.  The target must also be
$A_*$-comodule primitive.  The formulas for the $A_*$-comodule
structure on $H_*(ju; \F_p)$ imply the following formulas:
\begin{align*}
\nu(\sigma b) &= 1 \otimes \sigma b \\
\nu(\sigma\tilde\xi_1^p) &= 1 \otimes \sigma\tilde\xi_1^p
- \tau_0 \otimes \sigma b \\
\nu(\sigma\tilde\xi_2) &= 1 \otimes \sigma\tilde\xi_2 +
\bar\xi_1 \otimes \sigma\tilde\xi_1^p + \tau_1 \otimes \sigma b \\
\nu(\sigma\tilde\tau_2) &= 1 \otimes \sigma\tilde\tau_2 +
\bar\tau_0 \otimes \sigma\tilde\xi_2 + \bar\tau_1 \otimes
\sigma\tilde\xi_1^p - \tau_0\tau_1 \otimes \sigma b \,.
\end{align*}
The $\sigma\tilde\tau_k$ are $A_*$-comodule primitives for $k\ge3$, in
view of the relations $(\sigma\tilde\tau_k)^p = \sigma\tilde\tau_{k+1}$
for $k\ge2$, the formula for $\nu(\sigma\tilde\tau_2)$ and the fact
that $\bar\tau_0$, $\bar\tau_1$ and $\tau_0\tau_1$ all square to zero.
A calculation then shows that the simultaneous coalgebra- and comodule
primitives are equal to
$$
E(b) \otimes \F_p\{\sigma b, \sigma\tilde\tau_k \mid k\ge3\} \,.
$$
These live in degrees $pq$, $2pq-1$, $2p^k$ and $2p^k+pq-1$ for
$k\ge3$.  The image of a differential $d^r(\gamma_{p^k}(\sigma b))$ is
in degree $p^{k+1}q-1$ for $k\ge2$, which contains none of the possible
target classes.  Hence there are no further differentials, and the
spectral sequence collapses at the $E^p$-term.

{\rm(c)}\qua
By the K{\"u}nneth formula
$$
E^2_{**}(j) = HH_*(H_*(j; \F_2)) \cong
HH_*((A/\!/A_2)_*) \otimes HH_*(\F_2 \oplus \Sigma^7 K_*) \,.
$$
Here the first tensor factor was identified in the discussion of $\tmf$
in Section~6.  By the following Lemma~\ref{l7.14},
$$
HH_q(\F_2 \oplus \Sigma^7 K_*)
\cong [\Sigma^7 K_*^{\otimes q}]^{C_q}
\oplus [\Sigma^7 K_*^{\otimes (q+1)}]_{C_{q+1}} \,,
$$
and e.g.~$HH_1(\F_2 \oplus \Sigma^7 K_*)$ is not flat as an
$\F_2 \oplus \Sigma^7 K_*$-module.
\end{proof}

\begin{lem}
\label{l7.14}
Let $k$ be a field and $V$ a graded $k$-vector space.
The Hochschild homology of the split square-zero extension
$k \oplus V$, with unit $(1, 0)$ and multiplication
$(k_1, v_1) \cdot (k_2, v_2) = (k_1 k_2, k_1 v_2 + k_2 v_1)$, is
$$
HH_q(k \oplus V) \cong [V^{\otimes q}]^{C_q} \oplus
[V^{\otimes(q+1)}]_{C_{q+1}}
$$
where $[V^{\otimes q}]^{C_q} \subset V^{\otimes q}$ denotes the invariants
of the cyclic group $C_q$ of order $q$ acting by cyclic permutations on
$V^{\otimes q}$, and $[V^{\otimes q}]_{C_q}$ denotes the coinvariants
of this action.  When $\dim_k V \ge 2$ the Hochschild homology
$HH_*(k \oplus V)$ is not flat as a module over $k \oplus V$.
\end{lem}

\begin{proof}
We compute $HH_*(k \oplus V)$ as the homology of the normalized Hochschild
complex $NC_*(k \oplus V)$ with
$$
NC_q(k \oplus V) = (k \oplus V) \otimes V^{\otimes q} \cong V^{\otimes q} \oplus V^{\otimes(q+1)} \,.
$$
Here we remind the reader that the normalized Hochschild complex of a
$k$-algebra $\Lambda$ is given by $NC_q(\Lambda) = \Lambda \otimes
(\Lambda/k)^{\otimes q}$.  We are considering the case $\Lambda = k
\oplus V$, with $\Lambda/k = V$.

Since $V$ is a square-zero ideal, the Hochschild boundary $\partial$
is the direct sum over $q\ge1$ of the operators
$$
1 + (-1)^q t_q \: V^{\otimes q} \to V^{\otimes q} \,,
$$
i.e., $1+(-1)^q t_q$ on the $V^{\otimes q}$-summand and zero on the
$V^{\otimes(q+1)}$-summand, where $t_q(v_1 \otimes \dots \otimes v_q)
= (-1)^\epsilon v_q \otimes v_1 \otimes \dots \otimes v_{q-1}$ with
$\epsilon = |v_q|(|v_1| + \dots + |v_{q-1}|)$.  Let the generator $T \in
C_q$ act on $V^{\otimes q}$ as $(-1)^{q+1} t_q$, so $\partial$ is the sum
of the operators $1-T$ (and $T^q$ acts as the identity).  Clearly, then,
the Hochschild homology is the direct sum over $q\ge1$ of the kernels
$$
\ker(1-T) = [V^{\otimes q}]^{C_q}
$$
in degree $q$,
and the cokernels
$$
\cok(1-T) = [V^{\otimes q}]_{C_q}
$$
in degree $(q-1)$, plus the term $k = V^{\otimes 0}$ in degree~$0$.

For an example of the failure of flatness, let $V = k\{x, y\}$ with $x,
y$ in odd degrees and $q=1$.  Then $HH_1(k \oplus V) \cong V \oplus
V^{\otimes2}_{C_2} \cong k\{\sigma x, \sigma y, x \sigma x, x \sigma y
\equiv y \sigma x, y \sigma y\}$ is not flat over $k \oplus V$.
\end{proof}

\begin{thm}
\label{t7.15}
{\rm (a)}\qua
For $p=2$ there is an isomorphism
$$
H_*(\THH(ju); \F_2) \cong H_*(ju; \F_2) \otimes E(\sigma\bar\xi_1^4,
\sigma\bar\xi_2^2) \otimes P(\sigma\bar\xi_3) \otimes \Gamma(\sigma b)
$$
of $H_*(ju; \F_2)$-Hopf algebras.

{\rm (b)}\qua
For $p$ odd there is an isomorphism
$$
H_*(\THH(ju); \F_p) \cong H_*(ju; \F_p) \otimes E(\sigma\tilde\xi_1^p,
\sigma\tilde\xi_2) \otimes P(\sigma\tilde\tau_2) \otimes \Gamma(\sigma b)
$$
of $H_*(ju; \F_p)$-Hopf algebras.
\end{thm}

\begin{proof}
{\rm (a)}\qua
In view of Proposition~\ref{p7.13}(a) we must identify the possible algebra
extensions between the $E^\infty$-term
$$
E^\infty_{**}(ju) \cong H_*(ju; \F_2) \otimes E(\sigma\bar\xi_1^4,
\sigma\bar\xi_2^2, \sigma\bar\xi_k \mid k\ge3) \otimes \Gamma(\sigma b)
$$
and the abutment $H_*(\THH(ju); \F_2)$.  Here $H_*(ju; \F_2) \cong
(A/\!/A_1)_* \otimes E(b)$ by Proposition~\ref{p7.12}.

In $H_*(\THH(ju); \F_2)$ we have $(\sigma b)^2 = Q^4(\sigma b) = \sigma
Q^4(b) = 0$, $(\sigma\bar\xi_1^4)^2 = 0$, $(\sigma\bar\xi_2^2)^2 = 0$
and $(\sigma\bar\xi_k)^2 = \sigma\bar\xi_{k+1}$ for all $k\ge3$, by
Proposition~\ref{p5.9} and the statement about Dyer--Lashof operations
in Proposition~\ref{p7.12}(a).  It remains to prove that we can find
classes
$$
\gamma_{2^k} \in H_{4 \cdot 2^k}(\THH(ju); \F_2)
$$
that are represented by $\gamma_{2^k}(\sigma b)$ in $E^\infty_{**}(ju)$
and satisfy $\gamma_{2^k}^2 = 0$, for all $k\ge0$.  We have just seen
that we can take $\gamma_1 = \sigma b$.  So fix a number $k\ge1$, and
assume inductively that we have chosen classes $\gamma_{2^m}$ for $0
\le m < k$ that are represented by $\gamma_{2^m}(\sigma b)$ and satisfy
$\gamma_{2^m}^2 = 0$.

We shall prove below that a class $\gamma_{2^k}$ representing
$\gamma_{2^k}(\sigma b)$ can be chosen so that its square
$\gamma_{2^k}^2$ is both an $H_*(ju; \F_2)$-coalgebra primitive and an
$A_*$-comodule primitive.  We saw in the proof of
Proposition~\ref{p7.13}(a) that the simultaneous coalgebra- and
comodule primitives are
$$
E(b) \otimes
\F_2\{\sigma b, \sigma\bar\xi_1^4, \sigma\bar\xi_m \mid m\ge4\} \,.
$$
When $k\ge1$, the only such class in the degree of $\gamma_{2^k}^2$ is
$\sigma\bar\xi_{k+3} = (\sigma\bar\xi_{k+2})^2$.  So either
$\gamma_{2^k}^2 = 0$ or $\gamma_{2^k}^2 = (\sigma\bar\xi_{k+2})^2$.  In
the latter case we change $\gamma_{2^k}$ by subtracting
$\sigma\bar\xi_{k+2}$, which does not alter the representative at the
$E^\infty$-term.  Thereby we have achieved $\gamma_{2^k}^2 = 0$, which
will complete the inductive step.

To show that $\gamma_{2^k}^2$ can be arranged to be a coalgebra- and
comodule primitive, we make use of the maps of $E^\infty$-terms and
abutments induced by the commutative $S$-algebra homomorphism $\kappa \:
ju \to ku_{(2)}$.  The target $E^\infty$-term is
$$
E^\infty_{**}(ku) \cong H_*(ku; \F_2) \otimes E(\sigma\bar\xi_1^2,
\sigma\bar\xi_2^2, \sigma\bar\xi_k \mid k\ge3)
$$
with abutment
$$
H_*(\THH(ku); \F_2) \cong H_*(ku; \F_2) \otimes
E(\sigma\bar\xi_1^2, \sigma\bar\xi_2^2) \otimes P(\sigma\bar\xi_3)
\,.
$$
Here $H_*(ku; \F_2) \cong (A/\!/E_1)_*$.  Note that in degrees less than
that of $\gamma_{2^k}$, $\kappa$ maps $H_*(\THH(ju); \F_2)$ modulo classes
that square to zero, which by the inductive hypothesis is $(A/\!/A_1)_*
\otimes P(\sigma\bar\xi_3)$, injectively into $H_*(\THH(ku); \F_2)$
modulo classes that square to zero, which is $(A/\!/E_1)_* \otimes
P(\sigma\bar\xi_3)$.  We shall refer to this property as the
``near-injectivity of $\kappa$''.

So choose any class $\gamma_{2^k}$ in $H_*(\THH(ju); \F_2)$ that is
represented by $\gamma_{2^k}(\sigma b)$ in $E^\infty_{**}(ju)$.  We
shall arrange that its image $\kappa\gamma_{2^k}$ in $H_*(\THH(ku);
\F_2)$ squares to zero.  If not, we can write
$$
\kappa \gamma_{2^k} \equiv c \cdot (\sigma\bar\xi_3)^\ell
$$
with $c \in (A/\!/E_1)_*$ not equal to zero, modulo classes in
$H_*(\THH(ku); \F_2)$ that square to zero, and modulo similar terms with
$c$ of lower degree (or equivalently, with higher exponent $\ell$).  So
$c \cdot (\sigma\bar\xi_3)^\ell$ is the ``leading term'' of $\kappa
\gamma_{2^k}$ in $H_*(\THH(ku); \F_2)$ modulo classes that square to
zero.

We divide into three cases.  First, if $\ell=0$ and $\kappa\gamma_{2^k}
\equiv c$ with $c \in (A/\!/E_1)_*$, we can apply the Hopf algebra counit
$\epsilon \: \THH(R) \to R$ to see that $c = \epsilon(\kappa \gamma_{2^k})
= \kappa(\epsilon \gamma_{2^k})$ is in the image of $\kappa \: H_*(ju;
\F_2) \to H_*(ku; \F_2)$, so that in fact $c \in (A/\!/A_1)_*$.  We can
then replace $\gamma_{2^k}$ by $\gamma_{2^k} - c$ without altering its
representative in $E^\infty$, and thus eliminate the case $\ell=0$.

Second, if $\ell=1$ and $\kappa\gamma_{2^k} \equiv c \cdot
\sigma\bar\xi_3$ with $c \in (A/\!/E_1)_*$, we need to argue that $c
\in (A/\!/A_1)_*$.  If not, we can write $c = \bar\xi_1^2 \cdot d$ with
$d \in (A/\!/A_1)_*$.  Here $\nu(c) = \psi(c) \equiv c \otimes 1 + d
\otimes \bar\xi_1^2$ modulo terms with first tensor factor in degree
less than $|d|$.  Note from the structure of $H_*(\THH(ku); \F_2)$ that
$\kappa\gamma_{2^k} \equiv c \cdot \sigma\bar\xi_3$ modulo classes that
square to zero, and similar terms $e \cdot (\sigma\bar\xi_3)^\ell$ with
$\ell \ge 2$, so $|e| < |d|$.  So $\nu(\kappa \gamma_{2^k}) \equiv c
\otimes \sigma\bar\xi_3 + d \otimes \bar\xi_1^2 \sigma\bar\xi_3$ modulo
classes that square to zero, and terms with first factor in lower
degree.  This must be the image of $\nu(\gamma_{2^k})$ under $1 \otimes
\kappa$, so
$$
\nu(\gamma_{2^k}) \equiv c \otimes \sigma\bar\xi_3 + d \otimes x
$$
for some class $x \in H_*(\THH(ju); \F_2)$ with $\kappa(x) \equiv
\bar\xi_1^2 \sigma\bar\xi_3$.  But there exists no such class $x$, since
$H_*(\THH(ju); \F_2)$ modulo classes that square to zero is trivial in
degree~$10$.  (This is clear by inspection of the B{\"o}kstedt spectral
sequence, where the group that could cause a problem, $E^2_{2,8}(ju)$,
is in fact zero.)  This contradiction shows that $c \in (A/\!/A_1)_*$, and
we can alter $\gamma_{2^k}$ by $c \cdot \sigma\bar\xi_3$ without altering
the representative at $E^\infty$, and thus eliminate the case $\ell=1$.

Third, for the remaining cases $\ell\ge2$ we consider the coalgebra
coproduct.  We find
$$
\psi(\kappa\gamma_{2^k}) \equiv
\sum_{i+j=\ell} c \cdot (\sigma\bar\xi_3)^i \otimes (\sigma\bar\xi_3)^j
$$
in $H_*(\THH(ku); \F_2) \otimes_{H_*(ku; \F_2)} H_*(\THH(ku); \F_2)$, modulo
classes that square to zero and similar terms with $c$ of lower degree.
Since $\ell\ge2$ this sum includes some terms $c \cdot (\sigma\bar\xi_3)^i
\otimes (\sigma\bar\xi_3)^j$ with both $i$ and $j$ positive, so that $c
\cdot (\sigma\bar\xi_3)^i$ and $(\sigma\bar\xi_3)^j$ are both in degree
less than that of $\gamma_{2^k}$.  This is the image under $\kappa \otimes
\kappa$ of $\psi(\gamma_{2^k})$ in $H_*(\THH(ju); \F_2) \otimes_{H_*(ju;
\F_2)} H_*(\THH(ju); \F_2)$, so it follows by the near-injectivity of
$\kappa$ that $c$ in fact lies in $(A/\!/A_1)_* \subset (A/\!/E_1)_*$.

Then we can change the chosen $\gamma_{2^k}$ in $H_*(\THH(ju); \F_2)$
by subtracting $c \cdot (\sigma\bar\xi_3)^\ell$ from it, and thus
remove the ``leading'' term $c \cdot (\sigma\bar\xi_3)^\ell$ from
$\kappa\gamma_{2^k}$.  By repeating this process we can arrange that
$\gamma_{2^k}$ has be chosen so that $\kappa\gamma_{2^k}$ is zero modulo
classes that square to zero, i.e., that $\kappa\gamma_{2^k}^2 = 0$.

Then $\psi(\gamma_{2^k}) \equiv \gamma_{2^k} \otimes 1 +
1 \otimes \gamma_{2^k}$ modulo classes that square to zero.
For $\psi(\kappa\gamma_{2^k}^2) = 0$ so any other terms in
$\psi(\gamma_{2^k})$ must map under $\kappa$ to classes that square to
zero, hence square to zero themselves by the near-injectivity of $\kappa$.
Hence $\gamma_{2^k}^2$ is a coalgebra primitive.

Next consider the $A_*$-comodule coaction.  If $\nu(\gamma_{2^k})
\equiv 1 \otimes \gamma_{2^k}$ modulo classes that square to zero, then
$\nu(\gamma_{2^k}^2) = 1 \otimes \gamma_{2^k}^2$ and $\gamma_{2^k}^2$
is an $A_*$-comodule primitive, as desired.  Otherwise, we can write
$$
\nu(\gamma_{2^k}) \equiv a \otimes (\sigma\bar\xi_3)^{\ell}
$$
with $a \in A_*$ in positive degree, modulo classes in $A_* \otimes
H_*(\THH(ju); \F_2)$ that square to zero, and modulo similar terms with
$a$ of lower degree.  Then $\nu(\gamma_{2^k}^2) \equiv a^2 \otimes
(\sigma\bar\xi_3)^{2\ell}$ modulo terms with $a$ of lower degree.
Applying $\kappa$ yields $a^2 \otimes (\sigma\bar\xi_3)^{2\ell} = 0$,
since $\kappa\gamma_{2^k}^2 = 0$, so $a^2 = 0$.  This is impossible
for $a \ne 0$, so we conclude that $\gamma_{2^k}^2$ is indeed an
$A_*$-comodule primitive.

This completes the proof.

{\rm (b)}\qua
Also in the odd primary case we must identify the algebra extensions
between the $E^\infty$-term from Proposition~\ref{p7.13}(b)
$$
E^\infty_{**}(ju) = H_*(ju; \F_p) \otimes E(\sigma\tilde\xi_1^p,
\sigma\tilde\xi_2) \otimes P_p(\sigma\tilde\tau_k \mid k\ge2) \otimes
\Gamma(\sigma b)
$$
and the abutment $H_*(\THH(ju); \F_p)$.  In the latter we have $(\sigma
b)^p = Q^{pq/2}(\sigma b) = \sigma Q^{pq/2}(b) = 0$ and
$(\sigma\tilde\tau_k)^p = Q^{p^k}(\sigma\tilde\tau_k) = \sigma
Q^{p^k}(\tilde\tau_k)= \sigma\tilde\tau_{k+1}$, by
Propositions~\ref{p5.9} and~\ref{p7.12}(b).  The classes
$\sigma\tilde\xi_1^p$ and $\sigma\tilde\xi_2$ are in odd degree, and
therefore have square zero, since $p$ is odd.

It remains to prove that we can find classes $\gamma_{p^k} \in
H_{p^{k+1}q}(\THH(ju); \F_p)$ that are represented by $\gamma_{p^k}(\sigma
b)$ in $E^\infty_{**}(ju)$ and satisfy $\gamma_{p^k}^p = 0$, for all
$k\ge0$.  We have just verified this for $k=0$.

The remaining inductive proof follows exactly the same strategy as in
the $p=2$ case.  Instead of working modulo classes that square to zero we
work modulo classes with $p$-th power equal to zero.  The $S$-algebra map
$\kappa \: ju \to \ell$ induces an algebra homomorphism in homology that
has the required ``near--injectivity'' property, etc.  In the tricky
case when $\kappa \gamma_{p^k} \equiv c \cdot \sigma\tilde\tau_2$,
modulo classes in $H_*(\THH(\ell); \F_p)$ that have trivial $p$-th
power, and modulo similar terms with $c$ of lower degree, we have $c
\in (A/\!/E_1)_*$ and must argue that in fact $c \in (A/\!/A_1)_*$.
Writing $c = \bar\xi_1^f \cdot d$, with $d \in (A/\!/A_1)_*$ and $0 \le
f < p$ we assume $0 < f < p$ and reach a contradiction.  The coaction
$\nu(\kappa\gamma_{p^k})$ contains a term $f \bar\xi_1^{f-1} d \otimes
\bar\xi_1 \sigma\bar\tau_2$, and we must check that there is no class
$x \in H_*(\THH(ju); \F_p)$ with $\kappa x = \bar\xi_1 \sigma\bar\tau_2$.
This is again clear by the B{\"o}kstedt spectral sequence.  The rest of
the proof can safely be omitted.
\end{proof}

\section{ Topological $K$-theory revisited }

We now wish to pass from the homology of the spectra $\THH(ku)$ and
$\THH(ko)$ to their homotopy.  The first case is the analogue for $p=2$
of the discussion in \cite[\S\S5--7]{MS93}, where McClure and Staffeldt
compute the mod~$p$ homotopy groups $\pi_*(\THH(\ell); \Z/p)$ of the
Adams summand $\ell \subset ku_{(p)}$ for odd primes $p$.

The idea is to compute homotopy groups using the Adams spectral
sequence
$$
E_2^{s,t} = \Ext_{A_*}^{s,t}(\F_p, H_*(X; \F_p))
\Longrightarrow \pi_{t-s}(X)^\wedge_p \,,
$$
which is strongly convergent for bounded below spectra $X$ of
$\Z_p$-finite type.  This is difficult for $X = \THH(ku)$ or $\THH(ko)$ at
$p=2$, just as for $\THH(\ell)$ at $p$ odd, but becomes manageable after
introducing suitable finite coefficients, e.g.~after smashing $\THH(ku)$
with the mod~$2$ Moore spectrum $M = C_2$, or smashing $\THH(ko)$ with
the Mahowald spectrum $Y = C_2 \wedge C_\eta$ \cite{Mah82}.  Here $C_f$
denotes the mapping cone of a map $f$.  Neither of these finite CW spectra
are ring spectra, so an extra argument is needed to have products on
$\THH(ku) \wedge M$ and $\THH(ko) \wedge Y$.  We shall manage with the
following very weak version of a ring spectrum.

\begin{defn}
\label{d8.1}
A {\it $\mu$-spectrum} (or an $A_2$ ring spectrum) is a spectrum $R$ with a unit $\eta \: S \to
R$ and multiplication $\mu \: R \wedge R \to R$ that is left and right
unital, but not necessarily associative or commutative.
\end{defn}

The following result is similar to \cite[1.5]{Ok79}, but slightly easier.

\begin{lem}
\label{l8.2}
Let $R$ be a $\mu$-spectrum and let
$$
\xymatrix{
S^k \ar[r]^-{f} & S^0 \ar[r]^-{i} & C_f \ar[r]^-{\pi} & S^{k+1}
}
$$
be a cofiber sequence such that $id_R \wedge id_{C_f} \wedge f \:
R \wedge \Sigma^k C_f \to R \wedge C_f$ is null-homotopic.  Then there
exists a multiplication $\mu \: (R \wedge C_f) \wedge (R \wedge C_f)
\to (R \wedge C_f)$ that makes $R \wedge C_f$ a $\mu$-spectrum and $id_R
\wedge i \: R \to R \wedge C_f$ a map of $\mu$-spectra.
\end{lem}

\begin{proof}
A choice of null-homotopy provides a splitting $m \: R \wedge C_f \wedge C_f
\to R \wedge C_f$ in the cofiber sequence
$$
\xymatrix@C+11pt{
R \wedge \Sigma^k C_f \ar[r]^-{id \wedge id \wedge f} &
R \wedge C_f \ar[r]^-{id \wedge id \wedge i} &
R \wedge C_f \wedge C_f \ar[r]^-{id \wedge id \wedge \pi} &
R \wedge \Sigma^{k+1} C_f
}
$$
which satisfies $m (id \wedge id \wedge i) \simeq id$ (right
unitality).  The difference $id - m(id \wedge i \wedge id) \: R \wedge
C_f \to R \wedge C_f$ restricts trivially over $id \wedge i \: R \to R
\wedge C_f$, so extends over $id \wedge \pi \: R \wedge C_f \to R
\wedge S^{k+1}$ to a map $\phi \:  R \wedge S^{k+1} \to R \wedge C_f$.
There is a short exact sequence
$$
\xymatrix{
[R \wedge \Sigma^{k+1} C_f, R \wedge C_f] \ar[r]^-{i^*} &
[R \wedge S^{k+1}, R \wedge C_f] \ar[r]^-{f^*} &
[R \wedge S^{2k+1}, R \wedge C_f]
}
$$
where the homomorphism $f^* = (id \wedge f \wedge id)^*$ is trivial,
because it can be identified with $(id \wedge id \wedge f)^* \:  [R
\wedge C_f, R] \to [R \wedge \Sigma^k C_f, R]$ under Spanier--Whitehead
duality, and we have assumed that this map $id \wedge id \wedge f$ is
null-homotopic.  (The Spanier--Whitehead dual of $C_f$ is
$\Sigma^{-(k+1)} C_f$.)  Thus $i^* = (id \wedge \Sigma^{k+1} i)^*$ is
surjective and $\phi$ extends further over $id \wedge \Sigma^{k+1} i$
to a map $\psi \:  R \wedge \Sigma^{k+1} C_f \to R \wedge C_f$.  Then
we can subtract $\psi \circ (id \wedge id \wedge \pi)$ from $m$ so as
to make $m(id \wedge i \wedge id) \simeq id$ (left unitality), without
destroying the right unitality.  Once this is achieved, the required
multiplication is defined as the composite
$$
\xymatrix@C+6pt{
(R \wedge C_f) \wedge (R \wedge C_f) \ar[r]^-{id \wedge \tau \wedge id} &
R \wedge R \wedge C_f \wedge C_f \ar[d]^-{\mu \wedge id \wedge id} \\
& R \wedge C_f \wedge C_f \ar[r]^-{m} & (R \wedge C_f) \,,
}
$$
where $\tau$ denotes the twist map.
\end{proof}

In the case $f = 2 \: S^0 \to S^0$, the Moore spectrum $M = C_2$
has cohomology $H^*(M; \F_2) = E(Sq^1)$, which equals $E_1/\!/E(Q_1)$
as an $E_1$-module.  By the Cartan formula $Sq^2$ acts nontrivially on
$H^*(M \wedge M; \F_2)$, and therefore $M$ does not split off from $M
\wedge M$.  So $M$ is not a $\mu$-spectrum, and the map $2 \: M \to M$
is not null-homotopic.  It must therefore factor as the composite
\begin{equation}
\xymatrix{
M \ar[r]^-{\pi} & S^1 \ar[r]^-{\eta} & S^0 \ar[r]^-{i} & M \,.
}
\label{e8.3}
\end{equation}
See e.g.~\cite[1.1]{AT65}.

Similarly, in the case $f = \eta \: S^1 \to S^0$ the mapping cone
$C_\eta$ has cohomology $H^*(C_\eta; \F_2) = E(Sq^2)$.  By the Cartan
formula $Sq^4$ acts nontrivially on $H^*(C_\eta \wedge C_\eta; \F_2)$,
and therefore $C_\eta$ does not split off from $C_\eta \wedge C_\eta$.
In particular, the map $\eta \: \Sigma C_\eta \to C_\eta$ is not
null-homotopic, and must factor as the composite
\begin{equation}
\xymatrix{
\Sigma C_\eta \ar[r]^-{\pi} & S^3 \ar[r]^-{\nu} & S^0 \ar[r]^-{i} & C_\eta \,.
}
\label{e8.4}
\end{equation}
Here $\nu \in \pi_3(S)$ is the Hopf invariant one class.  The Mahowald
spectrum $Y = C_2 \wedge C_\eta$ has cohomology $H^*(Y; \F_2) = E(Sq^1,
Sq^2)$, which equals $A_1/\!/E(Q_1)$ as an $A_1$-module.

\begin{lem}
\label{l8.5}
{\rm (a)}\qua
$\THH(ku) \wedge M$ is a $\mu$-spectrum.

{\rm (b)}\qua
$\THH(ko) \wedge C_\eta$ and $\THH(ko) \wedge Y$ are $\mu$-spectra.
\end{lem}

\begin{proof}
{\rm (a)}\qua
Let $T = \THH(ku)$.  By~(\ref{e8.3}) the map $1 \wedge 2 \: T \wedge M \to T
\wedge M$ factors through $1 \wedge \eta \: \Sigma T \to T$, which in
turn factors as
$$
\xymatrix{
\Sigma T \ar[r]^-{1\wedge\eta} & T \wedge ku \to T \,,
}
$$
since $T = \THH(ku)$ is a $ku$-module spectrum.  But $\eta \in \pi_1(S)$
maps to zero in $\pi_1(ku)$, so this map is null-homotopic, and
Lemma~\ref{l8.2} applies.

{\rm (b)}\qua
Let $T = \THH(ko)$ and $R = \THH(ko) \wedge C_\eta$.  By~(\ref{e8.4}) the map
$1\wedge\eta \: \Sigma R \to R$ factors through $1\wedge\nu \: \Sigma^3
T \to T$, which in turn factors as
$$
\xymatrix{
\Sigma^3 T \ar[r]^-{1\wedge\nu} & T \wedge ko \to T \,,
}
$$
since $T = \THH(ko)$ is a $ko$-module spectrum.  But $\nu \in \pi_3(S)$
maps to zero in $\pi_3(ko)$, so $1\wedge\eta \: \Sigma R \to R$ is
null-homotopic, and Lemma~\ref{l8.2} applies again to prove that $R$
is a $\mu$-spectrum under $T$.

Using once more that $1 \wedge \eta \: \Sigma R \to R$ is null-homotopic
we get that $1 \wedge 2 \: R \wedge M \to R \wedge M$ is null-homotopic
by~(\ref{e8.3}), and so $R \wedge M = \THH(ko) \wedge Y$ is a $\mu$-spectrum
under $R$.
\end{proof}

\begin{lem}
\label{l8.6}
Let $B \subset A$ be a sub Hopf algebra and let $N$ be an $A_*$-comodule
algebra.  Then there is an isomorphism of $A_*$-comodule algebras
$$
(A/\!/B)_* \otimes N \cong A_* \cotensor_{B_*} N \,.
$$
Here $B_*$ is the quotient Hopf algebra of $A_*$ dual to $B$, and
$(A/\!/B)_* = A_* \cotensor_{B_*} \F_p$ is dual to $A \otimes_B \F_p$.
\end{lem}

\begin{proof}
This is analogous to the usual $G$-homeomorphism $G/H \times X \cong G
\times_H X$ for a $G$-space $X$ and subgroup $H \subset G$.  Let $i \:
(A/\!/B)_* \to A_*$ be the inclusion and $\nu \: N \to A_* \otimes N$
the coaction.  The composite homomorphism
$$
\xymatrix{
(A/\!/B)_* \otimes N \ar[r]^-{i \otimes \nu} & A_* \otimes A_* \otimes N
\ar[r]^-{\phi\otimes1} & A_* \otimes N
}
$$
equalizes the two maps to $A_* \otimes B_* \otimes N$, and hence factors
uniquely through $A_* \cotensor_{B_*} N$.  An explicit inverse can be
constructed using the Hopf algebra conjugation $\chi$ on~$A_*$.
\end{proof}

Recall the $n$-th connective Morava $K$-theory spectrum $k(n)$, with
homotopy $\pi_* k(n) = \F_p[v_n]$ where $|v_n| = 2(p^n-1)$, and cohomology
$H^*(k(n); \F_p) \cong A/\!/E(Q_n)$ \cite{BM72}.  Dually, $H_*(k(n);
\F_p) \cong (A/\!/E(Q_n))_* \subset A_*$.  In particular, for $n=1$
and $p=2$ we have $k(1) \simeq ku \wedge M$ with $H_*(k(1); \F_2) \cong
(A/\!/E(Q_1))_* \cong P(\bar\xi_1, \bar\xi_2^2, \bar\xi_k \mid k\ge3)$.
The $n$-th periodic Morava $K$-theory spectrum $K(n)$ is the telescope
$v_n^{-1} k(n)$ of the iterated maps $v_n \: k(n) \to \Sigma^{2(1-p^n)}
k(n)$, with homotopy $\pi_* K(n) = \F_p[v_n, v_n^{-1}]$.  These are
(non-commutative) $S$-algebras for all $p$ and $n$ (\cite{Rob89}, \cite{La03}).

\begin{prop}
\label{p8.7}
{\rm (a)}\qua
There are $A_*$-comodule algebra isomorphisms
\begin{align*}
H_*(\THH(ku) \wedge M; \F_2) &\cong (A/\!/E(Q_1))_* \otimes
E(\sigma\bar\xi_1^2, \sigma\bar\xi_2^2) \otimes P(\sigma\bar\xi_3) \\
& \cong H_*(k(1); \F_2) \otimes E(\lambda_1, \lambda_2) \otimes P(\mu) \,.
\end{align*}
Here $\nu(\sigma\bar\xi_1^2) = 1 \otimes \sigma\bar\xi_1^2$,
$\nu(\sigma\bar\xi_2^2) = 1 \otimes \sigma\bar\xi_2^2$, and
$\nu(\sigma\bar\xi_3) = 1 \otimes \sigma\bar\xi_3 + \bar\xi_1 \otimes
\sigma\bar\xi_2^2$.

The classes $\lambda_1 = \sigma\bar\xi_1^2$, $\lambda_2 =
\sigma\bar\xi_2^2$ and $\mu = \sigma\bar\xi_3 + \bar\xi_1 \cdot
\sigma\bar\xi_2^2$ (in degrees $3$, $7$ and $8$, respectively) are
$A_*$-comodule primitives.

{\rm (b)}\qua
There are $A_*$-comodule algebra isomorphisms
\begin{align*}
H_*(\THH(ko) \wedge Y; \F_2) &\cong (A/\!/E(Q_1))_* \otimes
E(\sigma\bar\xi_1^4, \sigma\bar\xi_2^2) \otimes P(\sigma\bar\xi_3) \\
& \cong H_*(k(1); \F_2) \otimes E(\lambda_1, \lambda_2) \otimes P(\mu) \,.
\end{align*}
Here $\nu(\sigma\bar\xi_1^4) = 1 \otimes \sigma\bar\xi_1^4$,
$\nu(\sigma\bar\xi_2^2) = 1 \otimes \sigma\bar\xi_2^2 + \bar\xi_1^2
\otimes \sigma\bar\xi_1^4$ and $\nu(\sigma\bar\xi_3) = 1 \otimes
\sigma\bar\xi_3 + \bar\xi_1 \otimes \sigma\bar\xi_2^2 + \bar\xi_2
\otimes \sigma\bar\xi_1^4$.

The exterior classes $\lambda_1 = \sigma\bar\xi_1^4$ and $\lambda_2
= \sigma\bar\xi_2^2 + \bar\xi_1^2 \cdot \sigma\bar\xi_1^4$ are
$A_*$-comodule primitives, while $\mu = \sigma\bar\xi_3 + \bar\xi_1
\cdot \sigma\bar\xi_2^2$ has
$$
\nu(\mu) = 1 \otimes \mu + \bar\xi_1^2 \otimes \bar\xi_1 \cdot \lambda_1 +
(\bar\xi_2 + \bar\xi_1^3) \otimes \lambda_1 \,.
$$
The squared class $\mu^2 = (\sigma\bar\xi_3)^2$ is $A_*$-comodule
primitive.
\end{prop}

In each case, the homology algebra on the left hand side has the unit
and product induced by the $\mu$-spectrum structure from
Lemma~\ref{l8.5}.  This product is in fact associative and graded
commutative, in view of the exhibited additive and multiplicative
isomorphism with the associative and graded commutative algebra on the
right hand side.

\begin{proof}
{\rm (a)}\qua
By Corollary~\ref{c5.14}(a) there is an $A_*$-comodule algebra isomorphism
$$
H_*(\THH(ku) \wedge M; \F_2) \cong (A/\!/E_1)_* \otimes (E_1/\!/E(Q_1))_*
\otimes E(\sigma\bar\xi_1^2, \sigma\bar\xi_2^2) \otimes P(\sigma\bar\xi_3)
$$
with the diagonal $A_*$-comodule structure on the first two tensor
factors, and the claimed coaction on the remaining generators.

Since $H_*(M; \F_2) \cong (E_1/\!/E(Q_1))_*$ is in fact an
$A_*$-comodule algebra, there is an $A_*$-comodule algebra isomorphism
$$
(A/\!/E_1)_* \otimes (E_1/\!/E(Q_1))_* \cong A_* \cotensor_{E_{1*}}
(E_1/\!/E(Q_1))_* \cong (A/\!/E(Q_1))_*
$$
by Lemma~\ref{l8.6}.  We are free to replace the polynomial generator
$\sigma\bar\xi_3$ by the primitive class $\mu$, since their difference
$\bar\xi_1 \cdot \sigma\bar\xi_2^2$ has square zero.

{\rm (b)}\qua
By Theorem~\ref{t6.2}(a) there is an $A_*$-comodule algebra isomorphism
$$
H_*(\THH(ko) \wedge Y; \F_2) \cong (A/\!/A_1)_* \otimes (A_1/\!/E(Q_1))_*
\otimes E(\sigma\bar\xi_1^4, \sigma\bar\xi_2^2) \otimes P(\sigma\bar\xi_3)
\,,
$$
with the claimed coaction on the exterior and polynomial generators.
Again, by Lemma~\ref{l8.6} there is an isomorphism of $A_*$-comodule algebras
$$
(A/\!/A_1)_* \otimes (A_1/\!/E(Q_1))_*
\cong A_* \cotensor_{A_{1*}} (A_1/\!/E(Q_1))_* \cong (A/\!/E(Q_1))_* \,.
$$
The classes $\lambda_1$, $\lambda_2$ and $\mu$ are defined as in the
statement of the proposition, and their coactions are obtained by direct
calculation.
\end{proof}

\begin{lem}
\label{l8.8}
{\rm (a)}\qua
The $E_2$-term of the Adams spectral sequence converging to
$\pi_*(\THH(ku) \wedge M)$ is
$$
E_2^{**} \cong P(v_1) \otimes E(\lambda_1, \lambda_2) \otimes P(\mu)
$$
with $\lambda_1$, $\lambda_2$, $\mu$ and $v_1 = [\xi_2]$ in
bidegrees $(0, 5)$, $(0, 7)$, $(0, 8)$ and $(1, 3)$, respectively.

{\rm (b)}\qua
The $E_2$-term of the Adams spectral sequence for
$\pi_*(\THH(ko) \wedge Y)$ is
$$
E_2^{**} \cong \bigl( P(v_1) \otimes E(\lambda_2, \lambda_3)
\oplus E(\lambda_2) \{\lambda_1\} \bigr) \otimes P(\mu^2)
$$
with $\lambda_1$, $\lambda_2$, $\lambda_3 = \lambda_1\mu$, $\mu^2$ and
$v_1 = [\xi_2]$ in bidegrees $(0, 5)$, $(0, 7)$, $(0, 13)$, $(0, 16)$
and $(1, 3)$, respectively.  It is the homology of the algebra
$$
P(v_1) \otimes E(\lambda_1, \lambda_2) \otimes P(\mu)
$$
with respect to the differential $d(\mu) = v_1 \lambda_1$, with cycles
$\lambda_1$, $\lambda_2$ and $v_1$.
\end{lem}

\begin{proof}
{\rm (a)}\qua
By change-of-rings, the $E_2$-term of the Adams spectral sequence
for $X = \THH(ku) \wedge M$ is
\begin{align*}
E_2^{**} &= \Ext_{A_*}^{**}(\F_2, (A/\!/E(Q_1))_* \otimes E(\lambda_1,
\lambda_2) \otimes P(\mu)) \\
&\cong \Ext_{E(Q_1)_*}^{**}(\F_2, E(\lambda_1, \lambda_2) \otimes P(\mu)) \\
&\cong P(v_1) \otimes E(\lambda_1, \lambda_2) \otimes P(\mu)
\end{align*}
as a graded algebra.  Here $E(Q_1)_* = E(\xi_2)$ and
$\Ext_{E(\xi_2)}^{**}(\F_2, \F_2) = P(v_1)$ with $v_1 = [\xi_2]$
in the cobar complex.

{\rm (b)}\qua
The $E_2$-term of the Adams spectral sequence for $X = \THH(ko) \wedge Y$
is
\begin{align*}
E_2^{**} &= \Ext_{A_*}^{**}(\F_2, (A/\!/E(Q_1))_*
\otimes E(\lambda_1, \lambda_2) \otimes P(\mu)) \\
&\cong \Ext_{E(Q_1)_*}^{**}(\F_2, E(\lambda_1, \lambda_2) \otimes P(\mu)) \\
&\cong \Ext_{E(Q_1)_*}^{**}(\F_2, \F_2\{1, \lambda_1, \lambda_2, \mu,
	\lambda_1\lambda_2, \lambda_1\mu, \lambda_2\mu,
	\lambda_1\lambda_2\mu\} \otimes P(\mu^2)) \,.
\end{align*}
Here
\begin{multline*}
\F_2\{1, \lambda_1, \lambda_2, \mu, \lambda_1\lambda_2, \lambda_1\mu,
	\lambda_2\mu, \lambda_1\lambda_2\mu\} \cong \\
\F_2\{1, \lambda_2, \lambda_1\mu, \lambda_1\lambda_2\mu\} \oplus
	E(Q_1)_*\{\lambda_1, \lambda_1\lambda_2\}
\end{multline*}
as $E(Q_1)_*$-comodules.   For $\lambda_1$, $\lambda_2$ and $\mu^2$ are
$A_*$-comodule primitives, while $\nu(\mu)$ maps to $1 \otimes \mu + \xi_2
\otimes \lambda_1$ in $E(Q_1)_* \otimes H_*(\THH(ko) \wedge Y; \F_2)$.  So
$$
E_2^{**} = \bigl( P(v_1) \otimes E(\lambda_2, \lambda_1\mu) \oplus
E(\lambda_2)\{\lambda_1\} \bigr) \otimes P(\mu^2) \,.
$$
We let $\lambda_3 = \lambda_1\mu$ to obtain the claimed formula.
\end{proof}

To determine the differentials in the Adams spectral sequence for
$\pi_*(\THH(ku) \wedge M)$ or $\pi_*(\THH(ko) \wedge Y)$ we first compute
the $v_1$-periodic homotopy.  This is in turn easy to derive from the
$K(1)$-homology.  Recall that $BP_* BP \cong BP_*[t_k \mid k \ge 1]$
with $|t_k| = 2(2^k-1)$ and $K(1)_* = \F_2[v_1, v_1^{-1}]$ with
$|v_1| = 2$.

\begin{lem}
\label{l8.9}
For $p=2$ there are isomorphisms of $K(1)_*$-algebras
$$
K(1)_*(ku) \cong K(1)_*[t_k\mid k\ge1]/(v_1 t_k^2 = v_1^{2^k} t_k)
\cong K(1)_*[u_k \mid k\ge1]/(u_k^2=u_k)
$$
and
$$
K(1)_*(ko) \cong K(1)_*[t_k \mid k\ge2]/(v_1 t_k^2 = v_1^{2^k} t_k)
\cong K(1)_*[u_k \mid k\ge2]/(u_k^2=u_k) \,.
$$
\end{lem}

\begin{proof}
We first follow the proof of \cite[5.3]{MS93}.  We have $K(1)_*(BP)
\cong K(1)_* \otimes_{BP_*} BP_*(BP) \cong K(1)_* [t_k \mid k\ge1]$,
where $|t_k| = 2(2^k-1)$.  The spectrum $ku_{(2)} = BP\langle1\rangle$
can be constructed from $BP$ by Baas--Sullivan cofiber sequences killing
the classes $v_{k+1}$ for $k\ge1$, which map to
$$
\eta_R(v_{k+1}) \equiv v_1 t_k^2 + v_1^{2^k} t_k \mod (\eta_R(v_2),
\dots, \eta_R(v_k))
$$
by \cite[6.1.13]{Ra04}.  So each $\eta_R(v_{k+1}) \in K(1)_*(BP)$ is not
a zero divisor mod $(\eta_R(v_2), \dots, \eta_R(v_k))$, and
$$
K(1)_*(ku) \cong K(1)_*[t_k \mid k\ge 1]/(v_1 t_k^2 = v_1^{2^k} t_k)
\,.
$$
Substituting $u_k = v_1^{1-2^k} t_k$, the relations become $u_k^2 = u_k$
for each $k\ge1$.

The cofiber sequence
$$
\xymatrix{
\Sigma ko \ar[r]^-{\eta} & ko \ar[r]^-{c} & ku \ar[r]^-{\partial} & \Sigma^2 ko
}
$$
induces a short exact sequence
$$
\xymatrix{
0 \to K(1)_*(ko) \ar[r]^-{c_*} & K(1)_*(ku) \ar[r]^-{\partial_*} & \Sigma^2 K(1)_*(ko) \to 0
}
$$
since multiplication by $\eta$ is zero in $K(1)$-homology.  The connecting
map $\partial$ right multiplies by $Sq^2$ in mod~$2$ cohomology, so
right ``comultiplies'' with the dual class $\bar\xi_1^2$ in mod~$2$
homology, which corresponds to $t_1$ in $BP_* BP$ \cite[p.~488]{Za72}.
From the coproduct formula for $\Delta(t_j)$ \cite[A2.1.27(e)]{Ra04}
it follows that $\partial_*(t_1) = \Sigma^2(1)$ while $\partial_*(t_k)
= 0$ for $k\ge2$.  Hence we can identify $K(1)_*(ko)$ with the claimed
subalgebra of $K(1)_*(ku)$, via the $K(1)_*$-algebra homomorphism $c_*$.
\end{proof}

\begin{lem}
\label{l8.10}
The unit maps $R \to \THH(R)$ induce isomorphisms
$$
\xymatrix{
K(1)_*(ku) \ar[r]^-{\cong} & K(1)_* \THH(ku)
}
$$
and
$$
\xymatrix{
K(1)_*(ko) \ar[r]^-{\cong} & K(1)_* \THH(ko) \,.
}
$$
\end{lem}

\begin{proof}
The proof of \cite[5.3]{MS93} continues as follows.  There is a
$K(n)$-based B{\"o}kstedt spectral sequence
$$
E^2_{**} = HH^{K(n)_*}_*(K(n)_*(R))
\Longrightarrow K(n)_* \THH(R)
$$
for every $S$-algebra $R$, derived like the one in~(\ref{e4.1}), but by
applying $K(n)$-homology to the skeleton filtration of $\THH(R)$.  The
identification of the $E^2$-term uses the K{\"u}nneth formula for
Morava $K$-theory.  When $K(n)_*(R)$ is concentrated in degrees $*
\equiv 0 \mod |v_n|$, then $K(n)_*(R) \cong K(n)_* \otimes_{\F_p}
K(n)_0(R)$ and we can rewrite the $E^2$-term as
$$
E^2_{**} \cong K(n)_* \otimes_{\F_p} HH^{\F_p}_*(K(n)_0(R)) \,.
$$
This is the case for $n=1$, $p=2$ and $R = ku$, when $K(1)_0(ku)$
is the colimit over~$m$ of the algebras $\F_2[u_k \mid 1 \le k \le
m]/(u_k^2=u_k) \cong \prod_{i=1}^{2^m} \F_2$.  Then for each $m$ the
unit map $\prod_{i=1}^{2^m} \F_2 \to HH^{\F_2}_*(\prod_{i=1}^{2^m} \F_2)$
is an isomorphism, so by passage to the colimit the unit map
$$
K(1)_0(ku) \to HH^{\F_2}_*(K(1)_0(ku))
$$
is an isomorphism.  Hence the $K(1)$-based B{\"o}kstedt spectral sequence
collapses at the edge $s=0$, and the unit map $ku \to \THH(ku)$ induces
the asserted isomorphism.

Likewise, $K(1)_0(ko)$ is the colimit of the algebras $\F_2[u_k \mid
2 \le k \le m]/(u_k^2 = u_k) \cong \prod_{i=1}^{2^{m-1}} \F_2$, and
the same argument shows that the unit map $ko \to \THH(ko)$ induces an
isomorphism in $K(1)$-homology.
\end{proof}

The mod~$2$ Moore spectrum $M$ admits a degree~$8$ self-map $v_1^4 \:
M \to \Sigma^{-8} M$ that induces multiplication by $v_1^4$ in $ku_*(M)
= k(1)_* = P(v_1)$.  Smashing with $C_\eta$ yields a self-map $v_1^4 \: Y
\to \Sigma^{-8} Y$, which admits a fourth root $v_1 \: Y \to \Sigma^{-2}
Y$ (up to nilpotent maps).  It induces multiplication by $v_1$ in
$ko_*(Y) = k(1)_* = P(v_1)$.  See \cite[1.2]{DM82}.

\begin{lem}
\label{l8.11}
Let $X$ be a spectrum such that $K(1)_*(X) = 0$.  Then $v_1^{-1} \pi_*(X
\wedge Y) = 0$ and $v_1^{-4} \pi_*(X \wedge M) = 0$.
\end{lem}

\begin{proof}
Recall that $j$ is the homotopy fiber of the map $\psi^3-1 \: ko_{(2)}
\to bspin_{(2)}$.  The unit map $e \: S \to j$ induces an equivalence
of mapping telescopes
$$
\xymatrix{
v_1^{-1} (S \wedge Y) \ar[r]^-{\simeq} & v_1^{-1}(j \wedge Y) \,.
}
$$
See e.g.~the case $n=0$ of \cite[1.4]{Mah82}.  Here $v_1^{-1}(ko_{(2)}
\wedge Y) \simeq v_1^{-1} k(1) = K(1)$ and likewise $v_1^{-1}(bspin_{(2)}
\wedge Y) \simeq K(1)$, so there is a cofiber sequence of spectra
$$
\xymatrix{
v_1^{-1} Y \to K(1) \ar[r]^-{\psi} & K(1) \,.
}
$$
Furthermore, $v_1^{-1} Y \simeq v_1^{-4} Y$ sits in a cofiber sequence
$$
\xymatrix{
v_1^{-4} \Sigma M \ar[r]^-{\eta} & v_1^{-4} M \to v_1^{-4} Y
}
$$
where $\eta$ is nilpotent ($\eta^4=0$).

From the first cofiber sequence it follows that if $K(1)_*(X) = 0$,
then $\pi_*(X \wedge v_1^{-1}Y) = v_1^{-1} \pi_*(X \wedge Y) = 0$.
From the second cofiber sequence it then follows that multiplication
by $\eta$ is an isomorphism on $\pi_*(X \wedge v_1^{-4}M) =
v_1^{-4} \pi_*(X \wedge M)$.  Since $\eta$ is nilpotent this implies
that $v_1^{-4} \pi_*(X \wedge M) = 0$.
\end{proof}
  
\begin{cor}
\label{c8.12}
The unit maps induce isomorphisms
$$
\xymatrix{
K(1)_* = v_1^{-4} \pi_*(ku \wedge M) \ar[r]^-{\cong} & v_1^{-4}
\pi_*(\THH(ku) \wedge M)
}
$$
and
$$
\xymatrix{
K(1)_* = v_1^{-1} \pi_*(ko \wedge Y) \ar[r]^-{\cong} & v_1^{-1} \pi_*(\THH(ko)
\wedge Y) \,.
}
$$
\end{cor}

\begin{proof}
Apply Lemmas~\ref{l8.10} and~\ref{l8.11} to the cofiber of the unit map
$R \to \THH(R)$, for $R = ku$ and $R = ko$, respectively.
\end{proof}

\begin{thm}
\label{t8.13}
Consider the Adams spectral sequence
$$
E_2^{**} = P(v_1) \otimes E(\lambda_1, \lambda_2) \otimes P(\mu)
$$
for $\pi_*(\THH(ku) \wedge M)$, with $\lambda_1 = \sigma\bar\xi_1^2$,
$\lambda_2 = \sigma\bar\xi_2^2$, $\mu = \sigma\bar\xi_3 + \bar\xi_1 \cdot
\sigma\bar\xi_2^2$ and $v_1 = [\xi_2]$ in bidegrees $(s, t) = (0, 3)$,
$(0, 7)$, $(0, 8)$ and $(1, 3)$, respectively.  Recursively define
$$
\lambda_n = \lambda_{n-2} \mu^{2^{n-3}}
$$
for $n\ge3$.  Likewise define $r(1) = 2$, $r(2) = 4$, $r(n) = 2^n +
r(n-2)$, $s(1) = 3$, $s(2) = 7$ and $s(n) = 2^n + s(n-2)$, for $n\ge3$.
So $\lambda_n$ has bidegree $(0, s(n))$ and $2 r(n) + s(n) = 2^{n+2} - 1$.

Then the Adams spectral sequence has differentials
generated by
$$
d^{r(n)}(\mu^{2^{n-1}}) = v_1^{r(n)} \lambda_n
$$
for all $n\ge1$.  This leaves the $E_\infty$-term
$$
E_\infty^{**} = P(v_1)\{1\} \oplus
\bigoplus_{n=1}^\infty P_{r(n)}(v_1) \{\lambda_n\} \otimes
E(\lambda_{n+1}) \otimes P(\mu^{2^n}) \,.
$$

Hence $\pi_*(\THH(ku) \wedge M) = \pi_*(\THH(ku); \Z/2)$ is generated
as a $P(v_1)$-module by elements $1$, $x_{n,m} = \lambda_n \mu^{2^n
m}$ and $x'_{n,m} = \lambda_n \lambda_{n+1} \mu^{2^n m}$ for $n\ge1$
and $m\ge0$.  Here $|x_{n,m}| = s(n) + 2^{n+3} m$ and $|x'_{n,m}| =
s(n) + s(n+1) + 2^{n+3} m$.  The module structure is generated by the
relations $v_1^{r(n)} x_{n,m} = 0$ and $v_1^{r(n)} x'_{n,m} = 0$ for
$n\ge1$, $m\ge0$.
\end{thm}

\begin{proof}
The classes $\lambda_1$, $\lambda_2$ and $\mu$ were introduced in
Proposition~\ref{p8.7}, and the Adams spectral sequence $E_2$-term was
found in Lemma~\ref{l8.8}. By Corollary~\ref{c8.12} the abutment
$\pi_*(\THH(ku) \wedge M)$ is all $v_1$-torsion, except the direct
summand $\pi_*(ku \wedge M) = P(v_1)$, which is included by the unit
map.  Hence every class $\lambda_n$ is $v_1$-torsion, so there is some
integer $r(n)$ such that $v_1^{r(n)} \lambda_n$ is hit by a
differential.

Suppose by induction that the $d^{r(k)}$-differentials for $1
\le k < n$ have been found, leaving the term
\begin{align*}
E_{r(n-1)+1}^{**} &= P(v_1) \otimes E(\lambda_n, \lambda_{n+1})
\otimes P(\mu^{2^{n-1}}) \\
&\qquad \oplus
\bigoplus_{k=1}^{n-1} P_{r(k)}(v_1)\{\lambda_k\} \otimes E(\lambda_{k+1})
\otimes P(\mu^{2^k}) \,.
\end{align*}
Consider the $P(v_1)$-module generated by $\lambda_n$.  Let $r$ be
minimal such that $v_1^r \lambda_n$ is a boundary.  The source $x$ of
such a differential cannot be divisible by $v_1$, since $r$ is minimal,
so $d^r(x) = v_1^r \lambda_n$ where $x$ has Adams filtration $s=0$ and even
total degree.  Furthermore, $v_1^r \lambda_n$ is not $v_1$-torsion at this
term, so $x$ cannot be $v_1$-torsion.  Likewise, $\lambda_{n+1}$ does
not annihilate $v_1^r \lambda_n$, so $\lambda_{n+1}$ cannot annihilate
$x$ either.  This forces $x \in P(\mu^{2^{n-1}})$.  By Lemma~\ref{l8.5},
$\THH(ku) \wedge M$ is a $\mu$-spectrum, so $d^r$ is a derivation and
the Leibniz rule shows that $x = \mu^{2^{n-1}}$, since $d^r$ on any
higher power of $\mu^{2^{n-1}}$ must be divisible by $\mu^{2^{n-1}}$.
Hence $d^r(\mu^{2^{n-1}}) = v_1^r \lambda_n$, and by a degree count we
must have $r = r(n)$.

To complete the induction step, we must compute the $E_{r(n)+1}$-term
of the Adams spectral sequence.  The $d^{r(n)}$-differential does not
affect the summands $P_{r(k)}(v_1)\{\lambda_k\} \otimes
E(\lambda_{k+1}) \otimes P(\mu^{2^k})$ for $1 \le k < n$, and is zero
on $E(\lambda_{n+1}) \otimes P(\mu^{2^n})$.  It acts on $P(v_1) \{1,
\lambda_n, \mu^{2^{n-1}}, \lambda_n \mu^{2^{n-1}}\}$, leaving
$P_{r(n)}(v_1) \{\lambda_n\} \oplus P(v_1) \otimes E(\lambda_{n+2})$,
where by definition $\lambda_{n+2} = \lambda_n \mu^{2^{n-1}}$.  This
shows that the term $P(v_1) \otimes E(\lambda_n, \lambda_{n+1}) \otimes
P(\mu^{2^{n-1}})$ at the $E_{r(n-1)+1}$-term gets replaced by the
direct sum of $P(v_1) \otimes E(\lambda_{n+1}, \lambda_{n+2}) \otimes
P(\mu^{2^n})$ and $P_{r(n)}(v_1) \{\lambda_n\} \otimes E(\lambda_{n+1})
\otimes P(\mu^{2^n})$.
\end{proof}

\begin{thm}
\label{t8.14}
Consider the Adams spectral sequence
$$
E_2^{**} = \bigl( P(v_1) \otimes E(\lambda_2, \lambda_3)
\oplus E(\lambda_2) \{\lambda_1\} \bigr) \otimes P(\mu^2)
$$
for $\pi_*(\THH(ko) \wedge Y)$, with $\lambda_1 = \sigma\bar\xi_1^4$,
$\lambda_2 = \sigma\bar\xi_2^2 + \bar\xi_1^2 \cdot \sigma\bar\xi_1^4$,
$\lambda_3 = \sigma\bar\xi_1^4 (\sigma\bar\xi_3 + \bar\xi_1 \cdot
\sigma\bar\xi_2^2)$,  $\mu^2 = (\sigma\bar\xi_3)^2$ and $v_1 = [\xi_2]$
in bidegrees $(s, t) = (0, 5)$, $(0, 7)$, $(0, 13)$, $(0, 16)$ and $(1,
3)$, respectively.  Recursively define
$$
\lambda_n = \lambda_{n-2} \mu^{2^{n-3}}
$$
for $n\ge4$.  Likewise define $r(1) = 1$, $r(2) = 4$, $r(n) = 2^n +
r(n-2)$, $s(1) = 5$, $s(2) = 7$ and $s(n) = 2^n + s(n-2)$, for $n\ge3$.
So $\lambda_n$ has bidegree $(0, s(n))$ and $2 r(n) + s(n) = 2^{n+2} - 1$.

Then the Adams spectral sequence has differentials
generated by
$$
d^{r(n)}(\mu^{2^{n-1}}) = v_1^{r(n)} \lambda_n
$$
for all $n\ge2$.  This leaves the $E_\infty$-term
$$
E_\infty^{**} = P(v_1)\{1\} \oplus
\bigoplus_{n=1}^\infty P_{r(n)}(v_1) \{\lambda_n\} \otimes
E(\lambda_{n+1}) \otimes P(\mu^{2^n}) \,.
$$

Hence $\pi_*(\THH(ko) \wedge Y) = \pi_*(\THH(ko); Y)$ is generated as
a $P(v_1)$-module by elements $1$, $x_{n,m} = \lambda_n \mu^{2^n m}$
and $x'_{n,m} = \lambda_n \lambda_{n+1} \mu^{2^n m}$ for $n\ge1$
and $m\ge0$.  Here $|x_{n,m}| = s(n) + 2^{n+3} m$ and $|x'_{n,m}| =
s(n) + s(n+1) + 2^{n+3} m$.  The module structure is generated by the
relations $v_1^{r(n)} x_{n,m} = 0$ and $v_1^{r(n)} x'_{n,m} = 0$ for
$n\ge1$, $m\ge0$.
\end{thm}

\begin{proof}
Starting with an imagined $E_1$-term
$$
E_1^{**} = P(v_1) \otimes E(\lambda_1, \lambda_2) \otimes P(\mu)
$$
and differential $d_1(\mu) = v_1 \lambda_1$, the proof is the same as
for Theorem~\ref{t8.13}.
\end{proof}

\Addresses\recd

\end{document}